\documentclass[a4paper,12pt]{scrreprt}

\usepackage[utf8x]{inputenc}

\usepackage{aky-layout}
\usepackage{aky-math}

\usepackage{makeidx}
\makeindex

\newcommand{\Dbcoh}{\ensuremath{\mathbf{D}^{\mathrm{b}}_{\mathrm{coh}}}\xspace}

\newcommand{\PfCorr}{\mathbf{PfCorr}} % perfect correspondences
\newcommand{\Pf}{\mathbf{Pf}} % derived cat of perfect complexes
\newcommand{\PfDG}{\mathbf{Pf}^\bull} % DG cat of perfect complexes
\newcommand{\piint}{\pi_{\mathrm{int}}}
\newcommand{\picof}{\pi_{\mathrm{cof}}}
\newcommand{\pifib}{\pi_{\mathrm{fib}}}
\newcommand{\HoDGQe}{\mathbf{HoDGQe}} % homotopy category wrt quasi-equivalences
\newcommand{\HoDGMo}{\mathbf{HoDGMo}} % homotopy category wrt Morita equivs
\newcommand{\SmPrDGQe}{\mathbf{SmPrDGQe}}
\newcommand{\SmPrDGMo}{\mathbf{SmPrDGMo}}
\newcommand{\GeoDGQe}{\mathbf{GeoDGQe}}
\newcommand{\GeoDGMo}{\mathbf{GeoDGMo}}
\newcommand{\QsiHom}{\mathbf{QsiHom}} % category of quasi-functors
\newcommand{\AddInv}{\mathbf{AddInv}} % additive invariants
\newcommand{\GrChMot}{\mathbf{GrChMot}}

\newcommand{\A}{\sA}

\renewcommand{\L}{\bL}
\newcommand{\M}{\sM}
\renewcommand{\H}{\bH}

\renewcommand{\AA}{{\sA^\bull}}
\newcommand{\BB}{{\sB^\bull}}
\newcommand{\CC}{{\sC^\bull}}
\newcommand{\Cdg}{{\Comp^\bull}}
\newcommand{\Ddg}{{\D^\bull}}
\newcommand{\RHomDG}{{\mathbf{RHom}^\bull}} % internal hom of HoDGQe
\newcommand{\ModDG}{{\mathbf{Mod}^\bull}} % DG category of DG modules
\newcommand{\IntDG}{{\mathbf{Int}^\bull}} % interior DG category
\newcommand{\QsiHomDG}{{\mathbf{QsiHom}^\bull}} % DG cat of quasi-functors

\newcommand{\Flag}{\ensuremath{\bF}\xspace} % flag variety
\newcommand{\inpr}{\cdot} % intersection product
\newcommand{\IdM}{\bI} % identity motive
\newcommand{\LefM}{\bL} % Lefschetz motive
\newcommand{\TateM}{\bT} % Tate motive
\newcommand{\FF}{\mathbf{FF}} % fully faithful functors
\makeatletter
\def\thickhrulefill{\leavevmode \leaders \hrule height 1pt\hfill \kern \z@}
\def\maketitle{%
  \null
  \thispagestyle{empty}%
  \vskip 1cm
  \vfil
  \begin{center}
    \Large \strut \@title \par
  \end{center}
  \par
  \vfil
  \begin{center}
    \emph{by } \normalfont\@author\par
  \end{center}
  \vfil
  \vfil
  \vfil
  \null
  }
\makeatother
\author{A. Kh. Yusufzai\\Freie Universit\"at Berlin\\\texttt{kadeel@gmail.com}}
\title{\uppercase{Perfect correspondences and\\Chow motives}}
\date{15 September 2013}

\begin{document}

\maketitle

\begin{center}
\textbf{Abstract}
\end{center}

It is an open conjecture of Orlov that the bounded derived category of coherent sheaves of a smooth projective variety determines its Chow motive with rational coefficients.
In this master's thesis we introduce a category of \emph{perfect correspondences}, whose objects are smooth projective varieties and morphisms $X \to Y$ are perfect complexes on $X \times Y$.
We show that isomorphism in this category is the same as equivalence of derived categories, and use this to show that the derived category determines the noncommutative Chow motive (in the sense of Tabuada) and, up to Tate twists, the commutative Chow motive with rational coefficients.
In particular, all additive invariants like K-theory and Hochschild or cyclic homology depend only on the derived category.

\tableofcontents

\chapter*{Introduction}
\label{chapter:introduction}

\section*{Overview}
\label{sec:overview}

\para
\emph{Derived categories.}
The derived category is a framework developed by Verdier as a general setting for homological algebra.
For any abelian category $\A$ there is an associated derived category $\D(\A)$ and a bounded version $\Db(\A)$.
Given a scheme $X$, the sheaves of $\O_X$-modules on $X$ form an abelian category $\Mod(\O_X)$.
When $X$ is noetherian, the coherent sheaves, which correspond locally to \emph{finitely generated} modules, form a full abelian subcategory $\Coh(\O_X)$.
Let us call the category $\Db(\Coh(\O_X))$ the \dfn[]{derived category of $X$}.
When $X$ is smooth, this is the same as the full subcategory $\Pf(X) \sub \D(\Mod(\O_X))$ of perfect complexes.

Let $\Var_K$ denote the category of smooth projective varieties over a field $K$.
We say that $X, Y \in \Var_K$ are \dfn[]{derived equivalent} if their derived categories are equivalent.
The derived category seems to lose just enough information so that derived equivalence becomes a geometrically significant invariant.
In fact, in the known examples of derived equivalence, the varieties always seem to have some deep geometric relationship.
For example, an abelian variety $X$ is derived equivalent to its dual $\hat{X}$ (\cit{mukai1981duality}), and varieties connected by a \emph{flop transformation} are derived equivalent (\cit{bridgeland2002flops}).
On the other hand, two K3 surfaces are derived equivalent if and only if their Mukai lattices are Hodge isometric (\cit{mukai1987moduli}), and further when the canonical sheaves are ample or anti-ample, derived equivalence is as strong as isomorphism (\cit{bondal2001reconstruction}).

Therefore, it is important for the classification of algebraic varieties to understand exactly which information is retained by the derived category.
In particular, there are various functors defined on $\Var_K$ like singular cohomology (for $K = \C$), $\ell$-adic cohomology (for $\chr(K) \ne \ell$), K-theory, Hochschild and cyclic homologies.
A natural question is, \emph{which of these functors are invariants of the derived category?}

\para
\emph{Cohomology of varieties.}
Let $\Var_K$ denote the category of smooth projective schemes over a field $K$.
A \dfn[]{Weil cohomology theory} (with coefficients in $\Q$) is a functor to the category of nonnegatively graded $\Q$-vector spaces, satisfying certain axioms (e.g. Poincar\'e duality, K\"unneth isomorphism).
There are several classical Weil cohomology theories for smooth projective varieties, all capturing different informations.
Though there is no universal cohomology theory, Grothendieck defined a universal functor, valued in a category of \emph{Chow motives}, through which every cohomology theory factors.
  \[ \begin{tikzcd}
    \Var_K^\circ \arrow{rr}{H^*}\arrow{rd}{M_\Q}
      &
      & \GrVec^+_\Q \\
      \
      & \CHM_K(\Q) \arrow[dotted]{ru}
      &
  \end{tikzcd} \]

When $K = \C$, it is known that that the derived category determines the singular cohomology with coefficients in $\Q$ (\cit{orlov2003derived}, 2.1.12).
Hence it is natural to ask whether the same holds for other cohomology theories.
This is the content of a conjecture made by Orlov in 2005.

\begin{paraconj}[\cit{orlov2005derived}]
Let $X$ and $Y$ be smooth projective varieties over a field $K$.
If their derived categories $\Pf(X)$ and $\Pf(Y)$ are equivalent, then their Chow motives with rational coefficients are isomorphic.
\end{paraconj}

\para
\emph{DG categories.}
Though we are interested in equivalence of derived categories, it has been widely recognized since the 80's (\cit{gelfand1988methods}) that the derived category is poorly behaved in many ways.
For example, the non-functoriality of cones of morphisms is one source of difficulty.
In general, the construction $\A \mapsto \D(\A)$ destroys most good properties $\A$ may have, like the existence of limits and colimits.
It also kills the local nature of sheaves on a scheme: for example, the derived category of the projective line is not determined by the derived categories of its affine parts.

One approach to the problem, going back to \cit{bondal1990enhanced}, is the theory of DG categories.
A \dfn[]{DG category} (over $K$) is a category where the sets of morphisms between objects form complexes of $K$-modules.
Every DG category $\AA$ can be flattened into an ordinary category $\H^0(\AA)$ called the \dfn[]{homotopy category}, where the morphisms are given by zeroth cohomologies of complexes of morphisms in $\AA$.
For a $K$-scheme $X$ there exists a DG category $\PfDG(X)$ whose homotopy category is the derived category $\Pf(X)$.
The crucial point is the existence of a homotopy theory of DG categories, and the homotopy category $\HoDGQe_K$ where an isomorphism $\PfDG(X) \isoto \PfDG(Y)$ is the same thing as a equivalence of derived categories $\Pf(X) \isoto \Pf(Y)$.
In particular, Orlov's conjecture can be reformulated as

\begin{paraconj}
Let $X$ and $Y$ be smooth projective varieties over a field $K$.
If their DG categories $\PfDG(X)$ and $\PfDG(Y)$ are homotopy equivalent, then their Chow motives with rational coefficients are isomorphic.
\end{paraconj}

\para
\emph{Noncommutative geometry.}
Kontsevich has characterized the DG categories which should be thought of as \dfn[]{noncommutative spaces}.
For example, the DG category $\PfDG(X)$ is a noncommutative space, when $X$ is a smooth proper scheme over $K$.
Other examples come from symplectic geometry (Fukaya categories) and complex analytic geometry (deformation quantization modules).

Following (\cit{tabuada2011chow}), we define the category $\NChMot_K(\Lambda)$ of \dfn[]{noncommutative Chow motives} (with coefficients in a commutative ring $\Lambda$).
This serves as the universal invariant with respect to additive invariants like K-theory and cyclic homology.

\para
\emph{Perfect correspondences.}
We introduce a category of perfect correspondences, where objects are smooth projective varieties of finite type over a field $K$, and morphisms $X \to Y$ are perfect complexes on $X \times Y$.
We show that

\begin{parathm}
Let $X$ and $Y$ be smooth projective varieties over a field $K$.
$X$ and $Y$ are derived equivalent if and only if they are isomorphic in $\PfCorr_K$.
\end{parathm}

In particular, a functor $\Var_K \to \sC$ is an invariant of the derived category if and only if can also be defined on $\PfCorr_K$.

First, we show how the Grothendieck-Riemann-Roch theorem implies the existence of a canonical functor
  \[ M : \PfCorr_K \too \CHM_K(\Q)/\TateM. \tag{\ensuremath{\ast}} \]
Here $\CHM_K(\Q)/\TateM$ is the category of \dfn[]{Chow motives modulo Tate twists}, where Chow motives that differ by a Tate twist become isomorphic.
Therefore we see

\begin{parathm}
Let $X$ and $Y$ be smooth projective varieties over a field $K$, of dimensions $m$ and $n$, respectively.
If their derived categories $\Pf(X)$ and $\Pf(Y)$ are equivalent, then there is an isomorphism $M(X) \shiso M(Y)$ in $\CHM_K(\Q)/\TateM$.
\end{parathm}

We also deduce as a consequence

\begin{parathm}[\cit{orlov2005derived}]
Let $X$ and $Y$ be smooth projective varieties over a field $K$, of dimensions $m$ and $n$, respectively.
Suppose that $F : \Pf(X) \isoto \Pf(Y)$ is an equivalence of triangulated categories for which the corresponding perfect complex $\sE^\bull$ has support of dimension $n = \dim(X) = \dim(Y)$.
Then there is an isomorphism $M(X) \shiso M(Y)$ in $\CHM_K(\Q)$.
\end{parathm}

In particular the corresponding statement holds for any Weil cohomology theory, e.g. singular, de Rham, $\ell$-adic or crystalline cohomology.

Secondly, we show that there is a canonical functor
  \[ NM : \PfCorr_K \too \NChMot_K(\Lambda). \tag{\ensuremath{\ast\ast}} \]
This means that derived equivalence implies isomorphism of noncommutative Chow motives.

\begin{parathm}
Suppose that $X$ and $Y$ are smooth projective varieties over a field $K$.
If their derived categories $\Pf(X)$ and $\Pf(Y)$ are equivalent, then their noncommutative Chow motives are isomorphic.
\end{parathm}

In particular, it follows that all additive invariants like K-theory, Hochschild homology and cyclic homology are invariants of the derived category.
This was already established for K-theory by \cit{neeman2001k}, and for Hochschild homology by \cit{caldararu2010mukai}.

Finally, following \cit{tabuada2011chow} we show that there is a canonical fully faithful functor
  \[ R : \CHM_K(\Q)/\TateM \hooklong \NChMot_K(\Q) \]
and that the canonical functors $(\ast)$ and $(\ast\ast)$ respect this embedding.

\begin{paraprop}
The diagram
  \[ \begin{tikzcd}
      \PfCorr_K \arrow{d}{M}\arrow{dr}{NM}
      & \\
      \CHM_K(\Q)/\TateM \arrow[hookrightarrow]{r}{R}
      & \NChMot_K(\Q)
  \end{tikzcd} \]
commutes.
\end{paraprop}

\section*{Outline}

I have tried to keep this exposition as self-contained as possible, in the hope of being readable to anyone familiar with the basics of the theories of categories (\cit{borceux20081}) and schemes (\citcust{EGAI}{EGA I}, \citcust{EGAII}{EGA II}).

In the first chapter we recall the properties of the category of complexes in abelian category.
We introduce in detail the language of model categories, as it is very important to the study of DG categories.
We introduce the derived category as the homotopy category of a symmetric monoidal model structure on the category of complexes.
Finally we introduce the triangulated category of perfect complexes on a scheme, and define the usual derived functors.

The second chapter is a short introduction to the theory of DG categories; see (\cit{keller2006differential}) or (\cit{toen2011lectures}) for more details.
We also define the DG category of perfect complexes on a scheme.

In the third chapter we summarize the construction of Chow groups and characteristic classes; see (\cit{grothendieck1958intersections}) or (\cit{fulton1998intersection}) for more details.
Then we introduce the category of Chow motives and some variants.
See (\cit{andre2004introduction}) for more details.
Finally we introduce the noncommutative Chow motives following \cit{tabuada2011guided}.

In the last chapter we introduce the category of perfect correspondences, and show that an isomorphism in this category is the same thing as a derived equivalence of varieties.
Then we show how to construct functors from this category to the category of Chow motives up to Tate twists and to the category of noncommutative Chow motives.

\section*{Conventions}

\begin{enumerate}
  \item We ignore all set-theoretic issues and assume every category is small.
  The interested reader may use Grothendieck universes to deal with these issues.
  \item All rings are assumed to be commutative and unital.
  \item When $X$ and $Y$ are schemes over a ring $K$, we write $X \times Y$ for the fibred product over $K$.
\end{enumerate}

\section*{Acknowledgements}

I would like to thank my advisor H\'el\`ene Esnault for her guidance and careful reading.
I would like to thank Marc Levine for pointing out mistakes and for very helpful suggestions.
I would also like to thank Sasha Efimov, Sergey Galkin, Shane Kelly, Fernando Muro, Dmitri Orlov, Kay R\"ulling, Timo Sch\"urg and Bertrand To\"en for helpful discussions.
Finally I would like to thank Niko Naumann for the invitation to speak in the GK-Kolloquium in Universit\"at Regensburg, where I discussed a preliminary version of this work.
\chapter{Complexes of sheaves}
\label{chap:complexes}

\section{Categorical generalities}
\label{sec:categorical generalities}

\para[def:opposite category]
Let $\sC$ be a category.
We let $\sC^\circ$ denote the \dfn[]{opposite category of $\sC$}, whose objects are the same as those of $\sC$, morphisms are given for all objects $X$ and $Y$ by
  \[ \Hom_{\sC^\circ}(X, Y) = \Hom_\sC(Y, X), \]
and the law of composition is induced from $\sC$.

\para[def:small category]
A category $\sC$ is called \dfn[]{small} if its collection of objects is a \emph{set}.
We write $\Iso(\sC)$ for the set of isomorphism classes of objects in a small category.
We let $\Cat$ denote the category of small categories and functors.

To avoid dealing with set-theoretic issues, all categories that appear in the sequel will be assumed to be small.

\para[def:monoidal category]
A \dfn[monoidal category]{monoidal category} is a category $\sV$ together with a bifunctor $\cdot \otimes_\sV \cdot : \sV \times \sV \to \sV$, called the \dfn[]{tensor product}, a \dfn[]{unit object} $I_\sV$, and
  \begin{enumerate}
    \item for every three objects $X$, $Y$ and $Z$ a functorial \dfn[]{associativity isomorphism}
      \[ (X \otimes Y) \otimes Z \isoto X \otimes (Y \otimes Z), \]
    \item for every object $X$ a functorial \dfn[]{left unit isomorphism}
      \[ I_\sV \otimes X \isoto X, \]
    \item and for every object $X$ a functorial \dfn[]{right unit isomorphism}
      \[ X \otimes I_\sV \isoto X, \]
  \end{enumerate}
subject to the standard associativity and unit axioms.
See (\cit{borceux20082}, 6.1.1).
By abuse of notation we will write simply $\sV$ for a monoidal category and $\cdot \otimes \cdot$ and $I$ for the tensor product and unit object, respectively, when there is no risk of confusion.

A \dfn[]{symmetric monoidal category} is a monoidal category $\sV$ together with functorial isomorphisms
  \[ X \otimes_\sV Y \isoto Y \otimes_\sV X \]
for all objects $X$ and $Y$, satisfying the standard associativity, unit and symmetry axioms.
See (\cit{borceux20082}, 6.1.2).

\para[def:closed monoidal category]
A monoidal category $\sV$ is called \dfn[]{closed} if for every object $X$ the functor $\cdot \otimes_\sV X : \sV \to \sV$ admits a right adjoint $[X, \cdot] : \sV \to \sV$.
When $\sV$ is symmetric monoidal, one gets a bifunctor
  \[ [\cdot, \cdot] : \sV^\circ \times \sV \longrightarrow \sV \]
called the \dfn[]{internal hom functor}.
See (\cit{borceux20082}, 6.1.7).

\para[def:rigid symmetric monoidal category]
Let $\sV$ be a symmetric monoidal category.
A \dfn[]{dual of an object $X$} is an object $X^*$ together with an evaluation morphism $\ev_X : X^* \otimes_\sV X \to I_\sV$ and coevaluation morphism $\coev_X : I_\sV \to X \otimes X^*$ satisfying the usual triangle identities.
If $X$ has a dual, then it is unique up to unique isomorphism.

If every object of $\sV$ has a dual, we say that $\sV$ is a \dfn{rigid symmetric monoidal category}.

\para[def:enriched category]
Let $\sV$ be a monoidal category with tensor product $\otimes$ and unit object $I$.
Recall that a \dfn[enriched category over a monoidal category]{category $\sC^*$ enriched over $\sV$} is the data of
  \begin{enumerate}
    \item for every two objects $X$ and $Y$, an object $\Hom_{\sC^*}(X,Y)$ of $\sV$,
    \item for every three objects $X$, $Y$ and $Z$, a composition morphism
      \[ \Hom_{\sC^*}(Y,Z) \otimes \Hom_{\sC^*}(X, Y) \longrightarrow \Hom_{\sC^*}(X, Z) \]
    in $\sV$,
    \item and for every object $X$, a unit morphism $I \to \Hom_{\sC^*}(X, X)$ in $\sV$,
  \end{enumerate}
satisfying the standard associativity and unity axioms.
If $\sC$ is a category, an \dfn[]{enrichment of $\sC$ over $\sV$} is a category $\sC^*$ enriched over $\sV$ whose objects are the same as those of $\sC$.
See (\cit{borceux20082}, 6.2.1).

\para[def:additive category]
A category $\A$ is called \dfn[additive category]{additive} if it is enriched over the symmetric monoidal category of abelian groups, and admits finite products.

More generally, for $K$ a commutative ring, a \dfn[]{$K$-linear category} is a category that is enriched over $\Mod(K)$ and admits finite products.

\para[def:localization of a category]
Recall that the \dfn[localization of a category]{localization of a category $\sC$ at a class of morphisms $\mathscr{S}$} is a category $\sC[\mathscr{S}^{-1}]$ together with a functor $\gamma : \sC \to \sC[\mathscr{S}^{-1}]$ that maps every morphism $s \in \mathscr{S}$ to an isomorphism $\gamma(s)$ in $\sC[\mathscr{S}^{-1}]$, such that if $F : \sC \to \mathscr{D}$ is another functor that also maps morphisms in $\mathscr{S}$ to isomorphisms in $\mathscr{D}$, then there exists a unique functor $G : \sC[\mathscr{S}^{-1}] \to \mathscr{D}$ such that $F = G \circ \gamma$.
Equivalently, for every category $\mathscr{D}$ the canonical functor
  \begin{equation}\label{eq:localization of a category}
    \Hom(\gamma, \mathscr{D}) : \Hom(\sC[\mathscr{S}^{-1}], \mathscr{D}) \longrightarrow \Hom(\sC, \mathscr{D})
  \end{equation}
defines an equivalence between the category of functors $\sC[\mathscr{S}^{-1}] \to \mathscr{D}$ and the full subcategory of functors $\sC \to \mathscr{D}$ that map morphisms in $\mathscr{S}$ to isomorphisms.
A localization of $\sC$ at any class of morphisms $\sS$ always exists and is unique up to equivalence.
See (\cit{gabriel1967calculus}, I, 1.1).

%\para[def:tensor product of functors]
%Let $\sV$ be a monoidal category with tensor product $\otimes_\sV$.
%Let $F : \sC^\circ \to \sV$ and $\sC^\circ \to \sV$ be functors.
%The \dfn[tensor product of functors]{tensor product of $F$ and $G$} is defined as the coend
%  \[ F \otimes_\sA G = \int^X F(X) \otimes_\sV G(X). \]
%See (Borceux, II, 6.6.8) for the definition of a \dfn[]{coend}.

\para[def:karoubian category]
Let $\sC$ be a category and $p : X \to X$ a projector in $\sC$ (i.e. $p \circ p = p$).
If the equalizer of the pair $(\id_X, p)$ exists, one calls it the \dfn[]{image} of $p$ and denotes it $\Im(p)$.
If every projector in $\sC$ admits an image, then $\sC$ is called \dfn[karoubian category]{karoubian}.

\para[rmk:projectors in a karoubian category admit kernels]
Suppose $\sC$ is a karoubian category enriched over the category of abelian groups.
Then every projector $p : X \to X$ also admits a kernel, which is defined to be the equalizer of the pair $(p, 0)$; this is identical to the equalizer of the pair $(\id_X, \id_X - p)$, which exists because it is the image of the projector $\id_X - p$.

\para[def:karoubian envelope]
For any category $\sC$, there is an associated karoubian category $\Kar(\sC)$ whose objects are pairs $(X, p)$, where $X$ is an object of $\sC$ and $p : X \to X$ is a projector, and morphisms $(X, p) \to (Y, q)$ are morphisms $f : X \to Y$ of $\sC$ such that $q \circ f \circ p = f$.
The canonical functor $\sC \to \Kar(\sC)$ mapping an object $X$ to the pair $(X, \id_X)$ is clearly fully faithful.
Further, for any karoubian category $\sC'$ the canonical functor
  \[ \bHom(\Kar(\sC), \sC') \isoto \bHom(\sC, \sC'), \]
mapping $F : \Kar(\sC) \to \sC'$ to the composite $\sC \hook \Kar(\sC) \to \sC'$, is an \emph{equivalence} of categories.
We call $\Kar(\sC)$ the \dfn{karoubian envelope} of $\sC$.

\para[def:orbit category]
Let $\A$ be a category enriched over the category of abelian groups (\ref{def:enriched category}) and let $T : \A \to \A$ be an autoequivalence.
The \dfn[]{orbit category of $\A$ with respect to $T$} is the category $\A/T$ whose objects are the same as those of $\A$ and whose morphisms are given by
  \[ \Hom_{\A/T}(X, Y) = \bigoplus_{i \in \Z} \Hom_\A(X, T^i(Y)) \]
for objects $X$ and $Y$.
We define the composition law as follows: let $f = (f^i)_i : X \to Y$ and $g = (g^j)_j : Y \to Z$ be two morphisms in $\A/T$; the $k$-th component of the composite $g \circ f$ is the sum
  \begin{equation} \label{eq:composition in orbit category}
    (g \circ f)^k = \sum_{i+j=k} T^i(g^j) \circ f^i.
  \end{equation}
Note that this sum is finite because, by the definition of the direct sum of abelian groups, only finitely many $f^i$ (resp. $g^j$) are nonzero, for $i \in \Z$ (resp. $j \in \Z$).

There is a canonical projection functor $\pi_{\A/T} : \A \to \A/T$ which maps a morphism $f$ in $\A$ to the morphism $(\ldots, 0, f, 0, \ldots)$ in $\A/T$ which is zero on all components except the zeroth one.
%The orbit category $\A/T$ is clearly additive, and $\pi_{\A/T}$ is an additive functor.

\para[prop:right adjoint to projection functor of orbit category]
\begin{paraprop}
Let $\A$ be an additive category and $T : \A \to \A$ an additive autoequivalence.
Suppose that $\A$ admits arbitrary coproducts.
Then the projection functor $\pi = \pi_{\A/T} : \A \to \A/T$ admits a right adjoint
  \[ \tau = \tau_{\A/T} : \A/T \to \A. \]
\end{paraprop}

\begin{paraproof}
Define the functor $\tau : \A/T \to \A$ by mapping an object $X$ to the direct sum
  \[ X \mapsto \bigoplus_{i\in\Z} T^i(X) \]
and a morphism $f : X \to Y$ to the morphism
  \[ \bigoplus_{i \in \Z} T^i(X) \longrightarrow \bigoplus_{j \in \Z} T^j(Y) \]
that is induced by the morphisms
  \[ T^i(f^{j-i}) : T^i(X) \longrightarrow T^i(T^{j-i}(Y)) = T^j(Y) \]
for all $i,j \in \Z$.
One verifies easily that this is a well-defined functor.
That we have functorial adjunction isomorphisms
  \[ \Hom_{\A/T}(\pi(X), Y) \isoto \Hom_\A(X, \tau(X)) \]
follows from the additivity of the autoequivalence $T$.
\end{paraproof}

\para[prop:isomorphism in A/T gives isomorphism in A sometimes]
\begin{paraprop}
Let $\A$ be an additive category and $T : \A \to \A$ an additive autoequivalence.
Suppose that $f : X \isoto Y$ and $g : Y \isoto X$ are mutually inverse morphisms in the orbit category $\A/T$.
Write $f^i : X \to T^i(Y)$ for the $i$-th component of the morphism $f$ and $g^j : Y \to T^j(X)$ for the $j$-th component of $g$, and suppose that $f^i = 0$ for $i < 0$ and $g^j = 0$ for $j < 0$.
Then the morphisms $f^0 : X \to Y$ and $g^0 : Y \to X$ in $\A$ are mutual inverses.
\end{paraprop}

\begin{paraproof}
Since $f$ and $g$ are mutual inverses, the zeroth component of the morphism $g \circ f$ is the identity of $X$ in $\A$.
But by definition (\ref{eq:composition in orbit category}), one has
  \[ (g \circ f)^0 = \sum_{i+j = 0} T^i(g^j) \circ f^i = g^0 \circ f^0. \]
Similarly $f^0 \circ g^0 = (f \circ g)^0 = \id_Y$, from which the claim follows.
\end{paraproof}

\section{Abelian categories}
\label{sec:abelian categories}

\para
Let $\A$ be an additive category (\ref{def:additive category}).
One sees that every finite product is a coproduct, and every finite coproduct is a product.
In particular, the terminal object (the empty product) is the same as the initial object (the empty coproduct), and we call this the \dfn[]{zero object}.
See (\cit{grothendieck1957tohoku}, I, 1.3).

We refer to both finite products and coproducts as \dfn[]{direct sums}, and use the symbol $\oplus$.

\para
Let $\A$ be an additive category and $f : X \to Y$ a morphism in $\A$.
The \dfn[]{kernel of $f$}, if it exists, is the equalizer of the pair $(f, 0) : X \rightrightarrows Y$, i.e. the ``largest'' object $\Ker(f)$ of $\A$ with a morphism $\Ker(f) \to X$ such that the composite $\Ker(f) \to X\stackrel{f}{\to} Y$ is zero (in the sense that if $X' \to X$ is another morphism whose composite with $f$ is zero, then there is a unique morphism $X' \to \Ker(f)$ through which $X' \to X$ factors).
Dually, the \dfn[]{cokernel of $f$} is the coequalizer of the pair $(f, 0) : X \rightrightarrows Y$, i.e. the ``smallest'' object $\Coker(f)$ of $\A$ with a morphism $Y \to \Coker(f)$ such that the composite $X \stackrel{f}{\to} Y \to \Coker(f)$ is zero.

It is immediate from the definitions that when $f$ admits a kernel, the morphism $\Ker(f) \to X$ is a monomorphism, and that when it admits a cokernel, the morphism $Y \to \Coker(f)$ is an epimorphism.

\para[def:abelian category]
Suppose that $\A$ is an additive category that admits all kernels and
cokernels.
We define the \dfn[]{image $\Im(f)$ of $f$} as the kernel of the canonical morphism $Y \to \Coker(f)$, and the \dfn[]{coimage $\Coim(f)$ of $f$} as the cokernel of the canonical morphism $\Ker(f)  \to X$.
  \[ \begin{tikzcd}
    \Ker(f) \arrow[hookrightarrow]{r}
      & X \arrow{r}{f}\arrow[twoheadrightarrow]{d}
      & Y \arrow[twoheadrightarrow]{r}
      & \Coker(f)
    \\
      & \Coim(f) \arrow[dotted]{ru}\arrow[dotted]{r}
      & \Im(f) \ar[hookrightarrow]{u}
      &
  \end{tikzcd} \]
By the definition of $\Coim(f)$, since the composite $\Ker(f) \to X \to Y$ is zero, there exists a unique morphism $\Coim(f) \to Y$ through which $f$ factors.
Then the composite $X \to \Coim(f) \to Y \to \Coker(f)$ is equal to the composite $X \to Y \to \Coker(f)$, i.e. zero, and since $X \to \Coim(f)$ is an epimorphism, it follows that the composite $\Coim(f) \to Y \to \Coker(f)$ is also zero.
Now by the definition of $\Im(f)$ there is a unique morphism $\Coim(f) \to \Im(f)$ making the diagram commute.

If for every morphism $f$ in $\A$, the canonical morphism $\Coim(f) \to \Im(f)$ is an \emph{isomorphism}, then the category $\A$ is called \dfn[abelian category]{abelian}.

\para[def:short exact sequence in an abelian category]
Let $\A$ be an abelian category.
A \dfn[]{short exact sequence in $\A$} is a diagram
  \[ 0 \longrightarrow X \stackrel{f}{\longrightarrow} Y \stackrel{g}{\longrightarrow} Z \longrightarrow 0 \]
with $f$ a monomorphism, $\Im(f) \shiso \Ker(g)$, and $g$ an epimorphism.
A morphism of short exact sequences is a commutative diagram
  \[ \begin{tikzcd}
    0 \arrow{r}
      & X \arrow{r}{f}\arrow{d}{a}
      & Y \arrow{r}{g}\arrow{d}{b}
      & Z \arrow{r}\arrow{d}{c}
      & 0 \\
    0 \arrow{r}
      & X' \arrow{r}{f'}
      & Y' \arrow{r}{g'}
      & Z' \arrow{r}
      & 0
  \end{tikzcd} \]
The short exact sequences in $\A$ clearly form a category.

\para[def:generators of a category]
Let $\sC$ be a category.
A \dfn[]{family of generators of $\sC$} is a family of objects $(X_\alpha)_\alpha$ such that for any object $X$ in $\sC$ and any subobject $Y \ne X$, there exists an index $\alpha$ and a morphism $u : X_\alpha \to X$ which does not factor through a morphism $X_\alpha \to Y$.

\para[def:grothendieck abelian category]
A \dfn[Grothendieck abelian category]{Grothendieck abelian category} is an abelian category $\A$ satisfying the following:
  \begin{enumerate}
    \item $\A$ admits arbitary coproducts,
    \item the category of short exact sequences in $\A$ (\ref{def:short exact sequence in an abelian category}) admits filtered colimits (i.e. inductive limits),
    \item $\A$ admits a family of generators (\ref{def:generators of a category}).
  \end{enumerate}

\para[def:Grothendieck group of an exact category]
A subcategory $\A$ of an abelian category $\A'$ is called \dfn[]{thick} if for every short exact sequence in $\A'$ (\ref{def:short exact sequence in an abelian category})
  \[ 0 \longrightarrow X \longrightarrow Y \longrightarrow Z \longrightarrow 0 \]
with $X$ and $Z$ in $\A$, the object $Y$ also belongs to $\A$.
A category $\A$ is called \dfn[]{exact} if there exists a fully faithful functor $\A \to \A'$ identifying it with an additive thick full subcategory of an abelian category $\A'$.
There is an obvious notion of short exact sequences in any such category.

The \dfn[Grothendieck group of an exact category]{Grothendieck group of an exact category $\A$}, denoted $K_0(\A)$, is the free abelian group generated by the set of isomorphism classes of objects in $\A$, modulo relations identifying a class $[X]$ with the sum $[X'] + [X'']$ whenever there exists an exact sequence
  \[ 0 \to X' \to X \to X'' \to 0 \]
in $\A$.

\section{Complexes in an abelian category}
\label{sec:categories of complexes}

\para[def:categories of complexes]
Let $\A$ be an additive category.
A \dfn[complex in an additive category]{complex in $\A$} is a sequence $(X^n, d^n)$ of objects $X^n$ of $\A$ and morphisms $d^n : X^n \to X^{n+1}$ $(n \in \Z)$, such that $d^n \circ d^{n-1} = 0$ for all $n \in \Z$.
The morphism $d^n$ is called the $n$-th differential morphism of the complex.
A morphism of complexes $f : (X^n, d^n) \to (Y^n, e^n)$ is a sequence of morphisms $(f^n : X^n \to Y^n)$ $(n \in \Z)$, such that the diagram
  \[ \begin{tikzcd}
    \cdots \arrow{r}
      & X^n \arrow{r}{d^n}\arrow{d}{f^n}
      & X^{n+1} \arrow{r}{d^{n+1}}\arrow{d}{f^{n+1}}
      & \cdots \\
        \cdots \arrow{r}
      & Y^n \arrow{r}{e^n}
      & Y^{n+1} \arrow{r}{e^{n+1}}
      & \cdots
  \end{tikzcd} \]
commutes.
By abuse of notation we will write $X^\bull$ for the complex $(X^n, d^n)$, and $d_{X^\bull}^n$ for the differential morphisms, when there is no possibility of confusion.

The complexes in $\A$ clearly form an additive category $\Comp(\A)$.

\para
A complex $X^\bull$ is called \dfn[]{bounded below} (resp. \dfn[]{bounded above}) if there exists $k \in \N$ such that $X^n = 0$ for all $n < k$ (resp. for all $n > k$); $X^\bull$ is called \dfn[]{bounded} if it is both bounded below and bounded above.
One sees that the bounded complexes (resp. complexes bounded below, complexes bounded above) form a full additive subcategory of $\Comp(\A)$ which we denote $\Compb(\A)$ (resp. $\Comp^+(\A)$, $\Comp^-(\A)$).

\para[rmk:embedding into category of complexes]
There is a canonical fully faithful functor $\A \to \Comp(\A)$ which maps an object $X$ of $\A$ to the complex $X^\bull$ defined by $X^0 = X$ and $X^n = 0$ for $n \ne 0$ (and all differential morphisms zero).
Hence we may identify $\A$ with its image in $\Comp(\A)$.

\para[rmk:category of complexes is abelian]
If $\A$ is an \emph{abelian} category, then one may show that $\Comp(\A)$ is also abelian (\cit{weibel1995introduction}, Theorem 1.2.3).

\para[rmk:symmetric monoidal structure on category of complexes]
Let $\A$ be a symmetric monoidal abelian category, with tensor product $\otimes_\A$ and unit object $I$.
We define a tensor product on $\Comp(\A)$ as follows: for two complexes $X^\bull$ and $Y^\bull$ the $n$-th component of the complex $X^\bull \otimes Y^\bull$ is given by
  \[ (X^\bull \otimes Y^\bull)^n = \bigoplus_{i+j = n} X^i \otimes_\A Y^j, \]
and the differential morphisms $d^n : (X^\bull \otimes Y^\bull)^n \to (X^\bull \otimes Y^\bull)^{n+1}$ are given by
  \[ x \otimes y \quad \mapsto \quad d_{X^\bull}^i(x) \otimes y + (-1)^i x \otimes d_{Y^\bull}^j(y) \]
where $x \in X^i$ and $y \in Y^j$.
This defines a symmetric monoidal structure on $\Comp(\A)$, where the unit object is the complex with the single object $I$ in degree zero.

\para[def:quasi-isomorphism of complexes]
Let $X^\bull$ be a complex in an abelian category $\A$.
For each $n \in \Z$, one obtains using the equality $d^n \circ d^{n-1} = 0$ a canonical monomorphism $\Im(d^{n-1}) \hook \Ker(d^n)$; we define the \dfn[]{$n$-th cohomology object of $X^\bull$}, denoted $H^n(X^\bull)$, as the cokernel:
  \[ H^n(X^\bull) = \Coker( \Im(d^{n-1}) \hook \Ker(d^n) ). \]

A morphism $f : X^\bull \to Y^\bull$ of complexes in $\A$ induces morphisms $H^n(f) : H^n(X^\bull) \to H^n(Y^\bull)$ on the cohomology objects.
If these morphisms are \emph{isomorphisms} for each $n \in \Z$, then $f$ is called a \dfn[quasi-isomorphism of complexes]{quasi-isomorphism}.

\para[rmk:complexes of modules]
Let $A$ be a unital commutative ring.
Recall that the category $\Mod(A)$ of $A$-modules is abelian and symmetric monoidal.
By (\ref{rmk:category of complexes is abelian}) and (\ref{rmk:symmetric monoidal structure on category of complexes}), the category $\Comp(\Mod(A))$ of complexes of $A$-modules is also abelian and symmetric monoidal.
We will abuse notation and write simply $\Comp(A)$ for $\Comp(\Mod(A))$.

%%%%%%%%%%%%%%%%%%%%%%%%%%%%%%%%%%%%%%%%%%%%%%%%%%%%%%%%

\section{Model categories}
\label{sec:model-categories}

\para[def:retract]
A morphism $f : X \to X'$ in a category $\sC$ is called a \dfn[retract of a morphism]{retract of a morphism $g : Y \to Y'$} if there exists a commutative diagram
  \[ \begin{tikzcd}
    X \arrow{r}{i}\arrow{d}{f}
      & Y \arrow{r}{r}\arrow{d}{g}
      & X \arrow{d}{f}
    \\
    X' \arrow{r}{i'}
      & Y' \arrow{r}{r'}
      & X'
  \end{tikzcd} \]
such that $r \circ i = \id_X$ and $r' \circ i' = \id_{X'}$.

\para[def:model category]
Let $\sC$ be a category.  A \dfn[model structure on a category]{model structure on $\sC$} consists of three classes of morphisms, called \dfn[]{weak equivalences}, \dfn[]{fibrations} and \dfn[]{cofibrations}, such that each class contains all identity morphisms and is closed under composition, and
\begin{enumerate}[labelsep=10pt]
  \item[MC-1] $\sC$ admits arbitrary limits and colimits;
  \item[MC-2] for morphisms $f : X \to Y$ and $g : Y \to Z$, if two of the morphisms $f$, $g$ and $g \circ f$ are weak equivalences, then so is the third;
  \item[MC-3] if a morphism $f : X \to X'$ is a retract of a morphism $g : Y \to Y'$ that is a weak equivalence (resp. fibration, cofibration), then so is $f$.
  \item[MC-4] if there is a commutative diagram
    \[ \begin{tikzcd}
      A \arrow{r}\arrow{d}{i}
        & X \arrow{d}{p}
      \\
      B \arrow{r}\arrow[dotted]{ru}
        & Y
    \end{tikzcd} \]
  of solid arrows in $\sC$, then the lift $B \to X$ exists when $i$ is a cofibration and $p$ is a fibration and weak equivalence, and also when $i$ is a cofibration and weak equivalence and $p$ is a fibration.
  \item[MC-5] any morphism $f : X \to Y$ can be factored as
    \[ X \stackrel{i}{\longrightarrow} X' \stackrel{p}{\longrightarrow} Y \]
  with $i$ a cofibration and weak equivalence and $p$ a fibration, and also as
    \[ X \stackrel{i'}{\longrightarrow} X'' \stackrel{p'}{\longrightarrow} Y \]
  with $i'$ a cofibration and $p$ a fibration and weak equivalence.
\end{enumerate}

By abuse of language we will call $\sC$ a \dfn[]{model category} when there is a model structure given, and there is no risk of confusion.
We will call a morphism that is both a cofibration (resp. fibration) and weak equivalence a \dfn[]{trivial cofibration} (resp. \dfn[]{trivial fibration}).

\para[def:cofibrant and fibrant objects]
Let $\sC$ be a model category.
Note that by (MC-1) $\sC$ contains an initial object $\inobj$ (resp. final object $\finobj$) as a colimit (resp. limit) of the diagram $\Empty \to \sC$, where $\Empty$ denotes the empty category with no objects.
As a consequence of the Yoneda lemma, $\inobj$ and $\finobj$ are determined up to unique isomorphism.

If for some object $X$ of $\sC$ the unique morphism $\inobj \to X$ (resp. $X \to \finobj$) is a cofibration (resp. fibration), then $X$ is called a \dfn[]{cofibrant object of $\sC$} (resp. \dfn[]{fibrant object of $\sC$}).
If $X$ is both cofibrant and fibrant, then it is called an \dfn[]{interior object} of $\sC$.

\para[def:cofibrant and fibrant replacement]
In general, by (MC-5) one may always factor $\inobj \to X$ as a composite
    \[ \begin{tikzcd}
      \inobj \arrow{rr}\arrow{rd}{i}
        &
        & X
      \\
        & Q(X) \arrow{ru}{p}
        &
    \end{tikzcd} \]
with $i$ a cofibration, $p$ a trivial fibration.
The object $Q(X)$ is then a cofibrant object and we call it a \dfn[]{cofibrant replacement of $X$}.
Dually, one may factor the morphism $X \to \finobj$ as a composite
    \[ \begin{tikzcd}
      X \arrow{rr}\arrow{rd}{i}
        &
        & \finobj
      \\
        & R(X) \arrow{ru}{p}
        &
    \end{tikzcd} \]
with $i$ a trivial cofibration and $p$ a fibration.
We call $R(X)$ a \dfn[]{fibrant replacement of $X$}.

%\para
%\begin{paraprop}
%In a model category $\sC$, the class of cofibrations (resp. trivial cofibrations) is stable under co-base change, and the class of fibrations (resp. trivial fibrations) is stable under base change.
%\end{paraprop}

\para[def:diagonal and codiagonal morphisms]
If $\sC$ is a category admitting finite colimits, let $A$ be an object of $\sC$ and consider the coproduct $A \sqcup A$.
There is a unique morphism $\nabla : A \sqcup A \to A$, called the \dfn[]{codiagonal morphism of $A$}, such that the diagram
    \[ \begin{tikzcd}
      \inobj \arrow{r}\arrow{d}
        & A \arrow{d}\arrow[bend left]{rdd}{\id_A}
        &
      \\
        A \arrow{r}\arrow[bend right]{rrd}{\id_A}
        & A \sqcup A \arrow[dotted]{rd}{\nabla}
        &
      \\
        &
        &
        A
    \end{tikzcd} \]
commutes.

Dually, when $\sC$ admits finite limits, for any object $X$ there is a unique morphism $\Delta : X \to X \times X$, called the \dfn[]{diagonal morphism of $X$}, such that the diagram
    \[ \begin{tikzcd}
      X \arrow[bend left]{rrd}{\id_X}\arrow[dotted]{rd}{\Delta}\arrow[bend right]{rdd}{\id_X}
        &
        &
      \\
        & X \times X \arrow{r}\arrow{d}
        & X \arrow{d}
      \\
        & X \arrow{r}
        & \finobj
    \end{tikzcd} \]
commutes.

\para[def:left and right homotopic morphisms]
Let $\sC$ be a model category.
Let $f$ and $g$ be morphisms $A \to B$ in $\sC$.
We say that $f$ is \dfn[]{left homotopic to $g$}, and write $f \sim^\ell g$, if there exists a commutative diagram
  \[ \begin{tikzcd}
    A \sqcup A \arrow{r}{f+g}\arrow{rd}{i}\arrow{d}{\nabla}
      & B
    \\
    A
      & A' \arrow{u}\arrow{l}{s}
  \end{tikzcd} \]
with $s$ a weak equivalence.
Dually $f$ is \dfn[]{right homotopic to $g$}, and write $f \sim^r g$ if there exists a commutative diagram
  \[ \begin{tikzcd}
    B' \arrow{rd}{p}
      & B \arrow{l}{s}\arrow{d}{\Delta}
    \\
    A \arrow{u}\arrow{r}{f+g}
      & B \times B
  \end{tikzcd} \]
with $s$ a weak equivalence.

\para[def:cylinder and path objects]
By (MC-5) there is a factorization of the codiagonal morphism $\nabla$ as $A \sqcup A \stackrel{i}{\to} A' \stackrel{s}{\to} A$ with $i$ a cofibration and $s$ a weak equivalence.
Such an object $A'$ is called a \dfn[cylinder object in a model category]{cylinder object for $A$} and we denote it by abuse of notation as $\Cyl(A)$ (it is unique up to weak equivalence).

Dually there exists a factorization of the diagonal morphism $\Delta : X \stackrel{s}{\to} X' \stackrel{p}{\to} X \times X$ by a weak equivalence $s$ and a fibration $p$.
The object $X'$ is called a \dfn[]{path object of $\sC$} and denoted by abuse of notation by $\Path(X)$.

\para[def:left and right homotopies]
If $f,g : A \rightrightarrows X$ are morphisms, a \dfn[]{left homotopy from $f$ to $g$} is a morphism $\varphi : \Cyl(A) \to X$ such that the diagram
    \[ \begin{tikzcd}
      A \sqcup A \arrow{rr}{f \sqcup g}\arrow{rd}{i}
        &
        & X
      \\
        & \Cyl(A) \arrow{ru}{\varphi}
        &
    \end{tikzcd} \]
commutes.

If $f,g : A \rightrightarrows X$ are morphisms, a \dfn[]{right homotopy from $f$ to $g$} is a morphism $\varphi : A \to \Path(X)$ such that the diagram
    \[ \begin{tikzcd}
      A \arrow{rr}{(f, g)}\arrow{rd}{\varphi}
        &
        & X \times X
      \\
        & \Path(X) \arrow{ru}{p}
        &
    \end{tikzcd} \]
commutes.

\para[rmk:left homotopic morphisms give a left homotopy]
Suppose $f$ and $g$ are morphisms $A \to B$ such that $f$ is left homotopic (resp. right homotopy) to $g$.
Then one can show that there exists a left homotopy $\varphi : \Cyl(A) \to B$ (resp. right homotopy $\varphi : A \to \Path(B)$).
(See \cit{quillen1967homotopical}, I, Lemma 1.)

\para[rmk:left and right homotopy for interior objects]
When $X$ is cofibrant (resp. $Y$ is fibrant), the left homotopy relation $\sim^\ell$ (resp. the right homotopy relation $\sim^r$) defines an equivalence relation on the set of morphisms $X \to Y$.
Also, if in this case two morphisms $f, g : X \rightrightarrows Y$ are left homotopic (resp. right homotopic), then they are also right homotopic (resp. left homotopic).
(See \cit{quillen1967homotopical}, I, Lemma 4 and Lemma 5, (i).)

\para[def:homotopy classes of morphisms of interior objects]
For two objects $X$ and $Y$, the right homotopy relation $\sim^r$ (resp. left homotopy relation $\sim^\ell$ defines an equivalence relation which is the intersection of all equivalence relations containing it.
We let $\pi^r(X, Y)$ (resp. $\pi^\ell(X, Y)$) be the set of equivalence classes of $\Hom_\sC(X, Y)$ with respect to this equivalence relation.
When $X$ and $Y$ are both cofibrant and fibrant objects, the relations $\sim^r$ and $\sim^\ell$ are already equivalence relations, and they coincide (\ref{rmk:left and right homotopy for interior objects}).
In this case we call the relation simply \dfn[]{homotopy} and we write $\pi(X, Y)$ for the set of equivalence classes of $\Hom_\sC(X, Y)$ with respect to homotopy.

\para[def:pi_int of a model category]
Let $\piint(\sC)$ denote the category whose objects are interior objects of $\sC$ (\ref{def:cofibrant and fibrant objects}) and whose morphisms are given by
  \[ \Hom_{\piint(\sC)}(X, Y) = \pi(X, Y) \]
for any two interior objects $X$ and $Y$.
One checks that the composition in $\sC$ induces a well-defined composition on the homotopy classes.

Similarly we let $\picof(\sC)$ (resp. $\pifib(\sC)$) be the category whose objects are cofibrant objects (resp. fibrant objects) of $\sC$ and whose morphisms are given by $\pi^\ell(X, Y)$ (resp. $\pi^r(X, Y)$).

\para[def:homotopy category of a model category]
The \dfn[homotopy category of a model category]{homotopy category of a model category $\sC$} is the localization (\ref{def:localization of a category}) at its class of weak equivalences, denoted $\gamma : \sC \to \Ho(\sC)$.

\para[rmk:functor mapping weak equivalences to isomorphisms identifies homotopic maps]
Let $\sC$ be a model category and $F : \sC \to \sD$ a functor to an arbitrary category $\sD$.
Suppose that $F$ maps all weak equivalences of $\sC$ into isomorphisms in $\sD$.
Then one sees that for any left or right homotopic morphisms $f,g : A \rightrightarrows B$ in $\sC$, the images $F(f)$ and $F(g)$ coincide in $\sD$.
(See \cit{quillen1967homotopical}, I, Lemma 8, (i))

\para[rmk:functors induced by cofibrant and fibrant replacement]
For every object $X$ of $\sC$, choose a cofibrant replacement $Q(X)$ and a fibrant replacement $R(X)$, so that there are trivial fibrations $p_X : Q(X) \to X$ and trivial cofibrations $i_X : X \to R(X)$.
(If $X$ is already cofibrant (resp. fibrant), then take $p_X$ (resp. $i_X$) to be the identity morphism.)
For a morphism $f : X \to Y$, the diagram
  \[ \begin{tikzcd}
    \inobj \arrow{r}\arrow{d}
      & Q(Y) \arrow{d}{p_Y}
    \\
    Q(X) \arrow[dotted]{ru}\arrow{r}{f \circ p_X}
      & Y
  \end{tikzcd} \]
commutes trivially and there is by (MC-4) a morphism $Q(f) : Q(X) \to Q(Y)$; it can be shown to be unique up to left homotopy.
It follows that $Q(\id_X) \sim^\ell \id_{Q(X)}$ and that for another morphism $g : Y \to Z$, one has $Q(g \circ f) \sim^\ell Q(g) \circ Q(f)$.
By (\ref{rmk:left and right homotopy for interior objects}) they are also right homotopic.
Hence one gets a functor
  \[ \bar Q : \sC \to \picof(\sC) \]
which maps $X \mapsto Q(X)$ and $f : X \to Y$ to the equivalence class $[Q(f)]$ in $\pi^r(Q(X), Q(Y))$.
Dually one gets a functor
  \[ \bar R : \sC \to \pifib(\sC). \]

\para[thm:fundamental theorem of model categories]
Let $\sC$ be a model category.
The localization functor $\gamma : \sC \to \Ho(\sC)$ induces by (\ref{rmk:functor mapping weak equivalences to isomorphisms identifies homotopic maps}) a canonical functor $\bar\gamma : \piint(\sC) \to \Ho(\sC)$ (\ref{def:pi_int of a model category}).

\begin{parathm}[fundamental theorem of model categories]
Let $\sC$ be a model category.
For all objects $X$ and $Y$, there are canonical isomorphisms
  \[ \Hom_{\Ho(\sC)}(\gamma(X), \gamma(Y)) \iso \pi(Q(R(X)), Q(R(Y))), \]
and the functor $\bar\gamma$ is an equivalence of categories
  \[ \piint(\sC) \isoto \Ho(\sC). \]
\end{parathm}

See (\cit{hovey2007model}, Theorem 1.2.10).

%\para[rmk:isomorphisms in homotopy category are weak equivalences]
%Let $\sC$ be a model category.
%It is straightforward to prove that if $f$ is a morphism in $\sC$ such that $\gamma(f)$ is an isomorphism in $\Ho(\sC)$, then $f$ is a weak equivalence (Hovey, Theorem 1.2.10, (iv)).

\para
Let $F : \sC \to \sD$ be a functor between two categories $\sC$ and $\sD$.
The \dfn[derived functor]{left derived functor of $F$ with respect to a functor $\gamma : \sC \to \sC'$} is a functor $\L^\gamma F : \sC' \to \sD$ such that there exists a morphism of functors $\epsilon : \L^\gamma F \circ \gamma \to F$ such that for any functor $G : \sC' \to \sD$ with a morphism $\zeta : G \circ \gamma \to F$, there exists a morphism $\theta : G \to \L^\gamma F$ such that the diagram
  \[ \begin{tikzcd}
    G \circ \gamma \arrow{d}{\zeta}\arrow{r}{\theta \ast \gamma}
      & \L^\gamma \circ \gamma \arrow{ld}{\epsilon}
    \\
    F
      &
  \end{tikzcd} \]
commutes, where $\theta \ast \gamma : G \circ \gamma \to \L^\gamma F \circ \gamma$ is the morphism canonically induced by $\theta$.

Similarly we define the \dfn[]{right derived functor of $F$ with respect to $\gamma : \sC \to \sC'$} as a functor $\R^\gamma F : \sC' \to \sD$ with a morphism $\eta : F \to \R^\gamma F \circ \gamma$ such that for any functor $G : \sC' \to \sD$ and morphism $\zeta : F \to G \circ \gamma$, there exists a morphism $\theta : \R^\gamma F \circ \gamma \to G \circ \gamma$ making the diagram
  \[ \begin{tikzcd}
    F \arrow{d}{\eta}\arrow{rd}{\zeta}
      &
    \\
    \R^\gamma F \circ \gamma \arrow{r}{\theta \ast \gamma}
      &
    G \circ \gamma
  \end{tikzcd} \]
commute.

\para
Let $\sC$ be a model category.
We will be interested in derived functors of functors  with respect to the canonical functor $\gamma : \sC \to \Ho(\sC)$; we will denote these simply by $\L F$ and $\R F$ for a functor $F : \sC \to \sD$.

\para[rmk:construction of left derived functor]
If $F : \sC \to \sD$ is a functor that maps weak equivalences between cofibrant objects to isomorphisms in $\sD$, then one may construct its left derived functor $\L F : \Ho(\sC) \to \sD$ as follows.
Let $\bar Q : \sC \to \picof(\sC)$ be the functor induced by cofibrant replacement (\ref{rmk:functors induced by cofibrant and fibrant replacement}).
Note that $F$ induces a functor $\bar F : \picof(\sC) \to \sD$ because by (\ref{rmk:functor mapping weak equivalences to isomorphisms identifies homotopic maps}), $F$ identifies left homotopic morphisms between cofibrant objects.
The composite functor $\bar F \circ \bar Q : \sC \to \sD$ maps weak equivalences in $\sC$ to isomorphisms in $\sD$ : in fact, $Q(s) : Q(X) \to Q(Y)$ fits by construction into a diagram
  \[ \begin{tikzcd}
      & Q(Y) \arrow{d}{p_Y}
    \\
    Q(X) \arrow{ru}{Q(s)}\arrow{r}{s \circ p_X}
      & Y;
  \end{tikzcd} \]
as $p_Y$ and $s \circ p_X$ are weak equivalences, so is $Q(s)$ (MC-2); hence $F(Q(s))$ is an isomorphism by assumption.
Hence by the definition of localization there exists a functor $\tilde F : \Ho(\sC) \to \sD$ such that the diagram
  \[ \begin{tikzcd}
    \sC \arrow{r}{\bar F \circ \bar Q}\arrow{d}{\gamma}
      & \sD
    \\
    \Ho(\sC) \arrow{ru}{\tilde F}
      &
    \end{tikzcd} \]
commutes.
One sees that $\tilde F$ is a left derived functor of $F$, with the morphism $\epsilon : \tilde F \circ \gamma \to F$ defined so that $\epsilon(X) : F(Q(X)) \to F(X)$ is the morphism $F(p_X)$ for every object $X$.

\para[def:total adjoint functors]
Let $F : \sC \to \sD$ be a functor of model categories.
We define the \dfn[total derived functor]{total derived functor of $F$} as the left-derived functor of the composite $\gamma_{\sD} \circ F : \sC \to \sD \to \Ho(\sD)$.
Similarly for a functor $G : \sD \to \sC$ we define the \dfn[]{total derived functor of $G$} as the right-derived functor of $\gamma_{\sC} \circ G : \sD \to \sC \to \Ho(\sC)$.

\para[thm:quillen adjunctions]
Let $\sC$ and $\sD$ be model categories and $(F, G) : \sC \rightleftarrows \sD$ be adjoint functors such that $F$ preserves cofibrations and weak equivalences, and $G$ preserves fibrations and weak equivalences.
Then one can prove that the total derived functors $\L F : \Ho(\sC) \to \Ho(\sD)$ and $\R G : \Ho(\sD) \to \Ho(\sC)$ are adjoint.
See (\cit{hovey2007model}, 1.3.10).

When $F$ and $G$ satisfy the above conditions, we call $F$ a \dfn[]{left Quillen functor}, $G$ a \dfn[]{right Quillen functor}, and the pair $(F, G)$ a \dfn[]{Quillen adjunction}.

\para[rmk:functoriality of derived functors]
Let $\sC$, $\sD$ and $\sE$ be model categories and let $F : \sC \to \sD$ and $F' : \sD \to \sE$ be left Quillen functors.
There is a canonical isomorphism of functors
  \[ \L F' \circ \L F \isoto \L (F' \circ F). \]
Dually, if $G : \sD \to \sC$ and $G' : \sE \to \sD$ are right Quillen functors, then there is a canonical isomorphism of functors
  \[ \R G \circ \R G' \isoto \R (G \circ G'). \]
See (\cit{hovey2007model}, Theorem 1.3.7).

\para[def:model category over a symmetric monoidal category]
Let $\sC$ be a model category and $\sV$ a symmetric monoidal model category, with tensor product $\otimes_\sV$, unit object $I_\sV$ and internal hom $[ \cdot, \cdot ]_\sV$.
A \dfn[model category enriched over a symmetric monoidal category]{model category enriched over $\sV$}, or a \emph{$\sV$-model category}, is a model category $\sC$ with an enrichment $\sC^*$ over $\sV$ (\ref{def:enriched category}) and two bi-functors $\cdot \otimes \cdot : \sV \times \sC \to \sC$ and $[\cdot, \cdot] : \sV \times \sC \to \sC$ such that
  \begin{enumerate}
    \item for objects $X$ of $\sC$ and objects $A$ and $B$ of $\sV$, there are functorial isomorphisms
  \[ (A \otimes_\sV B) \otimes X \isoto A \otimes (B \otimes X) \]
and
  \[ I_\sV \otimes X \isoto X \]
in $\sC$;
    \item for objects $X$ and $Y$ of $\sC$ and objects $A$ of $\sV$, there are functorial isomorphisms in $\sV$
  \[ \Hom_\sC(A \otimes X, Y) \isoto [A, \Hom_{\sC^*}(X, Y)]_\sV \isofrom \Hom_\sC(X, [A, Y]); \]
    \item the bi-functor $\cdot \otimes \cdot : \sV \times \sC \to \sC$ is \dfn[]{Quillen}, i.e. for every cofibration $A \to B$ in $\sV$ and cofibration $X \to Y$ in $\sC$, the induced morphism
	  \[ B \otimes X \bigsqcup_{A \otimes X} A \otimes Y \longrightarrow B \otimes Y \]
	is also a cofibration; further, it is a \emph{trivial} cofibration whenever $A \to B$ or $X \to Y$ is.
  \end{enumerate}
By abuse of notation we write $(\sC, \sC^*)$ or even just $\sC$ for the $\sV$-model category $(\sC, \sC^*, \cdot\otimes\cdot, [\cdot,\cdot])$ when there is no risk of confusion.

\para[rmk:homotopy category of a model category enriched over a symmetric monoidal category]
Let $\sC$ be a $\sV$-model category.
The homotopy category $\Ho(\sC)$ has a natural enrichment $\Ho^*(\sC)$ over $\Ho(\sV)$, defined by the formula
  \[ \Hom_{\Ho^*(\sC)} (X, Y) = \Hom_{\sC^*}(Q(X), R(Y)) \]
for objects $X$ and $Y$, where $Q$ and $R$ denote cofibrant and fibrant replacements, respectively.
See (\cit{hovey2007model}, Theorem 4.3.4).

As a consequence, when $\sV$ is the symmetric monoidal category $\C(K)$, one has
  \[ H^0(\Hom_{\sC^*}(X, Y)) \iso \Hom_{\Ho(\sC)}(X, Y) \]
for $X$ cofibrant and $Y$ fibrant.

\para[rmk:model structure on Cat]
Recall that $\Cat$ denotes the category of small categories and functors between them (\ref{def:small category}).
There is a model structure on $\Cat$ where the weak equivalences are equivalences of categories, and every category is both fibrant and cofibrant.
The cylinder object of a category $\sC$ with respect to this model structure is the product $\sC \times \mbf{I}$, where $\mbf{I}$ is the category with two objects 0 and 1 and an isomorphism $0 \to 1$.
See (\cit{rezk1996model}).

Let $F$ and $G$ be two functors $\sC \rightrightarrows \msc{D}$.
It is not difficult to see that an isomorphism of functors $F \isoto G$ is the same as a homotopy $\sC \times \mbf{I} \to \msc{D}$ between them.
Hence by the fundamental theorem of model categories (\ref{thm:fundamental theorem of model categories}), the homotopy category $\Ho(\Cat)$ is canonically equivalent to the category $[\Cat]$ whose morphisms are isomorphism classes of functors.

%%%%%%%%%%%%%%%%%%%%%%%%%%%%%%%%%%%%%%%%%%%%%%%%%%%%%%%%

\newpage

\section{Derived categories}
\label{sec:derived categories}

\para[def:triangles in a category]
Let $\sA$ be an additive category and $T : \sA \to \sA$  an additive autoequivalence.
A \dfn[]{triangle in $\sA$} (with respect to $T$) is a diagram
  \[ X \stackrel{u}{\longrightarrow}
     Y \stackrel{v}{\longrightarrow}
     Z \stackrel{w}{\longrightarrow}
     T(X). \]
We write such a diagram as a tuple $(X, Y, Z, u, v, w)$.
A morphism of triangles $(X, Y, Z, u, v, w) \to (X', Y', Z', u', v', w')$ is a tuple of morphisms $f : X \to X'$, $g : Y \to Y'$, $h : Z \to Z'$ such that the diagram
  \[ \begin{tikzcd}
    X \arrow{r}{u}\arrow{d}{f}
      & Y \arrow{r}{v}\arrow{d}{g}
      & Z \arrow{r}{w}\arrow{d}{h}
      & T(X) \arrow{d}{T(f)}
    \\
    X' \arrow{r}{u'}
      & Y' \arrow{r}{v'}
      & Z' \arrow{r}{w'}
      & T(X')
  \end{tikzcd} \]
commutes.
One obtains a category of triangles in $\sA$ with respect to $T$.

\para[def:triangulated category]
A \dfn[triangulated structure on a category]{triangulated structure on an additive category $\sA$} consists of an additive autoequivalence $T : \sA \to \sA$, called the \dfn[]{translation functor}, and a family of triangles, called the \dfn[]{distinguished triangles}, satisfying the following axioms.
\begin{enumerate}[labelsep=10pt]
  \item[TR-1]  Every triangle isomorphic to a distinguished triangle is distinguished.  For every morphism $u : X \to Y$ there exists some distinguished triangle $(X, Y, Z, u, v, w)$.  The triangle $(X, X, 0, \id_X, 0, 0)$ is distinguished.
  \item[TR-2]  A triangle $(X, Y, Z, u, v, w)$ is distinguished if and only if the triangle $(Y, Z, T(X), v, w, -T(u))$ is distinguished.
  \item[TR-3]  If $(X, Y, Z, u, v, w)$ and $(X', Y', Z', u', v', w')$ are distinguished, for every pair of morphisms $f : X \to X'$ and $g : Y \to Y'$ there exists a morphism $g : Z \to Z'$ such that $(f, g, h)$ is a morphism of triangles.
  \item[TR-4]  If $(X, Y, Z', u, i, j)$, $(Y, Z, X', v, j, s)$ and $(X, Z, Y', w, k, t)$ are distinguished triangles such that $w = v \circ u$, then there exist morphisms $f : Z' \to Y'$ and $g : Y' \to X'$ such that  $(\id_X, v, f)$ and $(u, \id_Z, g)$ are morphisms of triangles, and $(Z', Y', X', f, g, T(i) \circ s)$ is a distinguished triangle.
\end{enumerate}

By abuse of language we will simply say a category $\sA$ is \dfn[]{triangulated}, leaving the triangulated structure implicit when there is no risk of confusion.

\para[def:triangulated subcategory]
Let $\sA$ be a triangulated category.
A full subcategory $\sA' \sub \sA$ is called \dfn[]{triangulated} if every distinguished triangle of $\sA$ having two of its objects in $\sA'$ is isomorphic to a triangle where all three objects are in $\sA'$.

\para[def:triangulated functor]
Let $\sA$ and $\sB$ be two additive categories with fixed additive autoequivalences $T_\sA : \sA \to \sA$ and $T_\sB : \sB \to \sB$.
An additive functor $F : \sA \to \sB$ is called \dfn[]{graded} if there exists an isomorphism of functors $F \circ T_\sA \isoto T_\sB \circ F$.
If $\sA$ and $\sB$ are further triangulated, then $F$ is called \dfn[]{triangulated} if it is additive, graded and maps distinguished triangles of $\sA$ to distinguished triangles of $\sB$.

We let $\TriCat$ denote the category of triangulated categories and triangulated functors.

\para[def:compact object of a triangulated category]
Let $\sA$ be a triangulated category.
An object $X$ is called \dfn[]{compact} if the functor $\Hom_\sA(X, \cdot)$ commutes with arbitrary coproducts.

\para[def:translation and cones in category of complexes]
Let $\sA$ be an additive category and $X^\bull$ a complex in $\Comp(\sA)$.
We define the \dfn[]{$n$-th translation of $X^\bull$} as the complex $X^\bull[1]$ whose $n$-th component is $X^{n+1}$ and whose $n$-th differential morphism is the \emph{opposite} of the $(n+1)$-th differential morphism of $X^\bull$.
This clearly defines an autoequivalence $[1] : \Comp(\sA) \to \Comp(\sA)$.

For a morphism $f : X^\bull \to Y^\bull$ we define the \dfn[]{cone of $f$} as the complex $\Cone(f)$ whose $n$-th component is $X^{n+1} \oplus Y^n$ and whose differential morphisms $d_{\Cone(f)}^n : X^{n+1} \otimes Y^n \to X^{n+2} \otimes Y^{n+1}$ map a pair $(x, y')$ to $(f(x) + d_{Y^\bull}(y'), -d_{X^\bull}(x))$.
There are canonical morphisms $Y^\bull \to \Cone(f)$ and $\Cone(f) \to X^\bull[1]$ which are degree-wise inclusions and projections, respectively.

\para[def:homotopy category of complexes]
Let $\sA$ be an additive category.
Consider the category $\HComp(\sA)$, called the \dfn[]{homotopy category of complexes in $\sA$}, whose objects are complexes in $\sA$ and morphisms are homotopy classes of morphisms in $\Comp(\sA)$.
This category is also additive and there is an obvious functor $\Comp(\sA) \to \HComp(\sA)$.
The translation automorphism $[1] : \Comp(\sA) \to \Comp(\sA)$ induces an additive automorphism on $\HComp(\sA)$ which we denote again by $[1]$.

We also define the full subcategory $\HCompb(\sA) \sub \HComp(\sA)$ of \emph{bounded} complexes.
Clearly the translation automorphism induces a well-defined additive automorphism $[1] : \HCompb(\sA) \to \HCompb(\sA)$ as well.

\para[rmk:triangulated structure on homotopy category of complexes]
Let $\sA$ be an additive category.
Consider the family of triangles in $\Comp(\sA)$ of the form
  \[ X^\bull \stackrel{f}{\longrightarrow} Y^\bull \longrightarrow \Cone(f) \longrightarrow X^\bull[1] \]
for some morphism $f$.
Let $\sT$ denote the family of triangles in $\HComp(\sA)$ \emph{isomorphic} to a triangle in the image of this family by the canonical functor $\Comp(\sA) \to \HComp(\sA)$.
The additive automorphism $[1]$ and the family $\sT$ define a triangulated structure on $\HComp(\sA)$ (\cit{kashiwara1990sheaves}, Proposition 1.4.4).

It is not difficult to see that the category $\HCompb(\sA)$ is a \emph{triangulated} subcategory of $\HComp(\sA)$.

\para[rmk:localization of K(A) is triangulated]
Let $\sA$ be an abelian category.
Consider the localization (\ref{def:localization of a category}) of $\HComp(\sA)$ at the class of quasi-isomorphisms (\ref{def:quasi-isomorphism of complexes}).
The auto-equivalence $[1]$ on $\HComp(\sA)$ clearly induces an auto-equivalence on $\HComp(\sA)$ which we denote again by $[1]$.
One can prove that the data of this auto-equivalence together with the family of triangles that are isomorphic to the image of a distinguished triangle of $\HComp(\sA)$ define a triangulated structure on the localization, and that the localization functor is triangulated.

\para[rmk:localization of C(A) vs K(A)]
Let $\sA$ be an abelian category.
Let $\sQ$ denote the class of quasi-isomorphisms in $\Comp(\sA)$ and $\sQ'$ the image of this class under the canonical functor $\Comp(\sA) \to \HComp(\sA)$.
Consider the localizations $\gamma : \Comp(\sA) \to \Comp(\sA)[\sQ^{-1}]$ and $\gamma' : \HComp(\sA) \to \HComp(\sA)[\sQ'^{-1}]$.
The composite $\Comp(\sA) \to \HComp(\sA) \to \HComp(\sA)[\sQ'^{-1}]$ clearly induces a functor $\Comp(\sA)[\sQ^{-1}] \to \HComp(\sA)[\sQ^{-1}]$, by the definition of localization, such that the diagram
  \[ \begin{tikzcd}
    \Comp(\sA) \arrow{r}\arrow{d}{\gamma}
      & \HComp(\sA) \arrow{d}{\gamma'} \\
    \Comp(\sA)[\sQ^{-1}]  \arrow[dotted]{r}
      & \HComp(\sA)[\sQ'^{-1}]
  \end{tikzcd} \]
commutes.
Further we can define a canonical functor $\HComp(\sA) \to \Comp(\sA)[\sQ^{-1}]$ mapping a homotopy class of a morphism $f : X^\bull \to Y^\bull$ to its image $\gamma(f)$ in $\Comp(\sA)[\sQ^{-1}]$, since one can show that homotopic morphisms induce equal morphisms in the localization.
This induces a functor $\HComp(\sA)[\sQ'^{-1}] \to \Comp(\sA)[\sQ^{-1}]$.
Now it is clear that the composite $\Comp(\sA) \to \HComp(\sA) \to \HComp(\sA)[\sQ'^{-1}]$ is also a localization of $\Comp(\sA)$ at $\sQ$: any functor $F : \Comp(\sA) \to \sD$ mapping quasi-isomorphisms to isomorphisms factors through $\HComp(\sA)[\sQ^{-1}] \to \Comp(\sA)[\sQ^{-1}] \to \sD$.
Hence the localizations $\Comp(\sA)[\sQ^{-1}]$ and $\HComp(\sA)[\sQ'^{-1}]$ are canonically identified in such a way that the diagram above commutes.

\para[def:derived category]
For an abelian category $\sA$, the localization of $\Comp(\sA)$ with respect to the class of quasi-isomorphisms is called the \dfn[derived category of an abelian category]{derived category of $\sA$} and denoted
  \[ \gamma_\sA : \sA \to \D(\sA). \]
By (\ref{rmk:localization of C(A) vs K(A)}) and (\ref{rmk:localization of K(A) is triangulated}), it follows that there is a canonical triangulated structure on $\D(\sA)$ induced by the triangulated structure on $\HComp(\sA)$.
We will always consider $\D(\sA)$ as a triangulated category.

Similarly one gets the \dfn[]{bounded derived category of $\A$}, denoted $\Db(\A)$, by taking the localization of $\Compb(\A)$.
Again, the triangulated structure on $\HCompb(\A)$ induces a triangulated structure on $\Db(\A)$.

\para[def:Grothendieck group of a triangulated category]
Let $\A$ be a triangulated category.
The \dfn[Grothendieck group of a triangulated category]{Grothendieck group of $\A$}, denoted $K_0(\A)$, is the free abelian group generated by the set of isomorphism classes of objects in $\A$, modulo relations identifying a class $[X]$ with the sum $[X'] + [X'']$ whenever there exists a distinguished triangle
  \[ X' \to X \to X'' \to X'[1] \]
in $\A$.

%%%%%%%%%%%%%%%%%%%%%%%%%%%%%%%%%%%%%%%%%%%%%%%%%%%%%%%%

\section{Perfect complexes of sheaves}
\label{sec:perfect complexes of sheaves}

\para[rmk:category of O_X-modules]
Let $(X, \O_X)$ be a ringed space.
Recall that the category $\Mod(\O_X)$ of $\O_X$-modules is closed symmetric monoidal (\cit{godement1973topologie}, II, 2).
It is also a Grothendieck abelian category (\ref{def:grothendieck abelian category}), generated by the family $(i_U)_! (\O_X|_U)$ where $U$ ranges among open subsets of $X$, $i_U: U \hook X$ denotes the inclusion morphism).
(Recall that $(i_U)_! : \Mod(\O_X|_U) \to \Mod(\O_X)$ is the left adjoint  to the restriction functor $\Mod(\O_X) \to \Mod(\O_X|_U)$.)
See (\cit{grothendieck1957tohoku}, 3.1.1).

In the sequel we will write $\Comp(X)$ and $\D(X)$ for the categories $\Comp(\Mod(\O_X))$ and $\D(\Mod(\O_X))$, respectively.

\para[rmk:model structure on C(Mod(O_X))]
By (\cit{cisinski2009local}, Example 2.3 and Theorem 2.5), there exists a model structure on the category $\Comp(X)$ of complexes of $\O_X$-modules, where weak equivalences are quasi-isomorphisms.
In particular, the homotopy category of $\Comp(X)$ with respect to model structure is precisely the derived category $\D(X)$.

Further, this model structure is compatible with the symmetric monoidal structure (\ref{rmk:symmetric monoidal structure on category of complexes}) by (\cit{cisinski2009local}, Example 3.1 and Proposition 3.2).
In particular, the tensor product descends to a derived bifunctor
  \[ \cdot \otimes^\bL \cdot : \D(X) \times \D(X) \longrightarrow \D(X) \]
on the homotopy categories.

\para[rmk:derived functors on C(X)]
Let $f : (X, \O_X) \to (Y, \O_Y)$ be a morphism of ringed spaces.
By (\cit{cisinski2009local}, Theorem 2.14), the direct and inverse image functors form a Quillen adjunction
  \[ (f^*, f_*) : \C(Y) \rightleftarrows \C(X) \]
and descend by (\ref{thm:quillen adjunctions}) to an adjoint pair of derived functors
  \[ (\L f^*, \R f_*) : \D(Y) \rightleftarrows \D(X). \]

By (\ref{rmk:functoriality of derived functors}) one immediately gets canonical isomorphisms
  \[ \L f^* \circ \L g^* \isoto \L (g \circ f)^* \]
and
  \[ \R g_* \circ \R f_* \isoto \R (g \circ f)_*. \]

\para[def:Db(Coh(O_X))]
Let $X$ be a noetherian scheme.
Consider the full abelian subcategory $\Coh(\O_X) \sub \Mod(\O_X)$ of coherent $\O_X$-modules.
The canonical fully faithful functor
  \[ \Db(\Coh(\O_X)) \hooklong \Db(\Mod(\O_X)). \]
identifies $\Db(\Coh(\O_X))$ with the full subcategory $\Dbcoh(\Mod(\O_X)) \sub \Db(\Mod(\O_X))$ of complexes whose cohomology objects are coherent sheaves (\citcust{SGA6}{SGA VI}, Exp. II, Corollaire 2.2.2.1).

\para[def:perfect complexes]
Let $X$ be a scheme.
A \dfn[perfect complex on a scheme]{perfect complex on $X$} is a complex $\sE^\bull$ of $\O_X$-modules such that there exists an open cover $(U_\alpha)_\alpha$ of $X$ and for each $\alpha$ a quasi-isomorphism from $\sE^\bull|_{U_\alpha}$ to a bounded complex of free $\O_X|_{U_\alpha}$-modules of finite rank.

We let $\Pf(X)$ denote the full subcategory of $\D(X)$ of perfect complexes.
This is a triangulated subcategory (\ref{def:triangulated subcategory}) that is stable under derived tensor product, derived inverse image, and derived direct image of proper morphisms.
Also, note that one has an inclusion $\Pf(X) \sub \Dbcoh(\Mod(\O_X))$; when $X$ is smooth, one can further show that every bounded complex with coherent cohomology is quasi-isomorphic to a perfect complex, so $\Pf(X) \iso \Dbcoh(X)$.
See (\citcust{SGA6}{SGA 6}, VI, Exp. I, \S 4).
% for derived direct image, http://ens.math.univ-montp2.fr/~toen/sga6.pdf

Therefore, by (\ref{def:Db(Coh(O_X))}), when $X$ is smooth we may identify $\Pf(X)$, the triangulated category of perfect complexes, with $\Db(\Coh(\O_X))$, the bounded derived category of coherent sheaves.
\chapter{DG categories}
\label{chap:DG categories}

\section{DG categories}
\label{sec:DG categories}

\para[def:DG category]
Let $K$ be a commutative ring.
A \dfn[DG category]{differential graded category over $K$} (or simply \emph{DG category}) is a category $\AA$ enriched (\ref{def:enriched category}) over the symmetric monoidal category $\Comp(K)$ of complexes of $K$-modules (\ref{rmk:complexes of modules}).
By abuse of language we will often leave the base ring $K$ implicit when there is no risk of confusion.

\para[def:opposite DG category]
Let $\AA$ be a DG category.
The \dfn[]{opposite DG category of $\AA$} is the DG category $(\AA)^\circ$ whose objects are the same and complexes of morphisms are given by
  \[ \Hom_{(\AA)^\circ}(X, Y) = \Hom_\AA(Y, X) \]
for all objects $X$ and $Y$.

\para[def:DG category of complexes]
There is a DG category of complexes of $K$-modules, denoted $\Cdg(K)$ and defined as follows: for two complexes $X^\bull$ and $Y^\bull$, define $\Hom_{\Cdg(K)}(X^\bull, Y^\bull)$ as the complex whose $n$-th component is the $K$-module
  \[ \Hom_{\Cdg(K)}(X^\bull, Y^\bull)^n = \prod_{i \in \Z} \Hom(X^i, Y^{i+n}), \]
and whose differential morphisms $d^n : \Hom_{\Cdg(K)}(X^\bull, Y^\bull)^n \to \Hom_{\Cdg(K)}(X^\bull, Y^\bull)^{n+1}$ are defined by
  \[ (f^i)_{i \in \Z} \quad \mapsto \quad (d_{Y^\bull}^{i+n} \circ f^i - (-1)^n f^{i+1} \circ d_{X^{\bull}}^i)_{i \in \Z}. \]

\para[def:homotopy category of a DG category]
Let $\AA$ be a DG category.
The \dfn[cohomology category of a DG category]{cohomology category of $\AA$} is the category $\H^*(\AA)$ whose objects are the same as those of $\AA$, morphisms are given by
  \[ \Hom_{\H^*(\AA)}(X, Y) = H^*(\Hom_\AA(X, Y)) = \bigoplus_{i \in \Z} H^i(\Hom_\AA(X, Y)) \]
for any two objects $X$ and $Y$, and the composition law is induced from $\AA$.
Similarly the \dfn[homotopy category of a DG category]{homotopy category of $\AA$} is the category $\H^0(\AA)$ whose morphisms are given by
  \[ \Hom_{\H^0(\AA)}(X, Y) = H^0(\Hom_\AA(X, Y)) \]
for any two objects $X$ and $Y$, and composition law is again induced from $\AA$.

\para[def:DG morphism in a DG category]
A morphism $f : X \to Y$ in $\AA$ is called \dfn[]{closed} if its image under the differential of the complex $\Hom_\AA(X, Y)$ is zero.
It is called a \dfn[DG morphism]{DG morphism} if it is closed and of degree zero.
Note that a DG morphism induces a morphism in the homotopy category $\H^0(\AA)$.

\para[def:DG functor]
A \dfn{DG functor} $F : \AA \to \BB$ is the data of
  \begin{enumerate}
    \item for every object $X$ of $\AA$ an object $F(X)$ of $\BB$,
    \item for every two objects $X$ and $Y$ of $\AA$, a morphism
      \[ \Hom_\AA(X, Y) \longrightarrow \Hom_{\BB}(F(X), F(Y)) \]
    in $\Comp(\Mod(K))$,
  \end{enumerate}
subject to the following conditions:
  \begin{enumerate}
    \item for any three objects $X$, $Y$ and $Z$ of $\AA$, the diagram
      \[ \begin{tikzcd}
        \Hom_\AA(X, Y) \otimes \Hom_\AA(Y, Z) \arrow{r}\arrow{d}
      & \Hom(X, Z) \arrow{d}
        \\
        \Hom_{\BB}(F(X), F(Y)) \otimes \Hom_{\BB}(F(Y), F(Z)) \arrow{r}
          & \Hom_{\BB}(F(X), F(Z))
      \end{tikzcd} \]
    commutes (the horizontal morphisms are composition);
    \item for every object $X$ of $\AA$, the diagram
      \[ \begin{tikzcd}
        K \arrow{r}\arrow{rd}
	  & \Hom(X, X) \arrow{d}
	\\
	  & \Hom_{\BB}(F(X), F(X))
      \end{tikzcd} \]
    commutes (the morphisms with domain $K$ are the unit morphisms of $X$ and $F(X)$).
  \end{enumerate}

One gets a category $\DGCat_K$ of DG categories over $K$.
By abuse of notation we will write simply $\DGCat$ when there is no ambiguity.

\para[def:tensor product of DG categories]
Let $\AA$ and $\BB$ be DG categories over $K$.
The \dfn[tensor product of DG categories]{tensor product of $\AA$ and $\BB$}, denoted $\AA \otimes_K \BB$, is the category whose objects are pairs $(X, Y)$, with $X$ an object of $\AA$ and $Y$ an object of $\BB$, and complexes of morphisms are given by
  \[ \Hom_{\AA \otimes_K \BB}((X, Y), (X', Y')) = \Hom_{\AA}(X, Y) \otimes_K \Hom_{\BB}(X', Y') \]
for all objects $X$ and $X'$ of $\AA$ and $Y$ and $Y'$ of $\BB$.
By abuse of notation we will write $\AA \otimes \BB$ when there is no risk of confusion.

This tensor product makes $\DGCat_K$ a symmetric monoidal category.

\para[def:interior DG category]
Let $(\sC, \sC^*)$ be a model category enriched over $\C(K)$ (\ref{def:model category over a symmetric monoidal category}).
The enrichment $\sC^*$ is by a definition a DG category over $K$, which we call the \dfn[]{big DG category associated to $\sC$}.
Let $\IntDG(\sC, \sC^*)$ be the full sub-DG category of $\sC^*$ consisting of objects that are cofibrant and fibrant in $\sC$.
Then by (\ref{rmk:homotopy category of a model category enriched over a symmetric monoidal category}) we have the equivalence of categories
  \[ \H^0(\IntDG(\sC, \sC^*)) \iso \Ho(\sC). \]
We call $\IntDG(\sC, \sC^*)$ the \dfn[interior DG category]{interior DG category of $\sC$}.
By abuse of notation we write $\IntDG(\sC) = \IntDG(\sC, \sC^*)$ when there is no risk of confusion.

\para[def:homotopically flat DG category]
A DG category $\AA$ is called \dfn[]{homotopically flat} if the functor $\cdot \otimes_K \Hom_\AA (X, Y) : \C(K)  \to \C(K)$ preserves quasi-isomorphisms for all objects $X$ and $Y$.

%%%%%%%%%%%%%%%%%%%%%%%%%%%%%%%%%%%%%%%%%%%%%%%%%%%%%%%%%%%%%%%%%%%%%%%%%%%

\section{The derived category of a DG category}

\para[def:DG modules]
Let $\AA$ and $\BB$ be DG categories over $K$.
An \dfn[module over a DG category]{$\AA$-module} is a DG functor $(\AA)^\circ \to \Cdg(K)$.
A morphism of $\AA$-modules $M \to N$ is a morphism of DG functors.
We let $\Mod(\AA)$ denote the category of $\AA$-modules.

An \dfn[bimodule over DG categories]{$\AA$-$\BB$-bimodule} is a DG functor $\AA \otimes (\BB)^\circ \to \Cdg(K)$, i.e. a module over $(\AA)^\circ \otimes \BB$.

\para[def:Yoneda DG functor]
The canonical DG functor $\AA \to \ModDG(\AA)$ mapping an object $X$ to the $\AA$-module $\Hom_\AA(\cdot, X)$ is called the \dfn[]{Yoneda DG functor of $\AA$}.
The $\AA$-module $\Hom_\AA(\cdot, X)$ is called the \dfn[]{$\AA$-module represented by the object $X$}.

\para[def:composition of DG bimodules]
Let $M : \AA \otimes (\BB)^\circ \to \Cdg(K)$ be a $\AA$-$\BB$-bimodule and $N : \BB \otimes (\CC)^\circ \to \Cdg(K)$ a $\BB$-$\CC$-bimodule.
We define the composition $N \circ M$, a $\AA$-$\CC$-bimodule, as the equalizer of the two canonical morphisms
  \[ \bigsqcup_{Y', Y'' \in \BB} N(Y'', Z) \otimes \Hom_\BB(Y', Y'') \otimes M(X, Y') \rightrightarrows \bigsqcup_{Y \in \BB} N(Y, Z) \otimes M(X, Y). \]
See (\cit{borceux20082}, 6.2.11) for details.

\para[def:derived category of a DG category]
The \dfn[]{homology functor of an $\AA$-module $M$} is the functor $\H^*(M) : \H^*(\AA) \to \GrMod(K)$ mapping an object $X$ to the graded $K$-module $H^*(M(X))$.
A morphism $\varphi : M \to N$ induces a morphism of functors $\H^*(M) \to \H^*(N)$; if the latter is an isomorphism, $\varphi$ is called a \dfn[quasi-isomorphism of DG modules]{quasi-isomorphism of $\AA$-modules}.

The \dfn[derived category of a DG category]{derived category of a DG category $\AA$} is the localization (\ref{def:localization of a category}) of $\Mod(\AA)$ at the class of quasi-isomorphisms, denoted
  \[ \gamma_\AA : \Mod(\AA) \to \D(\AA). \]

\para[def:DG category of DG modules]
The category $\Mod(\AA)$ has a canonical enrichment over $\Comp(K)$, induced by the differential graded structure of $\Cdg(K)$.
Further there is a (projective) model structure on $\Mod(\AA)$ which is compatible with this enrichment (\cit{keller2006differential}, Theorem 3.2), so we get a $\Comp(K)$-model category and in particular a DG category $\IntDG(\Mod(\AA))$ whose homotopy category is equivalent to $\D(\AA)$ (\ref{rmk:homotopy category of a model category enriched over a symmetric monoidal category}).

\para[rmk:derived category of DG modules is triangulated]
Let $\AA$ be a DG category.
There exists a canonical triangulated structure (\ref{def:triangulated category}) on its derived category $\D(\AA)$ (\cit{keller2006differential}, Lemma 3.3).
Hence we will always view $\D(\AA)$ as a triangulated category.

%\para[rmk:derived category of DG modules is symmetric monoidal]
%The symmetric monoidal structure on $\Mod(\AA)$ (\ref{rmk:tensor product of DG modules}) induces a symmetric monoidal structure on the derived category $\D(\AA)$, with the derived tensor product $\cdot \otimes^\bL_\AA \cdot : \D(\AA) \to \D(\AA)$ defined by taking a cofibrant replacement first, as usual.

\para[def:quasi-functor]
Let $\AA$ and $\BB$ be DG categories.
Assume that $\AA$ is homotopically flat over $K$ (\ref{def:homotopically flat DG category}).
A \dfn[quasi-functor of DG categories]{quasi-functor $\AA \to \BB$} is an $\AA$-$\BB$-bimodule $M$ such that for every object $X$ of $\AA$, the $\BB$-module $M(X, \cdot)$ is quasi-isomorphic to $\Hom_\BB(\cdot, Y)$ for some object $Y$ of $\BB$.
%the functor $\cdot \otimes^\bL_{\AA} M : \D(\AA) \to \D(\BB)$ takes representable $\AA$-modules to objects isomorphic to representable $\BB$-modules.
Quasi-functors $\AA \to \BB$ define a full triangulated subcategory of $\D((\AA)^\circ \otimes \BB)$ which we denote by $\QsiHom(\AA, \BB)$.
We also let $\QsiHomDG(\AA, \BB)$ denote the full sub-DG category of quasi-functors in $\IntDG(\Mod((\AA)^\circ \otimes \BB))$.
By (\ref{def:DG category of DG modules}) one has a canonical equivalence of categories
  \begin{equation}\label{eq:QsiHom* and QsiHom}
    \H^0(\QsiHomDG(\AA, \BB)) \iso \QsiHom(\AA, \BB).
  \end{equation}

%In general, when $\AA$ is not homotopically flat we may take a cofibrant replacement of it and use the definitions above.

%\para[rmk:quasi-functors induce functors on homotopy categories]
%Any quasi-functor $M : (\AA)^\circ \otimes \BB \to \Cdg(K)$ induces a DG functor $\H^0(M) : \H^0(\AA) \to \H^0(\BB)$.
%This functor maps an object $X$ of $\AA$ to the object $M_X$ of $\BB$ representing the $\BB$-module isomorphic to $\Hom_\AA(X, \cdot) \otimes^\bL_\AA M$, and maps a morphism $f : X \to Y$ to the image of $\id_{M_Y}$ under the composite
%  \[ \Hom_\BB(M_Y, M_Y) \isoto M(Y, M_Y) \stackrel{f}{\longrightarrow} M(X, M_Y) \isoto \Hom_\BB(M_X, M_Y). \]
%One may verify that this is indeed a well-defined functor on the homotopy categories.

\para[def:composition of quasi-functors]
Let $M \in \QsiHom(\AA, \BB)$ and $N \in \QsiHom(\BB, \CC)$ be quasi-functors of DG categories.
It is not difficult to see that the composition $N \circ M$ as bimodules (\ref{def:composition of DG bimodules}) is a quasi-functor from $\AA$ to $\CC$.

%%%%%%%%%%%%%%%%%%%%%%%%%%%%%%%%%%%%%%%%%%%%%%%%%%%%%%%%%%%%%%%%%%%%%%%%%%%

\section{Quasi-equivalence of DG categories}
\label{sec:quasi-equivalence of DG categories}

\para[def:quasi-equivalence of DG categories]
A morphism $F : \AA \to \BB$ is called a \dfn{quasi-equivalence of DG categories} if for every two objects $X$ and $Y$ of $\AA$, the canonical morphism
  \[ \Hom_\AA(X, Y) \longrightarrow \Hom_{\BB}(F(X), F(Y)) \]
is a quasi-isomorphism of complexes (\ref{def:quasi-isomorphism of complexes}), and the induced functor
  \[ \H^0(F) : \H^0(\AA) \longrightarrow \H^0(\BB) \]
is an equivalence on the homotopy categories.

\para[rmk:quasi-equivalent model structure on DGCat]
There is a model structure on the category $\DGCat_K$ where the weak equivalences are quasi-equivalences (\ref{def:quasi-equivalence of DG categories}), and fibrations are morphisms $F : \AA \to \BB$ such that
  \begin{enumerate}
    \item for any two objects $X$ and $Y$ of $\AA$, the induced morphism of complexes
      \[ \Hom_\AA(X, Y) \longrightarrow \Hom_{\BB}(F(X), F(Y)) \]
    is surjective in every degree,
    \item for every object $X$ and morphism $g : F(X) \to Y'$ in $\BB$ that becomes an isomorphism in $\H^0(\BB)$, there exists a morphism $f : X \to Y$ in $\AA$ such that $F(f) = g$ and that becomes an isomorphism in $\H^0(\AA)$.
  \end{enumerate}
See (\cit{tabuada2005structure}, Th. 1.8 and Prop. 1.13).

\para[def:HoDGQe]
The homotopy category with respect to the model structure (\ref{def:homotopy category of a model category}), i.e. the localization of $\DGCat_K$ at the class of quasi-equivalences, is denoted $\HoDGQe_K$.
By the fundamental theorem of model categories (\ref{thm:fundamental theorem of model categories}), one may compute morphisms in $\HoDGQe_K$ by taking homotopy classes of morphisms between cofibrant replacements.
(It is clear that every DG category is fibrant.)

\para[rmk:derived H^0 functor]
Consider the functor $\H^0 : \DGCat_K \to \Cat$ assigning to a DG category $\AA$ its homotopy category $\H^0(\AA)$ (\ref{def:homotopy category of a DG category}).
By definition, all quasi-equivalences are mapped to isomorphisms, so it descends to a functor $\HoDGQe_K \to \Cat$.
Composing with the localization functor $\Cat \to \Ho(\Cat)$ (\ref{rmk:model structure on Cat}), we get a canonical functor on the homotopy categories, which we also denote $\H^0$ by abuse of notation.
  \[ \begin{tikzcd}
    \DGCat \arrow{r}{\H^0}\arrow{d} & \Cat \arrow{r} & \Ho(\Cat) \\
    \HoDGQe \arrow[dashed]{ru}\arrow[dashed,swap]{rru}{\H^0} & &
  \end{tikzcd} \]

\para[thm:morphisms in HoDGQe]
\begin{parathm}[To\"en]
For any two DG categories $\AA$ and $\BB$, there is a canonical isomorphism
  \[ \Iso(\QsiHom(\AA, \BB)) \isoto \Hom_{\HoDGQe}(\AA, \BB) \]
mapping a quasi-functor $M$ to the functor $X \mapsto \Hom_\AA(\cdot, X) \otimes^\bL_\AA M$.
It is functorial in the sense that composition of quasi-functors (\ref{def:composition of quasi-functors}) corresponds to composition in $\HoDGQe$.
\end{parathm}

See (\cit{toen2011lectures}, 4.1, Corollary 1).

\para[thm:internal hom on HoDGQe]
The tensor product of DG categories descends to a derived tensor product $\cdot\otimes^\bL\cdot$ on $\HoDGQe$ (given by first taking cofibrant replacements).
This gives a symmetric monoidal structure which is in fact closed (\ref{def:closed monoidal category}):

\begin{parathm}[To\"en]
The category $\HoDGQe$ admits an internal hom functor $\RHomDG : \HoDGQe \times \HoDGQe \to \HoDGQe$, so that one has functorial isomorphisms
  \[ \RHomDG(\AA, \cdot) \isoto \cdot \otimes^\bL \AA \]
for all DG categories $\AA$.
Further, if $\AA$ and $\BB$ are DG categories with $\AA$ homotopically flat (\ref{def:homotopically flat DG category}), then there is an equivalence of DG categories
  \[ \RHomDG(\AA, \BB) \iso \QsiHomDG(\AA, \BB). \]
\end{parathm}

See (\cit{toen2007homotopy}, 6.1).

\para[rmk:homotopy category of RHom]
As an immediate corollary of (\ref{thm:internal hom on HoDGQe}) one has canonical equivalences of categories
  \[ \H^0(\RHomDG(\AA, \BB)) \isoto \QsiHom(\AA, \BB) \]
by (\ref{eq:QsiHom* and QsiHom}), when $\AA$ is homotopically flat.

\para[def:quasi-equivalent DG categories]
Two DG categories $\AA$ and $\BB$ are called \dfn[]{quasi-equivalent} if there exist DG categories $\sC_1^\bull, \ldots, \sC_n^\bull$ and a chain of quasi-equivalences
  \[ \AA \leftarrow \sC_1^\bull \to \sC_2^\bull \leftarrow \cdots \leftarrow \sC_n^\bull \to \BB. \]

\begin{paraprop}
Two DG categories $\AA$ and $\BB$ are quasi-equivalent if and only if they are isomorphic in the homotopy category $\HoDGQe$.
\end{paraprop}

%%%%%%%%%%%%%%%%%%%%%%%%%%%%%%%%%%%%%%%%%%%%%%%%%%%%%%%%%%%%%%%%%%%%%%%%%%%

\section{Morita equivalence of DG categories}
\label{sec:Morita equivalence of DG categories}

\para
A DG functor $F : \AA \to \BB$ clearly induces a canonical functor $F_* : \Mod(\BB) \to \Mod(\AA)$.
One can show that $F_*$ admits a left adjoint $F_! : \Mod(\AA) \to \Mod(\BB)$, and further that the pair $(F_!, F_*)$ is a Quillen adjunction (\ref{thm:quillen adjunctions}).
In particular one gets an adjoint pair of derived functors
  \[ (\bL F_!, \bR F_*) : \D(\AA) \rightleftarrows \D(\BB). \]
See (\cit{toen2011lectures}, Section 3.2).

\para[def:Morita equivalence of DG categories]
Let $F : \AA \to \BB$ be a DG functor.
If the derived functor $\bR F_* : \D(\BB) \to \D(\AA)$ is an equivalence, then we call $F$ a \dfn[Morita equivalence of DG categories]{Morita equivalence}.

\para[def:Morita model structure on DGCat]
The category $\DGCat_K$ admits a model structure where the weak equivalences are Morita equivalences (\cit{tabuada2005invariants}, Theorem 5.3).
We let $\HoDGMo_K$ denote the homotopy category of $\DGCat_K$ with respect to this model structure.

\para[rmk:symmetric monoidal structure on HoDGMo]
The tensor product of DG categories induces a derived tensor product $\cdot \otimes^\bL_K \cdot$ on $\HoDGMo_K$, which is defined as usual by taking cofibrant replacements first.
This gives a symmetric monoidal structure (\cit{tabuada2005invariants}, Remark 5.11).

\para[rmk:quasi-equivalences are Morita equivalences]
When $F : \AA \to \BB$ is a quasi-equivalence of DG categories, one can show that the induced functor $\bR F_* : \D(\BB) \to \D(\AA)$ is an equivalence of categories (\cit{toen2007homotopy}, Prop 3.2).
In other words, every quasi-equivalence of categories is a Morita equivalence, and we have a commutative diagram of categories
  \[ \begin{tikzcd}
    \DGCat \arrow{r}\arrow{d}
      & \HoDGMo \\
    \HoDGQe \arrow[dotted]{ru}
  \end{tikzcd} \]

\para[rmk:HoDGMo is a full subcategory of HoDGQe]
By choosing for every DG category a fibrant replacement with respect to the Morita model structure (\ref{def:Morita equivalence of DG categories}), one gets a well-defined functor $\HoDGMo \to \HoDGQe$ (all fibrant replacement functors are identified in the homotopy category).

\begin{paraprop}
The canonical functor $\HoDGQe \to \HoDGMo$ (\ref{rmk:quasi-equivalences are Morita equivalences}) is left adjoint to the functor $\HoDGMo \to \HoDGQe$ induced by fibrant replacement.
Further, the latter functor is fully faithful and identifies $\HoDGMo$ with the full subcategory of $\HoDGQe$ whose objects are DG categories that are fibrant with respect to the Morita model structure.
\end{paraprop}

\para[rmk:morphisms in HoDGMo]
As an immediate consequence of (\ref{thm:morphisms in HoDGQe}) and (\ref{rmk:HoDGMo is a full subcategory of HoDGQe}) one has the corollary

\begin{paracor}
For any two DG categories $\AA$ and $\BB$, there is a canonical functorial isomorphism
  \[ \Iso(\QsiHom(\AA, \BB)) \isoto \Hom_{\HoDGMo}(\AA, \BB). \]
\end{paracor}

%%%%%%%%%%%%%%%%%%%%%%%%%%%%%%%%%%%%%%%%%%%%%%%%%%%%%%%%%%%%%%%%%%%%%%%%%%%

\section{Pretriangulated DG categories}
\label{sec:pretriangulated DG categories}

\para[def:pretriangulated DG category]
Let $\AA$ be a DG category over $K$.
The \dfn[]{$n$-translation of an object $X$} is an object $X[n]$ for which there are functorial isomorphisms
  \[ \Hom_\AA(Y, X[n]) \isoto \Hom_\AA(Y, X)[n] \]
in $\C(K)$ for every object $Y$ of $\AA$.
If $f : X \to X'$ is a DG morphism in $\AA$ (\ref{def:DG morphism in a DG category}), the \dfn[cone of a morphism in a DG category]{cone of $f$} is an object $\Cone(f)$ for which there are functorial isomorphisms
  \[ \Hom_\AA(Y, \Cone(f)) \isoto \Cone(\Hom_\AA(Y, X) \stackrel{f_*}{\to} \Hom_\AA(Y, X')) \]
in $\C(K)$ for every object $Y$ of $\AA$.
See (\ref{def:translation and cones in category of complexes}) for the definitions of translations and cones in the triangulated category $\C(K)$.

If $\AA$ admits a zero object $0$ and the objects $X[n]$ and $\Cone(f)$ exist for every object $X$, integer $n$, and morphism $f$, then the DG category $\AA$ is called \dfn[pretriangulated DG category]{pretriangulated}.
In this case, the homotopy category $\H^0(\AA)$ is triangulated (\cit{keller2006differential}, section 4.5).

\para[rmk:quasi-equivalence of pretriangulated DG categories]
It is not difficult to show that any DG functor $F : \AA \to \BB$ between pretriangulated categories commutes with translation, i.e. there are functorial isomorphisms $F(X[n]) \isoto F(X)[n]$ in $\BB$ for every object $X$ in $\AA$.
% In fact, one can construct canonical morphisms $F(X)[n] \to F(X[n])$ and $F(X[n]) \to F(X)[n]$ (using the functorial isomorphisms from the definitions of the translations), and do some diagram chasing to show that they are mutual inverses.
Similarly it preserves cones of DG morphisms.
In fact, one sees that the induced functor $\H^0(F) : \H^0(\AA) \to \H^0(\BB)$ is always triangulated.

\begin{paraprop}
Let $\AA$ and $\BB$ be pretriangulated DG categories.
A DG functor $F : \AA \to \BB$ is a quasi-equivalence if and only if the induced functor $\H^0(F) : \H^0(\AA) \to \H^0(\BB)$ is a triangulated equivalence.
\end{paraprop}

The claim obviously reduces to showing that if $\H^0(F)$ is an equivalence, then for each $n \in \Z$ the canonical morphism
  \[ H^n(\Hom_\AA(X, Y)) \longrightarrow H^n(\Hom_{\BB}(F(X), F(Y))) \tag{\ensuremath{\ast}}\label{eq:pretriangulated quasi-equivalence} \]
is an isomorphism.
In fact, by the pretriangulated assumptions one has isomorphisms
  \[ H^n(\Hom_\AA(X, Y)) = H^0(\Hom_\AA(X,Y)[n]) \iso H^0(\Hom_\AA(X[n],Y)) \]
and
  \begin{equation*}
    \begin{split}
      H^n(\Hom_{\BB}(F(X), F(Y))) &= H^0(\Hom_{\BB}(F(X),F(Y))[n])\\
        & \shiso H^0(\Hom_{\BB}(F(X)[n],F(Y))) \\
        & \shiso H^0(\Hom_{\BB}(F(X[n]),F(Y))),
    \end{split}
  \end{equation*}
under which the morphism (\ref{eq:pretriangulated quasi-equivalence}) is identified with the canonical morphism
  \[ H^0(\Hom_\AA(X[n], Y)) \longrightarrow H^0(\Hom_{\BB}(F(X[n]), F(Y))); \]
now this is an isomorphism since $\H^0(F)$ is fully faithful.

%%%%%%%%%%%%%%%%%%%%%%%%%%%%%%%%%%%%%%%%%%%%%%%%%%%%%%%%%%%%%%%%%%%%%%%%%%%

\section{DG categories of perfect complexes}
\label{sec:DG categories of perfect complexes}

\para[def:C(K)-model category associated to scheme]
Let $X$ be a scheme over a commutative ring $K$ and let $\psi : X \to Y = \Spec(K)$ be the structural morphism.
Recall that $\Comp(X)$ denotes the category of complexes of $\O_X$-modules (\ref{rmk:category of O_X-modules}).
It has a natural enrichment (\ref{def:enriched category}) over $\Comp(K)$ which can be described as follows: for two complexes $\sF^\bull$ and $\sG^\bull$, define the complex
  \[ \Hom_{\Cdg(X)}(\sF^\bull, \sG^\bull) = \Gamma(Y, \psi_*(\sHom(\sF^\bull, \sG^\bull))), \]
whose $n$-th component is the $K$-module of global sections of the direct image of the $\O_X$-module $\sHom(\sF^n, \sG^n)$, and differentials are defined in the obvious way.
It is clear that this gives a well-defined category enriched over $\Comp(K)$, i.e. a DG category over $K$, which we denote $\Cdg(X)$.
Further, this enrichment is compatible with the model structure on $\Comp(X)$ (\ref{rmk:model structure on C(Mod(O_X))}), so that one gets a $\Comp(K)$-model category (\ref{def:model category over a symmetric monoidal category}) with the bi-functors $\cdot\otimes\cdot$ and $[\cdot,\cdot]$ defined in the obvious ways.
% probably, one sets
%   (M^\bull \otimes \sF^\bull)^n = \psi^* (M^n)^\tilde \otimes_\O_X \sF^n
% and
%   [M^\bull, \sF^\bull]^n = [\psi^* (M^n)^\tilde, \sF^n]
% i.e. using the canonical functor Mod(A) -> Mod(O_Spec(A)) -> Mod(O_X)
% before applying the tensor product or internal hom of C(Mod(A)).

\para[def:Ddg(X)]
Let $X$ be a scheme over $K$ and consider the $\Comp(K)$-model category $(\Comp(X), \Cdg(X))$ associated to $X$ (\ref{def:C(K)-model category associated to scheme}).
We let $\Ddg(X) = \IntDG(\Comp(X), \Cdg(X))$ denote the associated interior DG category (\ref{def:interior DG category}), so that there is a canonical equivalence of categories
  \begin{equation} \label{eq:Ddg(X) is enrichment of D(X)}
    \H^0(\Ddg(X)) \isoto \D(X)
  \end{equation}
where $\D(X) = \D(\Mod(\O_X))$ denotes the derived category of the category of $\O_X$-modules (\ref{rmk:category of O_X-modules}).

\para[rmk:functoriality of Ddg(X)]
Let $f : X \to Y$ be a morphism of $K$-schemes.
The direct image functor $f_* : \Mod(\O_X) \to \Mod(\O_Y)$ (resp. inverse image functor $f^* : \Mod(\O_Y) \to \Mod(\O_X)$) clearly extends to a functor $f_* : \Cdg(X) \to \Cdg(Y)$ (resp. $f^* : \Cdg(Y) \to \Cdg(X)$).

We define a functor
  \[ f_* : \Ddg(X) \too \Ddg(Y) \]
by composing $f_* : \Cdg(X) \to \Cdg(Y)$ with a cofibrant replacement functor $Q_X$ on $\C(\O_X)$.
This indeed preserves objects which are both cofibrant and fibrant, because $f_*$ already preserves fibrant objects, being a right Quillen functor (\ref{rmk:model structure on C(Mod(O_X))}).
Though this definition depends on the choice of the functor $Q_X$, all cofibrant replacement functors are identified in $\HoDGQe_K$, so one gets a canonical morphism in $\HoDGQe_K$.

In the same way one defines canonical morphisms $f^* : \Ddg(Y) \too \Ddg(X)$ and $\cdot \otimes \cdot : \Ddg(X) \otimes \Ddg(X) \too \Ddg(X)$ in $\HoDGQe_K$.

\para[def:PfDG(X)]
Let $\PfDG(X) \sub \Ddg(X)$ be the full sub-DG category of perfect complexes (\ref{def:perfect complexes}).
It is clear that the equivalence (\ref{eq:Ddg(X) is enrichment of D(X)}) restricts to a canonical equivalence
  \[ \H^0(\PfDG(X)) \isoto \Pf(X). \]

\para[rmk:functoriality of PfDG(X)]
Let $f : X \to Y$ be a morphism of $K$-schemes.
Consider the morphisms $f^* : \Ddg(Y) \to \Ddg(X)$, $f_* : \Ddg(X) \to \Ddg(Y)$, and $\cdot \otimes \cdot : \Ddg(X) \otimes \Ddg(X) \to \Ddg(X)$ in $\HoDGQe_K$ (\ref{rmk:functoriality of PfDG(X)}).
The arguments of (\ref{def:perfect complexes}) show that $f^*$ and $\cdot \otimes \cdot$ preserve perfect complexes and that $f_* : \Ddg(X) \to \Ddg(Y)$ does also when $f$ is proper.
Hence one gets morphisms $f_* : \PfDG(X) \to \PfDG(Y)$, $f_* : \PfDG(Y) \to \PfDG(X)$ and $\cdot \otimes \cdot : \PfDG(X) \otimes \PfDG(X) \to \PfDG(Y)$ in $\HoDGQe_K$.

\para[rmk:Lpf(X) is pretriangulated]
If $X$ is a smooth proper scheme of finite type over a field $K$, then it is straightforward to verify that the DG category $\PfDG(X)$ is pretriangulated (\ref{def:pretriangulated DG category}).
\chapter{Chow motives}
\label{chap:motives}

%%%%%%%%%%%%%%%%%%%%%%%%%%%%%%%%%%%%%%%%%%%%%%%%%%%%%%%%
\section{Chow groups}
\label{sec:Chow groups}
%%%%%%%%%%%%%%%%%%%%%%%%%%%%%%%%%%%%%%%%%%%%%%%%%%%%%%%%

\para[def:swap morphism]
If $X$ and $Y$ are schemes, we will write $\sigma_{X,Y} : X \times Y \isoto Y \times X$ for the morphism that swaps the factors.

\para[def:algebraic cycles]
Let $X$ be a noetherian scheme.
The free abelian group generated by the closed integral subschemes of codimension $k$ is denoted $Z^k(X)$, and its elements are called \dfn[]{$k$-codimensional cycles on $X$}.
%Dually we write $\Z^k(X)$ for the \dfn[]{group of $k$-cocycles on $X$}, the free abelian group generated by the closed integral subschemes of codimension $k$.
An element of $Z^k(X)$ is called a \dfn[]{$k$-codimensional cycle on $X$} and is written
  \[ \alpha = \sum_{x \in X^{(k)}} n_x . \bar{\{x\}} \]
for some $n_x \in \Z$, where $X^{(k)} \sub X$ denotes the subset of $k$-codimensional points.
(By definition, only finitely many of the $n_x$ are nonzero.)

\para[def:Weil divisor]
Recall that a \dfn[]{Weil divisor on $X$} is by definition a 1-codimensional cycle.
Let $R(X)$ denote the ring of rational functions on $X$; recall that when $X$ is integral, there are canonical isomorphisms $R(X) \isoto \Frac(\O_{X,x})$ for every point $x \in X$ (\citcust{EGAI}{EGA I}, 7.1.5).
The \dfn[]{order of vanishing of a rational function $r$ at a point $x$} is defined as
  \[ \ord_x(r) = \length_{\O_{X,x}}(\O_{X,x}/(a_x)) - \length_{\O_{X,x}}(\O_{X,x}/(b_x)) \]
where $a_x/b_x \in \Frac(\O_{X,x})$ is the fraction corresponding to $r$, and $\length_A(M)$ denotes the length of an $A$-module $M$.
The \dfn[]{Weil divisor associated to $r$} is the 1-codimensional cycle
  \[ \div(r) = \sum_{x \in X^{(1)}} \ord_x(r) . \bar{\{x\}}, \]
where the sum is taken over the 1-codimensional points of $X$.
One verifies that this is well-defined because there are only finitely many points $x \in X^{(1)}$ with $\ord_x(r) \ne 0$.
See (\citcust{EGAIV4}{EGA IV\textsubscript{4}}, 21.6).

\para[def:Chow groups]
Two $k$-codimensional cycles $\alpha$ and $\beta$ on $X$ are called \dfn[rational equivalence of cycles]{rationally equivalent} if there exists a family $(Z_i)_i$ of $(k-1)$-codimensional closed integral subschemes and invertible rational functions $r_i \in R(Z_i)^*$ such that
  \[ \alpha - \beta = \sum_i (j_i)_* (\div(r_i)), \]
where $j_i$ denotes the inclusion morphism $Z_i \hook X$.
The \dfn[Chow groups]{Chow group of codimension $k$ of $X$}, denoted $\CH^k(X)$, is the quotient of $Z^k(X)$ by the subgroup of $k$-cocycles rationally equivalent to the zero cycle.
For a commutative ring $\Lambda$, we let $\CH^k(X, \Lambda)$ denote the $\Lambda$-module $\CH^k(X) \otimes_\Z \Lambda$, called the \dfn[]{Chow group of codimension $k$ with coefficients in $\Lambda$}.
%Dually we define $\CH_k(X)$ (resp. $\CH_k(X, \Lambda)$), the \dfn[]{Chow group of dimension $k$ (with coefficients in $\Lambda$)}.

\para[rmk:functoriality of Chow groups]
Let $X$ and $Y$ be smooth proper schemes of finite type over a field $K$.
Any morphism $f : X \to Y$ induces functorially a degree zero homomorphism $f^* : \CH^*(Y, \Lambda) \to \CH^*(X, \Lambda)$ of graded $\Lambda$-modules.
If $X$ and $Y$ are purely of dimensions $m$ and $n$, respectively, then $f$ also induces functorially a degree $n-m$ homomorphism $f_* : \CH^*(X, \Lambda) \to \CH^*(Y, \Lambda)$ of graded $\Lambda$-modules.

\para[def:transpose of a cycle]
Let $X$ and $Y$ be smooth proper schemes of finite type over $K$.
Recall that $\sigma_{X,Y} : X \times Y \isoto Y \times X$ denotes the swap morphism (\ref{def:swap morphism}).
The \dfn[]{transpose of a cycle} $\alpha$ on $X \times Y$ is defined as the cycle $\sigma_{Y,X}^*(\alpha)$ on $Y \times X$.

\para[def:cartesian product homomorphism]
There is a canonical degree zero homomorphism
  \[ \CH^*(X, \Lambda) \otimes_\Lambda \CH^*(Y, \Lambda) \longrightarrow \CH^*(X \times Y, \Lambda) \]
of graded $\Lambda$-modules, which is functorial with respect to inverse and direct images.
The image of $\alpha \otimes \beta$, denoted $\alpha \times \beta$, is called the \dfn[]{cartesian product of $\alpha$ and $\beta$}.
It is associative and commutative.

\para[def:degree homomorphism]
Suppose $X$ is purely of dimension $n$.
The \dfn[]{degree homomorphism}
  \[ \deg_X : \CH^n(X, \Lambda) \longrightarrow \Lambda \]
is defined as the composite
  \[ \CH^*(X, \Lambda) \stackrel{\varphi_*}{\longrightarrow} \CH^*(\Spec(K), \Lambda) \isoto \Lambda \]
where $\varphi : X \to \Spec(K)$ is the structural morphism.

%\para[def:augmentation morphism]
%For $X$ integral, the canonical \dfn[]{augmentation morphism}
%  \[ \epsilon : \CH^*(X, \Lambda) \longrightarrow \Lambda, \]
%which maps the cycle $[X]$ to the unit of $\Lambda$ and all $k$-codimensional cycles to zero $(k > 0)$, induces an isomorphism $\CH^0(X, \Lambda) \isoto \Lambda$.
%It is functorial with respect to inverse image and compatible with the cartesian product of cycles (\ref{def:cartesian product homomorphism}).

\para[def:intersection product]
Let $X$ be a smooth proper scheme of finite type over $K$.
For cycles $\alpha$ and $\beta$ of codimension $j$ and $k$, respectively, we define the \dfn[]{intersection product} $\alpha \inpr \beta$ by
  \[ \alpha \inpr \beta = \Delta_X^*(\alpha \times \beta) \in \CH^{j+k}(X, \Lambda) \]
where $\Delta_X : X \to X \times X$ denotes the diagonal morphism of $X$.
This gives $\CH^*(X)$ the structure of a graded commutative $\Lambda$-algebra.
For a morphism $f : X \to Y$, the homomorphism $f^* : \CH^*(Y) \to \CH^*(X)$ becomes a homomorphism of graded $\Lambda$-algebras.

\para[rmk:projection formula]
For a morphism $f : X \to Y$ and cycles $\alpha$ and $\beta$ on $X$ and $Y$, respectively, one has the projection formula
  \[ f_*(f^*(\beta) \inpr \alpha) = \beta \inpr f_*(\alpha). \]

\para[rmk:Weil divisor associated to invertible sheaf]
Let $X$ be an \emph{integral} smooth proper scheme of finite type over $K$.
Let $\Pic(X)$ denote the \dfn[]{Picard group of $X$}, the group of isomorphism classes of invertible $\O_X$-modules; the class of an invertible $\O_X$-module $\sL$ will be written $\cl_X(\sL)$.
There is a canonical isomorphism
  \[ p_X : \Pic(X) \isoto \CH^1(X) \]
defined as follows.
Given an invertible $\O_X$-module $\sL$, let $s \in \Gamma(X, \M_X(\sL))$ be a regular meromorphic section, where $\M_X(\sL)$ denotes the sheaf of meromorphic sections of $\sL$ on $X$ (\citcust{EGAIV4}{EGA IV\textsubscript{4}}, 20.1).
Via the canonical isomorphism of stalks $(\M_X(\sL))_x \iso \M_{X,x} \otimes_{\O_{X,x}} \sL_x$, the germ of $s$ at a point $x \in X$ corresponds to an element $f_x \otimes t_x$ with $f_x \in \M_{X,x}$ and $t_x \in \sL_x$ a generator (unique up to multiplication by an invertible element of $\O_{X,x}$).
Since $X$ is integral, the element $f_x \in \M_{X,x}$ corresponds canonically to an element $a_x/b_x$ of the fraction field $\Frac(\O_{X,x})$;
since $s$ is \emph{regular}, the fraction $a_x/b_x$ is further nonzero.
We define for any point $x \in X^{(1)}$ of codimension 1 the \dfn[]{order of vanishing of $s$ at the point $x$} by
  \[ \ord_x(s) = \ord_{\O_{X,x}}(f_x) = \length_{\O_{X,x}}(a_x) - \length_{\O_{X,x}}(b_x). \]
Then we define the image of $\cl_X(\sL)$ by $p_X$ as the class of the cycle
  \[ \div_\sL(s) = \sum_{x \in X^{(1)}} \ord_x(s) . \bar{\{x\}} \in Z^1(X), \]
which one checks is independent of $s$.

%%%%%%%%%%%%%%%%%%%%%%%%%%%%%%%%%%%%%%%%%%%%%%%%%%%%%%%%%%%%%%%
\section{Grothendieck groups}
%%%%%%%%%%%%%%%%%%%%%%%%%%%%%%%%%%%%%%%%%%%%%%%%%%%%%%%%%%%%%%%

\para[def:Grothendieck groups of a scheme]
Let $X$ be a scheme.
The category of locally free $\O_X$-modules of finite rank, being an additive thick full subcategory of $\Mod(\O_X)$, is exact, and we let
  \[ K^0(X) = K_0(\LocFr(\O_X)) \]
denote its Grothendieck group (\ref{def:Grothendieck group of an exact category}).
When $X$ is noetherian, the category $\Coh(\O_X)$ of coherent $\O_X$-modules is abelian and we let
  \[ K_0(X) = K_0(\Coh(\O_X)) \]
denote its Grothendieck group.

\para[rmk:ring and module structures on Grothendieck groups of X]
The tensor product $\otimes_{\O_X}$ of $\O_X$-modules induces the structure of a commutative ring on $K^0(X)$ (\citcust{SGA6}{SGA 6}, Exp IV, 2.7, b).
Similarly one gets a multiplication map
  \[ K^0(X) \otimes_\Z K_0(X) \longrightarrow K_0(X) \]
which gives $K_0(X)$ the structure of a $K^0(X)$-module.
See (\citcust{SGA6}{SGA 6}, Exp. IV, 2.10).

\para[def:functoriality of Grothendieck groups]
A morphism of schemes $f : X \to Y$ induces functorially a homomorphism of rings $f^* : K^0(Y) \to K^0(X)$ (\citcust{SGA6}{SGA 6}, Exp. IV, 2.7, b).

If $f$ is proper and $Y$ is noetherian, then it also induces a homomorphism of $K^0(Y)$-modules
  \[ f_* : K_0(X) \to K_0(Y), \]
where the $K^0(X)$-module $K_0(X)$ is viewed as a $K^0(Y)$-module via the homomorphism $f^* : K^0(Y) \to K^0(X)$.

\para[rmk:isomorphism of Grothendieck groups of a scheme]
If $X$ is a smooth separated noetherian scheme, then there is a canonical isomorphism of $K^0(X)$-modules
  \[ K^0(X) \isoto K_0(X) \label{eq:isomorphism of Grothendieck groups}. \]
See (\citcust{SGA6}{SGA 6}, Exp. IV, 2.5).

\para[def:Euler characteristic of perfect complex]
Let $X$ be a scheme.
Recall that $\Pf(X)$ is the triangulated category of perfect complexes on $X$ (\ref{def:perfect complexes}).
Given a perfect complex $\sF^\bull \in \Pf(X)$, we define its \dfn[]{Euler characteristic} $\chi_X(\sF^\bull) \in K_0(X)$ as the alternating sum of its cohomologies,
  \[ \chi_X(\sF^\bull) = \sum_{i \in \Z} (-1)^i \cdot [H^i(\sF^\bull)]; \]
this is well-defined because perfect complexes are bounded with coherent cohomology (\ref{def:perfect complexes}).

\para[rmk:Grothendieck groups of X and Pf(X) are isomorphic]
Let $X$ be a scheme.
Consider the Grothendieck group of the triangulated category $\Pf(X)$ (\ref{def:Grothendieck group of a triangulated category}).
Since the Euler characteristic (\ref{def:Euler characteristic of perfect complex}) is independent of the class in $K_0(\Pf(X))$, one gets a canonical homomorphism
  \[ \chi = \chi_X : K_0(\Pf(X)) \too K_0(X). \]

\begin{paraprop}
The homomorphism $\chi : K_0(\Pf(X)) \to K_0(X)$ defined by the Euler characteristic is bijective.
\end{paraprop}

In fact, its inverse is given by the homomorphism that maps the class of a coherent sheaf $\sF$ to the class of the complex $\sF^\bull$ concentrated in degree zero.

%%%%%%%%%%%%%%%%%%%%%%%%%%%%%%%%%%%%%%%%%%%%%%%%%%%%%%%%
\section{Characteristic classes}
\label{sec:characteristic classes}
%%%%%%%%%%%%%%%%%%%%%%%%%%%%%%%%%%%%%%%%%%%%%%%%%%%%%%%%

\para
Let $X$ be a scheme of finite type over a field $K$.
Recall that there is an additive, symmetric monoidal equivalence between the category of locally free $\O_X$-modules and the category of vector bundles on $X$, given by the assignment
  \[ \sE \mapsto \Spec(\bS(\sE)) \]
mapping a locally free $\O_X$-module $\sE$ to the affine spectrum of its symmetric $\O_X$-algebra (\citcust{EGAII}{EGA II}, 1.7).

\para[rmk:generators of Chow ring of projective bundle]
Let $\sE$ be a locally free $\O_X$-module of rank $p$.
Let $P = \P(\sE) = \bProj(\bS(\sE))$ be the projective bundle associated to $\sE$ and recall that there is a canonical invertible $\O_{\P(\sE)}$-module $\O_{\P(\sE)}(1)$ (\citcust{EGAII}{EGA II}, 4.1.1).

Let $f : P \to X$ be the projection morphism.
Recall that we have an isomorphism $p_P : \Pic(P) \isoto \CH^1(P)$ (\ref{rmk:Weil divisor associated to invertible sheaf}).
Let $\xi_\sE$ be the element of $\CH^1(P)$ corresponding to the class of $\O_P(1)$, i.e.
  \[ \xi_\sE = p_P(\cl_P(\O_{P}(1))). \]
The elements $(\xi_\sE)^i \in \CH^i(P)$ $(0 \le i \le p-1)$ generate $\CH^*(P)$, viewed as an $\CH^*(X)$-module via the homomorphism $f^* : \CH^*(X) \to \CH^*(P)$.
See (\cit{grothendieck1958intersections}, Proposition 4).

\para[def:Chern classes]
Considering the class $(\xi_\sE)^p \in \CH^p(P)$, one has by (\ref{rmk:generators of Chow ring of projective bundle}) unique classes $c_i(\sE) \in \CH^i(X)$ such that
  \[ (\xi_\sE)^p + c_1(\sE) (\xi_\sE)^{p-1} + \cdots + c_p(\sE) = 0 \]
and $c_0(\sE) = 1$, $c_i(\sE) = 0$ for $i > p$.
We call the class $c_i(\sE) \in \CH^i(X)$ the \dfn[Chern classes]{$i$-th Chern class of $\sE$} $(1 \le i \le p)$, and the sum
  \[ c(\sE) = \sum_{i \ge 0} c_i(\sE) \in \CH^*(X) \]
is called the \dfn[]{(total) Chern class of $\sE$}.

\para[thm:properties of Chern classes]
Chern classes are characterized completely by the following properties.
\begin{enumerate}
  \item For a morphism $f : X \to Y$ and a locally free $\O_Y$-module $\sE$, one has
    \[ c_i(f^*(\sE)) = f^*(c_i(\sE)). \]
  \item If $\sL$ is an invertible $\O_X$-module, then
    \[ c_1(\sL) = p_X(\cl_X(\sL)). \]
  \item If there is an exact sequence
    \[ 0 \to \sE' \to \sE \to \sE'' \to 0 \]
  of locally free $\O_X$-modules, then one has an equality
    \[ c(\sE) = c(\sE') \inpr c(\sE''). \]
\end{enumerate}
See (\cit{grothendieck1958intersections}, Th\'eor\`eme 1).

\para[def:flag variety]
Let $\sE$ be a locally free $\O_X$-module of rank $p$.
Consider the locally free sheaf $\sE^{(1)} = f^*(\sE) / \O_{\P(\sE)}(-1)$ of rank $p-1$ on $X^{(1)} = \P(\sE)$.
On $X^{(2)} = \P(\sE^{(1)})$ one has the locally free sheaf $\sE^{(2)} = (\sE^{(1)})^{(1)}$ of rank $p-2$.
Iterating this process one gets a sequence of locally free sheaves $\sE^{(i)}$ of rank $p-i$ on $X^{(i)}$ $(1 \le i \le p)$.
We call $\Flag(\sE) = X^{(p)}$ the \dfn[flag variety associated to a locally free sheaf]{flag variety associated to $\sE$}.
Let $f : \Flag(\sE) \to X$ be the projection morphism.
The inverse image $f^*(\sE)$ \dfn[]{splits completely}, i.e. there exists an increasing sequence of locally free sheaves
  \[ 0 = \sE_0 \sub \sE_1 \sub \cdots \sub \sE_p = f^*(\sE) \]
such that $\sE_i$ has rank $i$, and the successive quotients $\sE_i / \sE_{i-1}$ are of rank 1 (i.e. invertible sheaves).
See (\cit{grothendieck1958intersections}, \S 3).

\para[def:Chern roots]
Let $\sE$ be a locally free $\O_X$-module of rank $p$.
Let $f : \Flag(\sE) \to X$ be the flag variety associated to $\sE$.
By the properties of Chern classes (\ref{thm:properties of Chern classes}) one gets the formula
  \[ f^*(c(\sE)) = c(f^*(\sE)) = \prod_{i=1}^p (1 + p_{\Flag(\sE)}(\cl_{\Flag(\sE)}(\sE_i / \sE_{i-1}))). \]
Here $p_{\Flag(\sE)}(\cl_{\Flag(\sE)}(\sE_i / \sE_{i-1})) \in \CH^1(\Flag(\sE))$ denotes the 1-codimensional cycle on $\Flag(\sE)$ corresponding to the class of $\sE_i / \sE_{i-1}$ in $\Pic(\Flag(\sE))$ (\ref{rmk:Weil divisor associated to invertible sheaf}).

The homomorphism $f^* : \CH^*(X) \to \CH^*(\Flag(\sE))$ is injective, and the preimages $\alpha_i \in \CH^*(X)$ of the classes $p_{\Flag(\sE)}([\sE_i/\sE_{i-1}])$ on $\Flag(\sE)$ are called the \dfn[]{Chern roots of $\sE$}.

\para[def:Chern character]
Let $\sE$ be a locally free $\O_X$-module of rank $p$ and let $\alpha_i$ $(1 \le i \le p)$ be its Chern roots (\ref{def:Chern roots}).
The \dfn[Chern character]{Chern character of $\sE$} is the class $\ch(\sE) \in \CH^*(X, \Q)$ defined by the formula
  \[ \ch(\sE) = \sum_{i=1}^p \exp(\alpha_i) = \sum_{i=1}^p \left( \sum_{k=0}^\infty \frac{1}{k!} (\alpha_i)^k \right). \]
In terms of the Chern classes one computes
  \begin{equation*}
    \begin{split}
      \ch(\sE) = p
        & + c_1(\sE) + \frac{1}{2}(c_1(\sE)^2 - 2c_2(\sE)) \\
        & + \frac{1}{6}(c_1(\sE)^3 - 3c_1(\sE)c_2(\sE) + c_3(\sE)) \\
        & + \frac{1}{24}(c_1(\sE)^4 - 4c_1(\sE)^2 c_2(\sE) + 4c_1(\sE)c_3(\sE) + 2c_2(\sE)^2 - 4c_4(\sE))
          + \cdots
    \end{split}
  \end{equation*}
The \dfn[Todd class]{Todd class of $\sE$} is the class $\td(\sE) \in \CH^*(X, \Q)$ defined by
  \[ \td(\sE) = \sum_{i=1}^p \frac{\alpha_i}{1-e^{-\alpha_i}}. \]
Explicitly one computes
  \begin{equation*}
    \begin{split}
      \td(\sE) = 1
        & + \frac{1}{2}c_1(\sE) + \frac{1}{12}(c_1(\sE)^2 + c_2(\sE))
          + \frac{1}{24}c_1(\sE)c_2(\sE)
        \\
        & + \frac{1}{720}(-c_1(\sE)^4 + 4c_1(\sE)^2c_2(\sE) + 3c_2(\sE)^2 + c_1(\sE)c_3(\sE) - c_4(\sE))
          + \cdots
    \end{split}
  \end{equation*}

\para[rmk:functoriality of Chern character wrt inverse image]
Note that the Chern character and Todd class are functorial with respect to inverse image.
Indeed, for a morphism $f : X \to Y$, since $f^* : \CH^*(Y) \to \CH^*(X)$ is a homomorphism of graded rings, this follows from (\ref{thm:properties of Chern classes}, (i)).

\para[rmk:additivity of characteristic classes]
One verifies without difficulty that if there is an exact sequence
  \[ 0 \to \sE' \to \sE \to \sE'' \to 0 \]
of locally free sheaves, then one has the formulas
  \begin{equation} \label{eq:additivity of Chern character}
    \ch(\sE) = \ch(\sE') + \ch(\sE'')
  \end{equation}
and
  \begin{equation} \label{eq:additivity of Todd class}
    \td(\sE) = \td(\sE') \inpr \td(\sE'')
  \end{equation}

\para[rmk:Chern and Todd homomorphisms]
Let $X$ be a smooth separated scheme of finite type over a field $K$.
Recall that $K^0(X)$ and $K_0(X)$ denote the Grothendieck groups of locally free and coherent sheaves, respectively (\ref{def:Grothendieck groups of a scheme}).
By (\ref{eq:additivity of Chern character}) and (\ref{eq:additivity of Todd class}), the Chern character and Todd class induce homomorphisms of rings $K^0(X) \to \CH^*(X)$.
Via the isomorphism $K^0(X) \shiso K_0(X)$ (\ref{eq:isomorphism of Grothendieck groups}), the Chern character and Todd class induce homomorphisms
  \[ \ch, \td : K_0(X) \to \CH^*(X) \]
of rings.

\para[def:Chern character of a complex]
Let $X$ be a smooth separated scheme of finite type over a field $K$.
Recall that $\Pf(X)$ denotes the triangulated category of perfect complexes on $X$ (\ref{def:perfect complexes}).
Let $\sF^\bull \in \Pf(X)$ be a perfect complex and consider its Euler characteristic (\ref{def:Euler characteristic of perfect complex})
  \[ \chi(\sF^\bull) = \sum_{i \in \Z} (-1)^i [H^i(\sF^\bull)] \in K_0(X) \]
in the Grothendieck group of $X$ (\ref{def:Grothendieck groups of a scheme}).
We define the \dfn[]{Chern character $\ch(\sF^\bull)$} (resp. \dfn[]{Todd class $\td(\sF^\bull)$}) as the Chern character (resp. Todd class) of $\chi(\sF^\bull)$ (\ref{rmk:Chern and Todd homomorphisms}).

\para
Let $X$ be a smooth separated scheme of finite type over a field $K$.
We write $\Omega_{X/K}$ for the \dfn[]{sheaf of differentials on $X$ relative to $K$} and
  \[ \sT_{X/K} = \sHom_{\O_X}(\Omega_{X/K}, \O_X) \]
for its dual, \dfn[]{the tangent bundle on $X$ relative to $K$}, or just $\sT_{X}$ when there is no risk of confusion.
The latter is a coherent $\O_X$-module (\citcust{EGAIV4}{EGA IV\textsubscript{4}}, 16.5.7) and we write $\td_X = \td(\sT_{X/K})$ for its Todd class (\ref{rmk:Chern and Todd homomorphisms}).
We write
  \[ \sqrt{\td_{X}} = \exp\left(\frac{1}{2} \log(\td_{X})\right) \]
so that $\sqrt{\td_X} \cdot \sqrt{\td_X} = \td_X$.

\para[thm:GRR]
\begin{parathm}[Grothendieck-Riemann-Roch]
Let $X$ and $Y$ be smooth projective schemes of finite type over a field $K$.
For a morphism $f : X \to Y$, the diagram
  \[ \begin{tikzcd}[column sep=large]
    K_0(X) \arrow{r}{\ch(\cdot) \td_X}\arrow{d}{f_*}
      & \CH^*(X) \arrow{d}{f_*}
    \\
    K_0(Y) \arrow{r}{\ch(\cdot) \td_Y}
      & \CH^*(Y)
  \end{tikzcd} \]
commutes.
\end{parathm}

See (\citcust{SGA6}{SGA 6}, VI, Exp 0).

\para[thm:Chern character is an isomorphism]
\begin{parathm}
Let $X$ be a smooth projective scheme of finite type over a field $K$.
The Chern character induces a canonical isomorphism
  \[ K_0(X) \otimes \Q \isoto \CH^*(X, \Q). \]
\end{parathm}

This follows from (\ref{thm:GRR}), see (\cit{fulton1998intersection}, 15.2.16).

\para[rmk:Chern character of a coherent sheaf supported on subscheme]
Let $X$ be a smooth projective scheme of finite type over a field $K$.
Let $\sF$ be a coherent $\O_X$-module and write $\Supp(\sF) \sub X$ for its support, considered as a topological space.
Recall that there exists a closed subscheme $Z \sub X$ with underlying topological space $\Supp(\sF)$ such that there is a canonical isomorphism $\sF \isoto j_*(j^*(\sF))$, where $j : Z \hook X$ is the inclusion morphism (\citcust{EGAI}{EGA I}, 9.3.5).
Now by (\ref{thm:GRR}) the cycle $\ch(\sF) = \ch(j_*(j^*(\sF)))$ lies in the image of the homomorphism $j_* : A_*(Z) \to A_*(X)$; this maps cycles of codimension $i$ to cycles of codimension $n-m+i$, where $m = \dim(Z)$ and $n = \dim(X)$.
In particular one has

\begin{paraprop}
Let $X$ be an $n$-dimensional smooth projective scheme of finite type over $K$.
For a coherent sheaf $\sF$ on $X$ whose support has dimension $m$, the Chern character $\ch(\sF)$ has no components of codimension $0, \ldots, n - m$.
\end{paraprop}

Now suppose we have a perfect complex of sheaves.
Recall that the support of a complex is defined to be the union of the supports of its cohomology sheaves.
It is straightforward to see, arguing as above, that one also has

\begin{paraprop}
Let $X$ be an $n$-dimensional smooth projective scheme of finite type over $K$.
For a perfect complex $\sF^\bull$ on $X$ whose support has dimension $m$, the Chern character $\ch(\sF^\bull)$ has no components of codimension $0, \ldots, n - m$.
\end{paraprop}

%%%%%%%%%%%%%%%%%%%%%%%%%%%%%%%%%%%%%%%%%%%%%%%%%%%%%%%%
\section{Chow motives}
\label{sec:chow motives}
%%%%%%%%%%%%%%%%%%%%%%%%%%%%%%%%%%%%%%%%%%%%%%%%%%%%%%%%

\para[rmk:symmetric monoidal structure on SmProj]
Let $\Var_K$ denote the category of smooth projective varieties over a field $K$.
The fibered product over $K$ induces a symmetric monoidal structure on $\Var_K$.

\para[def:composition of cycles]
Let $X, Y, Z \in \Var_K$ of dimension $m$, $n$ and $p$, respectively.
Given homogeneous cycles $\alpha \in \CH^i(X \times Y)$ and $\beta \in \CH^j(Y \times Z)$, we define their composition as follows.
Let $p_{XY}$, $p_{XZ}$ and $p_{YZ}$ be the three projections from $X \times Y \times Z$ onto $X \times Y$, $X \times Z$ and $Y \times Z$, respectively.
Via the inverse image homomorphisms $p_{XY}^* : \CH^i(X \times Y) \to \CH^i(X \times Y \times Z)$ and $p_{YZ}^* : \CH^j(Y \times Z) \to \CH^j(X \times Y \times Z)$, one gets cycles on $X \times Y \times Z$.
Their intersection product $p_{XY}^*(\alpha) \inpr p_{YZ}^*(\beta)$ lies in $\CH^{i+j}(X \times Y \times Z)$.
Then we define the \dfn[]{composition of the homogeneous cycles $\alpha$ and $\beta$} as the cycle
  \[ \beta \circ \alpha = (p_{XZ})_*(p_{XY}^*(\alpha) \inpr p_{YZ}^*(\beta)) \]
which has codimension $i+j-n$ and lies in $\CH^{i+j-n}(X \times Z)$.
Note that the diagonal $[\Delta_X(X)] \in \CH^m(X \times X)$ is the identity of $X$ with respect to this composition law.

For mixed cycles $\alpha \in \CH^*(X \times Y)$ and $\beta \in \CH^*(Y \times Z)$ we define the composition as the cycle whose $k$-th component is given by the sum
  \[ (\beta \circ \alpha)^k = \sum_{i+j = k+n} \beta^j \circ \alpha^i \]
where the indices $i$ and $j$ range over the nonnegative integers not greater than $m$ and $n$, respectively, whose sum is equal to $k+n$.
Note that the identity of $X$ is still given by the cycle $[\Delta_X(X)] \in \CH^*(X \times X)$.

\para[def:Chow correspondences]
Let $X, Y \in \Var_K$ of dimension $m$ and $n$, respectively.
A \dfn[Chow correspondence]{Chow correspondence between $X$ and $Y$ of degree $d$} is an element of the direct sum
  \[ \CorrGp^d_K(X, Y) = \bigoplus_i \CH^{m_i+d}(X_i \times Y) \]
where $(X_i)_i$ is the family of irreducible components of $X$ and $m_i = \dim(X_i)$.
Using the law of composition described in (\ref{def:composition of cycles}), we may define a category $\Corr_K$ whose objects are smooth projective varieties over $K$ and morphisms are Chow correspondences of degree zero.
This is called the \dfn[]{category of correspondences over $K$}.
For a commutative ring $\Lambda$, we also define the \dfn[]{category $\Corr_K(\Lambda)$ of Chow correspondences over $K$ with coefficients in $\Lambda$}, where the $\Lambda$-modules of morphisms are
  \[ \Hom_{\Corr_K(\Lambda)}(X, Y) = \CorrGp^0_K(X, Y; \Lambda) = \bigoplus_i \CH^{m_i}(X_i \times Y, \Lambda). \]

\para
There is a canonical functor
  \[ M_{0,\Lambda} : \Var_K^\circ \longrightarrow \Corr_K(\Lambda) \]
mapping a morphism $f : X \to Y$ to the class of $\sigma_{Y,X}^*[\Gamma_f(X)]$, the transpose (\ref{def:transpose of a cycle}) of its graph.
By abuse of notation we will write $M_0 = M_{0,\Lambda}$ when there is no risk of confusion.

\para[rmk:additive and symmetric monoidal structure on Corr]
The category $\Corr_K(\Lambda)$ is $\Lambda$-linear, with direct sums given by
  \[ M_0(X) \oplus M_0(Y) = M_0(X \sqcup Y). \]

The product in $\Var_K$ also induces a symmetric monoidal structure on $\Corr_K(\Lambda)$: we define $M_0(X) \otimes M_0(Y) = M_0(X \times Y)$, and we define homomorphisms
  \[ \Hom(M_0(X), M_0(X')) \otimes \Hom(M_0(Y), M_0(Y'))
    \longrightarrow \Hom(M_0(X) \otimes M_0(Y), M_0(X') \otimes M_0(Y')) \]
by mapping $\alpha \otimes \beta$ to the morphism corresponding to $p^*(\alpha) \inpr q^*(\beta)$ via the canonical isomorphism
  \[ \CH^{m'+n'}(X \times X' \times Y \times Y', \Lambda) \iso \Hom_{\Corr_K(\Lambda)}(M_0(X \times Y), M_0(X' \times Y')), \]
where $\alpha \in \Hom_{\Corr_K(\Lambda)}(M_0(X), M_0(X'))$, $\beta \in \Hom_{\Corr_K(A)}(M_0(Y), M_0(Y'))$.
Here $p$ and $q$ denote the projection morphisms from $X \times X' \times Y \times Y'$ to $X \times X'$ and $Y \times Y'$, respectively, and $m'$ and $n'$ are the dimensions of $X'$ and $Y'$, respectively.
By abuse of notation we denote $p^*(\alpha) \inpr q^*(\beta)$ again by $\alpha \otimes \beta$.
The functor $M_0$ is clearly compatible with the symmetric monoidal structures.

\para[rmk:M_0(P^1)]
Consider the image $M_0(\P^1_K)$ of the projective line in $\Corr_K(\Lambda)$.
One can prove that the $\Lambda$-module of its endomorphisms, i.e. $\CH^1(\P^1_K \times \P^1_K, \Lambda)$, is isomorphic to $\Lambda \oplus \Lambda$, and is generated by the classes $\alpha = [\{\infty\} \times \P^1_K]$ and $\beta = [\P^1_K \times \{\infty\}]$.
The structural morphism $\P^1_K \to \Spec(K)$ induces a morphism $M_0(\Spec(K)) \to M_0(\P^1)$, and the inclusion morphism $\Spec(K) \to \P^1_K$ of the point $\infty$ induces a morphism $M_0(\P^1) \to M_0(\Spec(K))$; one sees that the composite $M_0(\Spec(K)) \to M_0(\P^1) \to M_0(\Spec(K))$ is the identity, and the other composite is $\alpha$.
See (\cit{andre2004introduction}, 4.1.2.1).

\para[def:category of effective Chow motives]
We define the \dfn[effective Chow motives]{category $\CHM^+_K(\Lambda)$ of effective Chow motives over $K$ with coefficients $\Lambda$} as the karoubian envelope (\ref{def:karoubian envelope}) of the category $\Corr_K(\Lambda)$.
An object of $\CHM^+_K(\Lambda)$ is a pair $(X, \alpha)$, where $X \in \Var_K$ and $\alpha \in \CorrGp^0(X \times X, \Lambda)$ with $\alpha \circ \alpha = \alpha$.
Morphisms $(X, \alpha) \to (Y, \beta)$ are of the form
  \[ \beta \circ \gamma \circ \alpha \in \CH^n(X \times Y, \Lambda) \]
for some $\gamma \in \CorrGp^0(X, Y; \Lambda)$.

\para[rmk:CHM+ is additive]
The category $\CHM^+_K(\Lambda)$ is $\Lambda$-linear, with direct sums defined by
  \[ (X, \alpha) \oplus (Y, \beta) = (X \sqcup Y, (\alpha, \beta)) \]
where we abuse notation and write $(\alpha, \beta)$ for the class in $\CH^*(X \sqcup Y, \Lambda)$ corresponding to $(\alpha, \beta) \in \CH^*(X, \Lambda) \oplus \CH^*(Y, \Lambda)$ under the canonical isomorphism.

\para[rmk:CHM+ is symmetric monoidal]
The category $\CHM^+_K(\Lambda)$ is symmetric monoidal, with the tensor product of two objects $(X, \alpha)$ and $(Y, \beta)$ given by
  \[ (X, \alpha) \otimes (Y, \beta) = (X \times Y, \alpha \otimes \beta) \]
where $\alpha \otimes \beta$ is defined in (\ref{rmk:additive and symmetric monoidal structure on Corr}).
The unit object is $\IdM^+ = (\Spec(K), [\Delta])$ where $[\Delta] \in \CH^1(\Spec(K) \times \Spec(K), \Lambda)$ denotes the class of the diagonal.

\para[def:effective motive functor]
There is a canonical functor
  \[ M^+_\Lambda : \Var_K^\circ \longrightarrow \CHM^+_K(\Lambda) \]
which is the composition of $M_{0,\Lambda} : \Var_K^\circ \to \Corr_K(\Lambda)$ with the canonical fully faithful functor $\Corr_K(\Lambda) \hook \CHM^+_K(\Lambda)$ (\ref{def:karoubian envelope}).
By abuse of notation we write $M^+ = M^+_\Lambda$ when there is no risk of confusion.
One verifies without difficulty that the functor $M^+$ is additive and symmetric monoidal (i.e. compatibile with the structures defined in (\ref{rmk:CHM+ is additive}) and (\ref{rmk:CHM+ is symmetric monoidal})).

\para[def:Lefschetz effective motive]
Let $\LefM^+$ be the object $(\P^1, \beta)$ in $\CHM^+_K(\Lambda)$, where $\beta = [\P^1 \times \{\infty\}] \in \CH^1(\P^1 \times \P^1, \Lambda)$ is the class of the closed subscheme $\P^1 \times \{\infty\}$.
By the discussion in (\ref{rmk:M_0(P^1)}), one sees that there is a canonical decomposition
  \begin{equation}
    M^+(\P^1) \iso \IdM^+ \oplus \LefM^+ \label{eq:Lefschetz decomposition}
  \end{equation}
of the image of the projective line.

%The functor $\cdot \otimes \LefM^+ : \CHM^+_K(\Lambda) \to \CHM^+_K(\Lambda)$ is fully faithful;
%in other words, for any objects $M$ and $N$ of $\CHM^+_K(\Lambda)$, the canonical morphism
%  \[ \Hom_{\CHM^+_K(\Lambda)}(M, N) \longrightarrow \Hom_{\CHM^+_K(\Lambda)}(M \otimes \LefM^+, N \otimes \LefM^+) \]
%is an isomorphism of $\Lambda$-modules (Riou, 5.8).

\para[def:category of Chow motives]
Let $\CHM_K(\Lambda)$ be the karoubian envelope (\ref{def:karoubian envelope}) of the category whose objects are pairs $(X, r)$ with $X \in \Var_K$ and $r \in \Z$, and morphisms are given by
  \[ \Hom_{\CHM_K(\Lambda)}((X, r), (Y, s)) = \CorrGp_K^{s-r}(X, Y; \Lambda) \]
with composition as defined in (\ref{def:composition of cycles}).
This category is called the \dfn[Chow motives]{category of Chow motives over $K$ with coefficients in $\Lambda$}.
Note that by definition, objects of $\CHM_K(\Lambda)$ are tuples $(X, r, \alpha)$ with $X \in \Var_K$, $r \in \Z$, and $\alpha \in \CorrGp^0_K(X, X; \Lambda)$ with $\alpha \circ \alpha$.

\para[rmk:morphisms in CHM]
Let $X, Y \in \Var_K$.
By definition, morphisms between the Chow motives $M(X)$ and $M(Y)$ are given by
  \[ \Hom_{\CHM_K(\Lambda)}(M(X), M(Y)) = \CorrGp^0(X, Y; \Lambda) \]
in $\CHM_K(\Lambda)$.
In particular, there is a canonical fully faithful functor
  \[ \Corr_K \hooklong \CHM_K(\Lambda) \]
which maps a variety $X$ to its Chow motive $M(X)$.

\para[def:Chow motive functor]
The canonical functor
  \[ \CHM^+_K(\Lambda) \hooklong \CHM_K(\Lambda) \]
mapping an object $(X, \alpha)$ to the pair $(X, 0, \alpha)$ is obviously fully faithful.
Composing with the functor $M^+_\Lambda : \Var_K^\circ \to \CHM^+_K(\Lambda)$ gives a canonical functor
  \[ M_\Lambda : \Var_K^\circ \too \CHM_K(\Lambda). \]
For $X \in \Var_K$, the image $M_\Lambda(X)$ is called the \dfn[]{Chow motive of $X$} (with coefficients in $\Lambda$).
By abuse of notation we write $M = M_\Lambda$ when there is no risk of confusion.

\para[rmk:CHM is symmetric monoidal]
The category $\CHM_K(\Lambda)$ is symmetric monoidal with the tensor product
  \[ (X, r, \alpha) \otimes (Y, s, \beta) = (X \times Y, r+s, \alpha \otimes \beta); \]
the unit object is the image $\IdM$ of $\IdM^+$, called the \dfn[]{identity motive}.
Every object $(X, r, \alpha) \in \CHM_K(\Lambda)$ has a dual
  \[ (X, r, \alpha)^\vee = (X, m - r, \sigma_{X,X}^*(\alpha)), \]
where $\sigma_{X,X}^*(\alpha)$ is the transpose of $\alpha$ (\ref{def:transpose of a cycle}) and $m = \dim(X)$.
It follows that this symmetric monoidal structure is \emph{rigid} (\ref{def:rigid symmetric monoidal category}).
See (\cit{andre2004introduction}, 4.1.4).

\para[def:Tate motive]
The \dfn[]{Lefschetz motive} $\LefM$ is the image of $\LefM^+$ (\ref{def:Lefschetz effective motive}) in $\CHM_K(\Lambda)$.
%The \dfn{Tate motive} $\TateM$ is defined to be its dual $\LefM^\vee$ (\ref{rmk:category of Chow motives is rigid symmetric monoidal}).
The \dfn{Tate motive $\TateM$} is defined to be the Chow motive $(\Spec(K), -1, \alpha)$ where $\alpha$ is the transpose of the class of the diagonal of $\Spec(K)$.
One can prove that there is a canonical isomorphism $\TateM \shiso \LefM^\vee$ in $\CHM_K(\Lambda)$.

\para[def:Tate twist]
For a Chow motive $M \in \CHM_K(\Lambda)$ and for an integer $i \in \Z$, we write
  \[ M(i) = M \otimes \TateM^{\otimes (-i)} \]
and call this the \dfn[Tate twist]{$i$-th Tate twist of $M$}.

For any motive $(X, r, \alpha)$, note that one has $(X, r, \alpha)(i) = (X, r-i, \alpha)$, by the definition of the tensor product (\ref{rmk:CHM is symmetric monoidal}).

\para[def:CHM is additive]
The category $\CHM_K(\Lambda)$ is $\Lambda$-linear, with the direct sum defined as follows.
Let $(X, r, \alpha)$ and $(Y, s, \beta)$ be objects of $\CHM_K(\Lambda)$.
Assuming without loss of generality that $r \le s$, there are canonical isomorphisms
  \begin{align*}
    (X, r, \alpha) &\shiso (X, s, \alpha) \otimes \TateM^{\otimes (r-s)} \\
      &\shiso (X, s, \alpha) \otimes \LefM^{\otimes (s-r)} \\
      &\shiso (X, s, \alpha) \otimes ((\P^1_K)^{s-r}, 0, \beta^{s-r}) \\
      &\shiso (X \times (\P^1_K)^{s-r}, s, (\alpha, \beta^{s-r}))
  \end{align*}
by (\ref{def:Tate motive}) and (\ref{def:Tate twist}).
Hence it is sufficient to define the direct sum in the case $r = s$.
We define
  \[ (X, s, \alpha) \oplus (Y, s, \beta) = (X \sqcup Y, s, (\alpha, \beta)). \]
See (\cit{scholl1994classical}, 1.14).

\para[rmk:morphisms in ChMot]
For varieties $X, Y \in \Var_K$, one has canonical isomorphisms
  \begin{align*}
    \Hom_{\CHM_K(\Lambda)}(M(X)(i), M(Y)(j))
      &\shiso \Hom_{\CHM_K(\Lambda)}((X, -i, \alpha), (Y, -j, \beta)) \\
      &\shiso \Corr_K^{i-j}(X, Y, \Lambda)
  \end{align*}
In particular if $X$ is purely of dimension $n$, then
  \[ \Hom_{\CHM_K(\Lambda)}(M(X)(i), M(Y)(j)) \iso \CH^{m+i-j}(X \times Y, \Lambda). \]

%%%%%%%%%%%%%%%%%%%%%%%%%%%%%%%%%%%%%%%%%%%%%%%%%%%%%%
\section{Weil cohomology theories}
%%%%%%%%%%%%%%%%%%%%%%%%%%%%%%%%%%%%%%%%%%%%%%%%%%%%%%

\para
Let $\Lambda$ be a field and $\GrVec^+_\Lambda$ the category of nonnegatively graded finite dimensional vector spaces over $\Lambda$.
The tensor product $\cdot \otimes_\Lambda \cdot$ induces a symmetric monoidal structure on $\GrVec^+_\Lambda$;
for any $V^*, W^* \in \GrVec^+_\Lambda$, one has the canonical isomorphism $V^* \otimes_\Lambda W^* \isoto W^* \otimes_\Lambda V^*$ given by
  \[ v \otimes w \mapsto (-1)^{d+e} \cdot w \otimes v \]
for $v \in V^d$, $w \in V^e$.
As usual we write $V^d$ for the $d$-th component of $V^*$ $(d \in \Z)$.

\para[def:Weil cohomology theory]
Let $\Lambda$ be a field of characteristic zero.
A \dfn[Weil cohomology theory]{Weil cohomology theory with coefficients in $\Lambda$} is the data of
  \begin{enumerate}[labelsep=10pt]
    \item  a symmetric monoidal functor $H^* : \Var_K^\circ \to \GrVec^+_\Lambda$ (we write $f^*$ for $H^*(f)$);
    \item  for every $X \in \Var_K$ purely of dimension $n$, a homomorphism of $\Lambda$-vector spaces
      \[ \tr_X : H^{2n}(X)(n) \longrightarrow \Lambda \]
    called the \dfn[]{trace homomorphism of $X$}, where for any $V^* \in \GrVec^+_\Lambda$ we write $V^*(i)$ for the tensor product $V^* \otimes_\Lambda H^2(\P^1_K)^{\otimes(-i)}$;
    \item  for every $X \in \Var_K$ and integer $0 \le i \le n = \dim(X)$, a homomorphism of abelian groups
      \[ \gamma_X^i : \CH^i(X, \Lambda) \longrightarrow H^{2i}(X) \]
    called the \dfn[]{$i$-th cycle class homomorphism of $X$};
  \end{enumerate}
satisfying the following axioms:
  \begin{enumerate}[labelsep=10pt]
    \item[WC-1]  the $K$-vector space $H^2(\P^1_K)$ has dimension 1;
    \item[WC-2]  if $X \in \Var_K$ is geometrically connected, the trace homomorphism is an isomorphism;
    \item[WC-3]  for $X, Y \in \Var_K$, the trace homomorphism on $X \times Y$ is identified with $\tr_X \otimes \tr_Y$ under the obvious isomorphisms;
    \item[WC-4]  if $X \in \Var_K$ is purely of dimension $n$, the homomorphism
      \[ H^i(X) \otimes_\Lambda H^{2n-i}(X)(n) \longrightarrow H^{2n}(X) \stackrel{\td_X}{\longrightarrow} K \]
    is a perfect pairing, where the left-hand morphism is induced by the composite
      \[ H^i(X) \otimes H^{2n-i}(X) \isoto H^{2n}(X \times X) \stackrel{\Delta_X^*}{\longrightarrow} H^{2n}(X); \]
    \item[WC-5]  the homomorphisms $\gamma^i_X$ are functorial in $X$;
    \item[WC-6]  for $\alpha \in \CH^i(X, \Lambda)$ and $\beta \in \CH^j(Y, \Lambda)$, the element
      \[ \gamma^{i+j}_{X \times Y}(\alpha \times \beta) \in H^{2i+2j}(X \times Y)(i+j) \]
    is identified with
      \[ \gamma^i_X(\alpha) \otimes \gamma^j_Y(\beta) \in H^{2i}(X)(i) \otimes_\Lambda H^{2j}(Y)(j) \]
under the canonical isomorphisms;
    \item[WC-7]  if $X \in \Var_K$ is purely of dimension $n$, the composite
      \[ \CH^n(X, \Lambda) \stackrel{\gamma^n_X}{\longrightarrow} H^{2n}(X)(n) \stackrel{\tr_X}{\longrightarrow} \Lambda \]
    is the degree homomorphism (\ref{def:degree homomorphism}).
  \end{enumerate}

By abuse of notation we write simply $H^*$ for a Weil cohomology theory $(H^*, (\tr_X)_X, (\gamma^i_X)_{i,X})$.

%Let $\GrAlg_\Lambda$ be the category of finite-dimensional graded anti-commutative $\Lambda$-algebras.
%For a Weil cohomology theory $H^*$, the diagonal morphisms $\Delta_X : X \to X \times X$ induce multiplication homomorphisms
%  \[ H^*(X) \otimes H^*(X) \isoto H^*(X \times X) \stackrel{\Delta_X^*}{\longrightarrow} H^*(X) \]
%which give $H^*(X)$ the structure of a graded anti-commutative $\Lambda$-algebra.
%Hence $H^*$ in fact defines a functor with values in $\GrAlg_\Lambda$.

\para[def:realization functors]
The functor $M : \Var_K^\circ \to \CHM_K(\Lambda)$ induces a canonical bijection between the set of Weil cohomology theories $H^*$ with coefficients in $\Lambda$, and the set of symmetric monoidal functors $\bar H^* : \CHM_K(\Lambda) \to \GrVec^+_\Lambda$ satisfying $\bar H^i(\LefM) = 0$ for $i \ne 2$.
See (\cit{andre2004introduction}, 4.2.5.1).

  \[ \begin{tikzcd}
    \Var_K^\circ \arrow{r}{M}\arrow{d}{H^*}
      & \CHM_K(\Lambda) \arrow[dashed]{dl}{\bar H^*}\\
    \GrVec^+_\Lambda
      &
  \end{tikzcd} \]

Let $H^*$ be a Weil cohomology theory.
The corresponding functor $\bar H^* : \CHM_K(\Lambda) \to \GrVec^+_\Lambda$ is called the \dfn[realization functor of a Weil cohomology theory]{realization functor of $H^*$}.

%%%%%%%%%%%%%%%%%%%%%%%%%%%%%%%%%%%%%%%%%%%%%%%%%%%%%%%
\section{Chow motives modulo Tate twists}
%%%%%%%%%%%%%%%%%%%%%%%%%%%%%%%%%%%%%%%%%%%%%%%%%%%%%%%

\para[def:Chow motives modulo Tate twists]
Let $K$ be a field.
We write $\CHM_K(\Q)/\TateM$ for the orbit category (\ref{def:orbit category}) of $\CHM_K(\Q)$ with respect to the autoequivalence $\cdot \otimes \TateM$, where $\TateM$ is the Tate motive (\ref{def:Tate motive}).
We call $\CHM_K(\Q)/\TateM$ the \dfn[Chow motives modulo Tate twists]{category of Chow motives modulo Tate twists}.
If $M, N \in \CHM_K(\Q)$ are Chow motives whose images by the projection functor $\pi : \CHM_K(\Q) \to \CHM_K(\Q)/\TateM$ are isomorphic, then we say $M$ and $N$ are \dfn[]{isomorphic up to Tate twists}.

\para[def:GrCorr]
Let $X, Y \in \Var_K$.
We define a \dfn[graded Chow correspondence]{graded (Chow) correspondence between $X$ and $Y$} as an element of the direct sum
  \[ \bigoplus_{d \in \Z} \CorrGp^{d}(X, Y; \Lambda). \]
We let $\GrCorr_K(\Lambda)$ denote the \dfn[]{category of graded (Chow) correspondences} (over $K$, with coefficients in $\Lambda$), with composition as defined in (\ref{def:composition of cycles}).

\para[def:GrChMot]
Let $\GrChMot_K(\Lambda)$ denote the karoubian envelope (\ref{def:karoubian envelope}) of the category $\GrCorr_K(\Lambda)$.
We call this the \dfn[]{category of graded Chow motives over $K$ with coefficients in $\Lambda$}.
We let $M : \Var_K \to \GrChMot_K(\Lambda)$ denote the functor associating to $X \in \Var_K$ its graded Chow motive $M(X) = (X, \delta_X)$, where $\delta_X \in \CorrGp^0(X, X; \Lambda) \hook \bigoplus_d \CorrGp^d(X, X; \Lambda)$ is the tuple $(\Delta_i(X_i))_i$ where $(X_i)_i$ is the family of irreducible components of $X$ and $\Delta_i$ is the canonical morphism $X_i \to X_i \times X$.

\para[rmk:CHM/Tate = Kar(GrCorr)]
\begin{paraprop}
Let $K$ be a field.
There is a canonical equivalence of categories
  \[ \CHM_K(\Q)/\TateM \isoto \GrChMot_K(\Q). \]
\end{paraprop}

\begin{paraproof}
By (\ref{rmk:morphisms in ChMot}), morphisms in $\CHM_K(\Q)/\TateM$ are given by
  \begin{align*}
    \Hom_{\CHM_K(\Q)/\TateM}((X, r, \alpha), (Y, s, \beta))
      &= \bigoplus_{i \in \Z} \Hom_{\CHM_K(\Q)}((X, r, \alpha), (Y, s+i, \beta)) \\
      &\shiso \beta \circ \bigoplus_{i \in \Z} \Hom_{\CHM_K(\Q)}((X, r, \alpha), (Y, s+i, \beta)) \circ \alpha \\
      &\shiso \beta \circ \bigoplus_{i \in \Z} \CorrGp^{s-r+i}(X, Y, \Q) \circ \alpha \\
      %&\shiso \beta \circ \CH^*(X \times Y, \Q) \circ \alpha \\
      &\shiso \Hom_{\GrChMot_K(\Q)}((X, \alpha), (Y, \beta))
  \end{align*}
for two Chow motives $(X, r, \alpha)$ and $(Y, s, \beta)$.
Hence we have a canonical functor defined on objects by $(X, r, \alpha) \mapsto (X, \alpha)$ and on morphisms by the canonical identifications above.
It is obvious that it is an equivalence.
\end{paraproof}

%%%%%%%%%%%%%%%%%%%%%%%%%%%%%%%%%%%%%%%%%%%%%%%%%%%%%%%%%%%%%%%%%%%%
\section{K-motives}
%%%%%%%%%%%%%%%%%%%%%%%%%%%%%%%%%%%%%%%%%%%%%%%%%%%%%%%%%%%%%%%%%%%%

\para[def:K-correspondences]
Let $K$ be a field and let $X$ and $Y$ be smooth projective varieties over $K$.
A \dfn[K-correspondence]{K-correspondence between $X$ and $Y$} (over $K$) is an element of the Grothendieck group $K_0(X \times Y)$ (\ref{def:Grothendieck groups of a scheme}).
The composition of two K-correspondences $[\sE] \in K_0(X \times Y)$ and $[\sE'] \in K_0(Y \times Z)$ is defined as
  \[ \sE' \circ \sE = \chi(\sE \otimes^\bL \sE') = \sum_{i \ge 0} (-1)^i \cdot [\sTor_i^{\O_Y}(\sE, \sE')]. \]
Note that the derived tensor product $\sE \otimes^\bL \sE'$, which lives \emph{a priori} in the category $\Pf(X \times Y \times Y \times Z)$ (\ref{def:perfect complexes}), can in fact be viewed as a perfect complex on $X \times Z$ in the obvious way; hence its Euler characteristic (\ref{def:Euler characteristic of perfect complex}) is a well-defined class in $K_0(X \times Z)$.
One verifies that this composition law is associative, and that the class of the structure sheaf of the graph of $\id_X : X \to X$ is the identity of $X$ with respect to this composition.

\para[def:KCorr]
Let $\Lambda$ be a commutative ring.
We let $\KCorr_K(\Lambda)$ denote the \dfn[]{category of K-correspondences over $K$ with coefficients in $\Lambda$}, where objects are smooth projective varieties over $K$, morphisms are given by
  \[ \Hom_{\KCorr_K(\Lambda)}(X, Y) = K_0(X \times Y) \otimes_\Z \Lambda \]
for varieties $X$ and $Y$, and composition is as defined in (\ref{def:K-correspondences}).

\para[rmk:functor SmPrVar to KCorr]
Let $\Var_K$ denote the category of smooth projective varieties over $K$.
There is a canonical functor
  \[ \Var_K^\circ \too \KCorr_K(\Lambda) \]
which is the identity on objects and which maps a morphism $f : X \to Y$ to the correspondence from $Y$ to $X$ which is given by the class of the structure sheaf of the graph of $f$ (under the canonical identification $K_0(X \times Y) \shiso K_0(Y \times X)$).

\para[def:K-motives]
The karoubian envelope (\ref{def:karoubian envelope}) of the category $\KCorr_K(\Lambda)$ is denoted
  \[ \KMot_K(\Lambda) = \Kar(\KCorr_K(\Lambda)) \]
and called the \dfn[K-motives]{category of K-motives over $K$ with coefficients in $\Lambda$}.
The functor (\ref{rmk:functor SmPrVar to KCorr}) induces a canonical functor
  \[ KM : \Var_K \too \KMot_K(\Lambda) \]
which one sees is $\Lambda$-linear and symmetric monoidal.
The image $KM(X)$ for a variety $X \in \Var_K$ is called the \dfn[]{K-motive of $X$ (with coefficients in $\Lambda$)}.

\para[rmk:KCorr(Q) = GrCorr(Q)]
Let $K$ be a field.
Recall that $\GrCorr_K$ and $\GrChMot_K$ are the categories of graded Chow correspondences (\ref{def:GrCorr}) and graded Chow motives (\ref{def:GrChMot}), respectively, over $K$.

\begin{paraprop}
There is a canonical equivalence of categories
  \[ \KCorr_K(\Q) \isoto \GrCorr_K(\Q), \]
and therefore
  \[ \KMot_K(\Q) \isoto \GrChMot_K(\Q), \]
\end{paraprop}

\begin{paraproof}
Recall that for all $X, Y \in \Var_K$, the Chern character induces isomorphisms $K_0(X \times Y) \otimes \Q \isoto \CH^*(X \times Y, \Q)$ (\ref{thm:Chern character is an isomorphism}).
By Grothendieck-Riemann-Roch (\ref{thm:GRR}) it follows that the map
  \[ [\sE] \mapsto \ch(\sE) \inpr \sqrt{\td_{X \times Y}} \]
is also an isomorphism that further maps composition of K-correspondences to composition of graded Chow correspondences.
Hence the result follows.
\end{paraproof}

%%%%%%%%%%%%%%%%%%%%%%%%%%%%%%%%%%%%%%%%%%%%%%%%%%%%%%%%%%%%%%%%%%%%%%
\section{Noncommutative Chow motives}
\label{sec:noncommutative Chow motives}
%%%%%%%%%%%%%%%%%%%%%%%%%%%%%%%%%%%%%%%%%%%%%%%%%%%%%%%%%%%%%%%%%%%%%%

\para[def:smooth proper DG categories]
Let $\AA$ be a DG category over a commutative ring $K$.
Recall that $\D(\AA)$ denotes the derived category of DG modules over $\AA$ (\ref{def:derived category of a DG category}).
Recall also that in a triangulated category, an object $X$ is called compact if the functor $\Hom(X, \cdot)$ commutes with arbitrary coproducts (\ref{def:compact object of a triangulated category}).

The DG category $\AA$ is called \dfn[smooth DG category]{smooth} if the $\AA$-$\AA$-bimodule $\Hom_\AA(\cdot, \cdot)$ is a compact object of $\D((\AA)^\circ \otimes \AA)$.
We call $\AA$ \dfn[proper DG category]{proper} if for all objects $X$ and $Y$ the complex of $K$-modules $\Hom_\AA(X, Y)$ is a compact object of $\D(\Mod(K))$.

\para[rmk:PfDG(X) is smooth and proper]
If $X$ is a smooth proper scheme over $K$, then it is possible to prove that the DG category $\PfDG(X)$ is smooth and proper.
See (\cit{toen2007moduli}, Lemma 3.27).

\para
Recall that $\HoDGMo_K$ denotes the homotopy category of $\DGCat_K$ with respect to Morita equivalence (\ref{def:Morita model structure on DGCat}).
Given a commutative ring $\Lambda$, let $\HoDGMo_{0,K}(\Lambda)$ denote the category with the same objects as $\HoDGMo_{K}$ and morphisms given by
  \[ \Hom_{\HoDGMo_{0,K}}(\AA, \BB) = K_0(\QsiHom(\AA, \BB)) \otimes_\Z \Lambda \]
for $\AA$ and $\BB$ in $\HoDGMo_{0,K}$.
Here $K_0(\QsiHom(\AA, \BB))$ denotes the Grothendieck group (\ref{def:Grothendieck group of a triangulated category}) of the triangulated category $\QsiHom(\AA, \BB)$ of quasi-functors from $\AA$ to $\BB$ (\ref{def:quasi-functor}).
This is a $\Lambda$-linear category with symmetric monoidal structure inherited from $\HoDGMo_K$.
There is a canonical symmetric monoidal functor
  \[ \HoDGMo_K \longrightarrow \HoDGMo_{0,K}(\Lambda) \]
which is the identity on objects and maps a morphism $\AA \to \BB$ to the class of the corresponding quasi-functor (\ref{thm:morphisms in HoDGQe}).

\para[def:GeoDGQe and SmPrDGQe]
Let $\SmPrDGQe_K$ (resp. $\SmPrDGMo_K$, $\SmPrDGMo_{0,K}(\Lambda)$) denote the full subcategory of $\HoDGQe_K$ (resp. $\HoDGMo_K$, $\HoDGMo_{0,K}(\Lambda)$) consisting of smooth and proper DG categories.
Similarly let $\GeoDGMo_K$ (resp. $\GeoDGMo_{0,K}(\Lambda)$) denote the full subcategory of $\SmPrDGMo_K$ (resp. $\SmPrDGMo_{0,K}(\Lambda)$) consisting of the DG categories $\PfDG(X)$ for some $X \in \Var_K$.

\para[rmk:embedding KCorr into GeoDGMo_0]
Let $X, Y \in \Var_K$.
By (\ref{eq:QsiHom of geometric DG categories}) and (\ref{rmk:isomorphism of Grothendieck groups of a scheme}) one has canonical isomorphisms
  \begin{align*}
    \Hom_{\GeoDGMo_{0,K}(\Lambda)}(\PfDG(X), \PfDG(Y))
      &= K_0(\QsiHom(\PfDG(X), \PfDG(Y))) \otimes_\Z \Lambda\\
      &\shiso K_0(\Pf(X \times Y)) \otimes_\Z \Lambda
        \tag{\ref{eq:QsiHom of geometric DG categories}} \\
      &\shiso K_0(X \times Y) \otimes_\Z \Lambda.
        \tag{\ref{rmk:isomorphism of Grothendieck groups of a scheme}}
  \end{align*}

In particular, one has a canonical fully faithful functor
  \[ \KCorr_K(\Lambda) \hooklong \GeoDGMo_{0,K}(\Lambda) \]
induced by the assignment $X \mapsto \PfDG(X)$, where $\KCorr_K(\Lambda)$ denotes the category of $K$-correspondences (\ref{def:KCorr}).

\para[def:noncommutative Chow motives]
Let $\NChMot_K(\Lambda)$, the \dfn[noncommutative Chow motives]{category of noncommutative Chow motives over $K$ with coefficients in $\Lambda$}, be the karoubian envelope (\ref{def:karoubian envelope}) of $\SmPrDGMo_{0,K}(\Lambda)$.
Explicitly, its objects are pairs $(\AA, [\varphi])$ with $\AA$ a smooth proper DG category and $[\varphi] \in K_0(\QsiHom(\AA, \AA)) \otimes \Lambda$ the class of a quasi-functor such that $[\varphi \circ \varphi] = [\varphi]$.
Its morphisms are given by
  \[ \Hom_{\NChMot_K(\Lambda)}((\AA, [\varphi]), (\BB, [\psi])) = [\psi] \circ \left(K_0(\QsiHom(\AA, \BB)) \otimes_\Z \Lambda \right) \circ [\varphi] \]
for two objects $(\AA, [\varphi])$ and $(\BB, [\psi])$.

%There are canonical fully faithful functors
%  \[ \Kar(\GeoDGMo_{0,K}(\Lambda)) \hooklongrightarrow \NChMot_K(\Lambda) \hooklongrightarrow \Kar(\HoDGMo_{0,K}(\Lambda)), \]
%so we identify $\NChMot_K(\Lambda)$ with a full subcategory of $\Kar(\HoDGMo_{0,K}(\Lambda))$.
Note that $\NChMot_K(\Lambda)$ may be identified with a full subcategory of $\Kar(\HoDGMo_{0,K}(\Lambda))$ since $\SmPrDGMo_{0,K}(\Lambda)$ is a full subcategory of $\HoDGMo_{0,K}(\Lambda)$.
It is clear that the composite
  \[ \HoDGMo_K \longrightarrow \HoDGMo_{0,K}(\Lambda) \hooklongrightarrow \Kar(\HoDGMo_{0,K}(\Lambda)) \]
restricts on the full subcategory $\SmPrDGMo_K \sub \HoDGMo_K$ to a functor
  \begin{equation} \label{eq:noncommutative Chow motive functor}
    U_\Lambda : \SmPrDGMo_K \longrightarrow \NChMot_K(\Lambda).
  \end{equation}
For a smooth proper DG category $\AA$, the image $U_\Lambda(\AA)$ is called the \dfn[]{noncommutative Chow motive of $\AA$} (with coefficients in $\Lambda$).
For a variety $X \in \Var_K$, its image by the functor
  \[ NM_\Lambda : \Var_K \stackrel{\PfDG}{\too} \SmPrDGMo_K \stackrel{U_\Lambda}{\too} \NChMot_K(\Lambda) \]
is called the \dfn[]{noncommutative Chow motive of $X$} (with coefficients in $\Lambda$).
By abuse of notation we will write $NM = NM_\Lambda$ when there is no risk of confusion.

\para[rmk:morphisms of NC motives]
Since $\GeoDGMo_{0,K}(\Lambda)$ is a full subcategory of $\SmPrDGMo_{0,K}(\Lambda)$, there is a canonical fully faithful functor $\Kar(\GeoDGMo_{0,K}(\Lambda)) \hook \NChMot_K(\Lambda)$.
Recall also that there is a fully faithful functor $\KCorr_K(\Lambda) \hook \GeoDGMo_{0,K}(\Lambda)$ (\ref{rmk:embedding KCorr into GeoDGMo_0}).
Hence we have a sequence of fully faithful functors
  \[ \KCorr_K(\Lambda) \hook \GeoDGMo_{0,K}(\Lambda) \hook \Kar(\GeoDGMo_{0,K}(\Lambda)) \hook \NChMot_K(\Lambda). \]

In particular, for varieties $X, Y \in \Var_K$, one has canonical functorial isomorphisms
  \[ \Hom_{\NChMot_K(\Lambda)}(NM(X), NM(Y)) \isoto K_0(X \times Y) \otimes_\Z \Lambda. \]

\para[def:additive invariant]
Let $\AA$ be a DG category.
Consider the DG category $\sT^\bull(\AA)$ whose objects are pairs $(i, X)$ with $i \in \{ 1, 2 \}$ and $X$ an object of $\AA$, and morphisms given by
  \[ \Hom_{\sT^\bull(\AA)}((i, X), (j, Y)) = \Hom_{\AA}(X, Y) \]
for $i \le j$ and 0 otherwise.
There are two canonical inclusion DG functors
  \begin{equation} \label{eq:additive inclusions}
    I_1, I_2 : \AA \rightrightarrows \sT^\bull(\AA).
  \end{equation}

Let $\sD$ be an additive category.
An \dfn[additive invariant]{additive invariant with values in $\sD$} is a functor $F : \DGCat_K \to \sD$ such that
  \begin{enumerate}
    \item  $F$ maps Morita equivalences (\ref{def:Morita equivalence of DG categories}) to isomorphisms in $\sD$;
    \item  for all DG categories $\AA$, the morphism
      \[ F(\AA) \oplus F(\AA) \longrightarrow F(\sT^\bull(\AA)) \]
    induced by $I_1$ and $I_2$ (\ref{eq:additive inclusions}) is an isomorphism.
  \end{enumerate}
We let $\AddInv(\sD)$ denote the category of additive invariants with values in $\sD$.

\para[rmk:universal property of NChMot]
\begin{parathm}[Tabuada]
Let $\sD$ be an additive category.
The canonical morphism
  \[ \Hom_{\mathbf{AddCat}}(\HoDGMo_{0,K}, \sD) \longrightarrow \AddInv(\sD), \]
induced by precomposing with $\DGCat \to \HoDGMo_{0,K}$, is an isomorphism.
\end{parathm}

See (\cit{tabuada2005invariants}).
\chapter{Perfect correspondences}
\label{chap:perfect correspondences}

\section{Perfect correspondences}
\label{sec:perfect correspondences}

\para[rmk:composition of perfect complexes]
Let $K$ be a field.
Recall that for a variety $X \in \Var_K$, we write $\Pf(X)$ for the triangulated category of perfect complexes on $X$ (\ref{def:perfect complexes}).
Let $X, Y, Z \in \Var_K$.
Given two perfect complexes $\sE^\bull \in \Pf(X \times Y)$ and $\sE'^\bull \in \Pf(Y \times Z)$, we define their composite $\sE'^\bull \circ \sE^\bull \in \Pf(X \times Z)$ by the formula
  \[ \sE'^\bull \circ \sE^\bull = \R (p_{XZ})_*(\L (p_{XY})^*(\sE^\bull) \otimes^\L \L (p_{YZ})^*(\sE'^\bull)) \]
where $p_{XY}$, $p_{YZ}$ and $p_{XZ}$ denote the projections from $X \times Y \times Z$ onto $X \times Y$, $Y \times Z$ and $X \times Z$, respectively.
The identity of $X$ with respect to this composition law is given by the complex with the single object
  \[ \R (\Delta_X)_*(\O_X) \]
concentrated in degree zero, where $\Delta_X : X \to X \times X$ is the diagonal morphism.

\para[def:PfCorr]
A \dfn[perfect correspondence]{perfect correspondence from $X$ to $Y$} is defined as a perfect complex $\sE^\bull \in \Pf(X \times Y)$.
We let $\PfCorr_K$ denote the category whose objects are those of $\Var_K$, morphisms are isomorphism classes of perfect complexes on the product, i.e.
  \[ \Hom_{\PfCorr_K}(X, Y) = \Iso(\Pf(X \times Y)) \]
for $X, Y \in \Var_K$, and composition is as defined in (\ref{rmk:composition of perfect complexes}).

\para[rmk:functor of derived categories represented by perfect complex]
Let $X, Y \in \Var_K$.
Any perfect correspondence $\sE^\bull$ from $X$ to $Y$ induces canonically a triangulated functor $\Phi(\sE^\bull) : \Pf(X) \to \Pf(Y)$ defined by
  \[ \Phi(\sE^\bull) = \R (p_Y)_*(\sE^\bull \otimes^\L \L (p_X)^*(\cdot)) \]
where $p_X$ and $p_Y$ are the projections from $X \times Y$ to $X$ and $Y$, respectively.
Similarly there is a triangulated functor $\Psi(\sE^\bull)$ in the other direction given by the formula
  \[ \Psi(\sE^\bull) = \R (p_X)_*(\sE^\bull \otimes^\L \L (p_Y)^*(\cdot)). \]

We will say a triangulated functor $\Pf(X) \to \Pf(Y)$ (resp. $\Pf(Y) \to \Pf(X)$) is \dfn[]{represented (resp. corepresented) by $\sE^\bull$} if it is isomorphic to the functor $\Phi(\sE^\bull)$ (resp. $\Psi(\sE^\bull)$).

\para[rmk:functor PfCorr to TriCat]
One verifies directly that the construction $\sE^\bull \mapsto \Phi(\sE^\bull)$ is functorial, i.e. there is a canonical isomorphism of functors
  \[ \Phi(\sE'^\bull) \circ \Phi(\sE^\bull) \iso \Phi(\sE'^\bull \circ \sE^\bull). \]
See (\cit{mukai1981duality}, Prop. 1.3).
Hence the assignment $\sE^\bull \mapsto \Phi(\sE^\bull)$ defines a functor
  \begin{equation} \label{eq:functor PfCorr to TriCat}
    \Pf: \PfCorr_K \longrightarrow \TriCat
  \end{equation}
to the category of triangulated categories (\ref{def:triangulated functor}).

\para
Let $\sE^\bull$ be a perfect correspondence from $X$ to $Y$ and consider the associated functor $\Phi(\sE^\bull)$.
It follows from Serre duality that the functors corepresented (\ref{rmk:functor of derived categories represented by perfect complex}) by the complexes
  \begin{equation}
    \R\sHom(\sE^\bull, \O_{X \times Y}) \otimes^\L \L (p_Y)^*(\omega_Y)[n] \label{eq:complex corepresenting left adjoint}
  \end{equation}
and
  \begin{equation}
    \R\sHom(\sE^\bull, \O_{X \times Y}) \otimes^\L \L (p_X)^*(\omega_X)[m] \label{eq:complex corepresenting right adjoint}
  \end{equation}
are left and right adjoint, respectively, to $\Phi(\sE^\bull)$, where $m = \dim(X)$ and $n = \dim(Y)$.

\para[thm:FF functors are represented on the product]
\begin{parathm}[Orlov]
Let $X$ and $Y$ be smooth projective varieties over a field $K$.
For any fully faithful triangulated functor $F : \Pf(X) \to \Pf(Y)$, there exists a perfect complex $\sE^\bull$ on $X \times Y$, unique up to isomorphism, such that $F$ is isomorphic to $\Phi(\sE^\bull)$.
\end{parathm}

For $F$ admitting a left adjoint, this was proved in (\cit{orlov1997equivalences}, Theorem 2.2).
In fact, the existence of a left adjoint is automatic by (\cit{bondal2003generators}, Theorem A.1).

Let $\FF(\Pf(X), \Pf(Y)) \sub \bHom_\TriCat(\Pf(X), \Pf(Y))$ denote the subclass of fully faithful functors.
One gets a functorial injective morphism of sets
  \begin{equation} \label{eq:morphism from triangulated functors to perfect correspondences}
    \Iso(\FF(\Pf(X), \Pf(Y))) \hooklongrightarrow \Iso(\Pf(X \times Y))
  \end{equation}
which maps (the isomorphism class of) any fully faithful triangulated functor $F : \Pf(X) \to \Pf(Y)$ to the isomorphism class of the complex $\sE^\bull$ representing it.

\para[thm:DG functors are represented on the product]
\begin{parathm}[\cit{toen2007homotopy}, Theorem 8.15]
Let $X$ and $Y$ be smooth proper schemes over $K$.
There is a canonical functorial isomorphism in $\HoDGQe_K$
  \[ \PfDG(X \times Y) \isoto \RHomDG(\PfDG(X), \PfDG(Y)). \]
\end{parathm}

Recall that $\HoDGQe_K$ denotes the homotopy category of $\DGCat_K$ with respect to quasi-equivalences (\ref{def:HoDGQe}), and that $\RHomDG$ denotes its internal hom functor (\ref{rmk:homotopy category of RHom}).

\para[rmk:bijection between perfect complexes and morphisms in HoDGQe]
Recall that the homotopy category of the DG category $\PfDG(X \times Y)$ is equivalent to the triangulated category $\Pf(X \times Y)$ (\ref{def:PfDG(X)}), and that the homotopy category of $\RHomDG(\PfDG(X), \PfDG(Y))$ is equivalent to the triangulated category of quasi-functors $\QsiHom(\PfDG(X), \PfDG(Y))$ (\ref{rmk:homotopy category of RHom}).
Hence the isomorphism (\ref{thm:DG functors are represented on the product}) induces on homotopy categories a triangulated equivalence
  \begin{equation} \label{eq:QsiHom of geometric DG categories}
    \Pf(X \times Y) \isoto \QsiHom(\PfDG(X), \PfDG(Y)).
  \end{equation}
Since isomorphism classes of quasi-functors are in bijection with morphisms in $\HoDGQe_K$ (\ref{thm:morphisms in HoDGQe}), one gets a bijection of sets
  \begin{equation} \label{eq:morphisms of geometric DG categories in HoDGQe}
    \Iso(\Pf(X \times Y)) \isoto \Hom_\HoDGQe(\PfDG(X), \PfDG(Y)).
  \end{equation}
The proof of (\ref{thm:DG functors are represented on the product}) shows that this bijection can be described explicitly as the morphism that maps a perfect complex $\sE^\bull$ to the morphism
  \[ \sF^\bull \mapsto (p_Y)_*(p_X^*(\sF^\bull) \otimes \sE^\bull) \]
in $\HoDGQe_K$.

\para[prop:derived equivalence lifts to DG quasi-equivalence]
Let $X$ and $Y$ be smooth projective varieties over $K$.
By (\ref{thm:FF functors are represented on the product}) and (\ref{rmk:bijection between perfect complexes and morphisms in HoDGQe}) one has a canonical functorial morphism
  \[ \ell_{XY} : \Iso(\FF(\Pf(X), \Pf(Y))) \hook \Iso(\Pf(X \times Y)) \isoto \Hom_{\HoDGQe}(\PfDG(X), \PfDG(Y)) \]
lifting every fully faithful triangulated functor to the DG categories.

\begin{paraprop}
Let $X$ and $Y$ be smooth projective varieties over $K$.
The triangulated categories $\Pf(X)$ and $\Pf(Y)$ are equivalent if and only if the DG categories $\PfDG(X)$ and $\PfDG(Y)$ are quasi-equivalent.
\end{paraprop}

\begin{paraproof}
If $\PfDG(X)$ and $\PfDG(Y)$ are quasi-equivalent, then by definition there is an equivalence of homotopy categories $\H^0(\PfDG(X)) \iso \H^0(\PfDG(Y))$.
Since these categories are canonically identified with $\Pf(X)$ and $\Pf(Y)$, respectively (\ref{def:PfDG(X)}), the conclusion follows.

Conversely, suppose $F : \Pf(X) \isoto \Pf(Y)$ is a triangulated equivalence and let $G : \Pf(Y) \isoto \Pf(X)$ be a quasi-inverse.
Let $\tilde F : \PfDG(X) \to \PfDG(Y)$ and $\tilde G : \PfDG(Y) \to \PfDG(X)$ be the morphisms $\ell_{XY}(F)$ and $\ell_{YX}(G)$, respectively.
Since $G \circ F$ and $F \circ G$ are isomorphic to the identities and $\ell$ is functorial, it follows that $\tilde F$ and $\tilde G$ are isomorphisms in $\HoDGQe_K$.
Therefore $\PfDG(X)$ and $\PfDG(Y)$ are isomorphic in $\HoDGQe_K$, and by (\ref{def:quasi-equivalent DG categories}) they are thus quasi-equivalent.
\end{paraproof}

\para[def:GeoDGQe]
Let $\GeoDGQe_K \sub \HoDGQe_K$ denote the full subcategory of $\HoDGQe_K$ whose objects are DG categories $\PfDG(X)$ for some variety $X \in \Var_K$.

Note that there is a canonical equivalence of categories
  \[ \PfCorr_K \isoto \GeoDGQe_K \]
which is defined on objects by $X \mapsto \PfDG(X)$ and on morphisms by the canonical isomorphisms (\ref{eq:morphisms of geometric DG categories in HoDGQe})
  \[ \Iso(\Pf(X \times Y)) \isoto \Hom_{\GeoDGQe_K}(\PfDG(X), \PfDG(Y)). \]

\para[rmk:derived equivalence is the same as isomorphism in PfCorr]
\begin{parathm}
Let $X$ and $Y$ be smooth projective varieties over a field $K$.
The triangulated categories $\Pf(X)$ and $\Pf(Y)$ are equivalent if and only if $X$ and $Y$ are isomorphic in $\PfCorr_K$.
\end{parathm}

\begin{paraproof}
By (\ref{prop:derived equivalence lifts to DG quasi-equivalence}), an equivalence of $\Pf(X)$ and $\Pf(Y)$ is the same thing as an isomorphism of $\PfDG(X)$ and $\PfDG(Y)$ in $\GeoDGQe_K$.
Since $\GeoDGQe_K$ is canonically identified with $\PfCorr_K$ (\ref{def:GeoDGQe}), the claim follows.
\end{paraproof}

%%%%%%%%%%%%%%%%%%%%%%%%%%%%%%%%%%%%%%%%%%%%%%%%%%%%%%%%%%%%%%%%%%%%%%%%%%%%%%
\section{Perfect correspondences and Chow motives}
\label{sec:Perfect correspondences and Chow motives}
%%%%%%%%%%%%%%%%%%%%%%%%%%%%%%%%%%%%%%%%%%%%%%%%%%%%%%%%%%%%%%%%%%%%%%%%%%%%%%

\para[def:cycle associated to complex]
Let $X,Y \in \Var_K$.
To each perfect correspondence $\sE^\bull$ between $X$ and $Y$ we can associate the cycle
  \[ \mu_{X,Y}(\sE^\bull) = \ch(\sE^\bull) \inpr \sqrt{\td_{X \times Y}} \in \CH^*(X \times Y, \Q). \]
% p_X^*\left(\sqrt{\td_X}\right) \inpr \ch(\sE^\bull) \inpr p_Y^*\left(\sqrt{\td_Y}\right)

\begin{paraprop}
The morphisms $\mu_{X,Y} : \Pf(X \times Y) \to \CH^*(X \times Y, \Q)$ are functorial in the sense that they respect composition of perfect correspondences.
\end{paraprop}

\begin{paraproof}
This can be verified directly using the functoriality of the Chern character and Todd classes with respect to inverse image (\ref{rmk:functoriality of Chern character wrt inverse image}), and the Grothendieck-Riemann-Roch theorem (\ref{thm:GRR}).
\end{paraproof}

\para[rmk:functor mu]
Recall that $\GrCorr_K(\Q)$ denotes the category of graded Chow correspondences with rational coefficients, where morphisms are given by $\Hom_{\GrCorr_K(\Q)}(X, Y) = \CH^*(X \times Y, \Q)$ for $X, Y \in \Var_K$ (\ref{def:GrCorr}).
The construction (\ref{def:cycle associated to complex}) defines a functor
  \[ \mu : \PfCorr_K \too \GrCorr_K(\Q). \]

\para[rmk:embedding GrCorr into CHM/Tate]
Recall from (\ref{def:Chow motives modulo Tate twists}) that $\CHM_K(\Q)/\TateM$, the category of Chow motives with rational coefficients modulo Tate twists, is by definition the orbit category (\ref{def:orbit category}) of $\CHM_K(\Q)$ with respect to the autoequivalence $\cdot \otimes \TateM$.
%Let $\pi : \CHM_K(\Q) \to \CHM_K(\Q)/\TateM$ be the canonical functor.
Since $\CHM_K(\Q)/\TateM$ is equivalent to $\GrChMot_K(\Q)$, the karoubian envelope of $\GrCorr_K(\Q)$ (\ref{rmk:CHM/Tate = Kar(GrCorr)}), one has a canonical fully faithful functor
  \[ \GrCorr_K(\Q) \hooklong \CHM_K(\Q)/\TateM \]
mapping $X \in \Var_K$ to its Chow motive $M_\Q(X)$.

\para[def:functor PfCorr to CHM/Tate]
There is a canonical functor
  \[ M_\Lambda: \PfCorr_K \stackrel{\mu}{\too} \GrCorr_K(\Q) \stackrel{\ref{rmk:embedding GrCorr into CHM/Tate}}{\hooklong} \CHM_K(\Q)/\TateM \]
where $\mu$ is as defined in (\ref{rmk:functor mu}).
By abuse of notation we will write $M = M_\Lambda$ when there is no risk of confusion.

\para[thm:derived equivalence implies isomorphism up to Tate twists]
\begin{parathm}
Let $X$ and $Y$ be smooth projective varieties over a field $K$.
If their triangulated categories $\Pf(X)$ and $\Pf(Y)$ are equivalent, then their Chow motives $M(X)$ and $M(Y)$ are isomorphic modulo Tate twist (that is, in the category $\CHM_K(\Q)/\TateM$).
\end{parathm}

\begin{paraproof}
Recall that equivalence of $\Pf(X)$ and $\Pf(Y)$ is the same as isomorphism in $\PfCorr_K$ (\ref{rmk:derived equivalence is the same as isomorphism in PfCorr}).
In particular $X$ and $Y$ are isomorphic in $\PfCorr_K$, so the functor $M : \PfCorr_K \to \CHM_K(\Q)/\TateM$ (\ref{def:functor PfCorr to CHM/Tate}) induces an isomorphism $M_\Q(X) \shiso M_\Q(Y)$ in $\CHM_K(\Q)/\TateM$.
\end{paraproof}

\para[thm:isomorphism of motives when support of dimension n]
\begin{parathm}[Orlov]
Let $X$ and $Y$ be smooth projective varieties over $K$ of dimension $n$.
Suppose there is a triangulated equivalence $F : \Pf(X) \isoto \Pf(Y)$ such that the corresponding perfect complex $\sE^\bull$ (\ref{thm:FF functors are represented on the product}) has support of dimension $n$.
Then their motives $M_\Q(X)$ and $M_\Q(Y)$ are isomorphic in $\CHM_K(\Q)$.
\end{parathm}

\begin{paraproof}
Consider the left adjoint $G : \Pf(Y) \isoto \Pf(X)$ to $F$ and let $\sE'^\bull$ be the complex corepresenting $G$ (\ref{rmk:functor of derived categories represented by perfect complex}).
Since $F$ and $G$ are mutual quasi-inverses, it follows by functoriality (\ref{eq:morphism from triangulated functors to perfect correspondences}) that $\sE^\bull$ and $\sE'^\bull$ are mutual inverses in $\PfCorr_K$.
The functor $\mu : \PfCorr_K \to \GrCorr_K(\Q)$ (\ref{def:cycle associated to complex}) induces mutually inverse correspondences $\alpha = \mu(\sE^\bull) \in \bigoplus_d \CorrGp^d(X, Y; \Q)$ and $\beta = \mu(\sE'^\bull) \in \bigoplus_d \CorrGp^d(Y,X; \Q)$ in $\GrCorr_K(\Q)$.
Since $\sE^\bull$ has support of dimension $n$, it follows that the correspondence $\alpha$ has no components in degree less than $0$ (\ref{rmk:Chern character of a coherent sheaf supported on subscheme}).
The complex $\sE'^\bull$ can also be seen to have support of dimension $n$, using the explicit description (\ref{eq:complex corepresenting left adjoint}), so $\beta$ also has no components in degree less than $0$.

Now consider the mutually inverse morphisms $f : M_\Q(X) \isoto M_\Q(Y)$ and $g : M_\Q(Y) \isoto M_\Q(X)$ in $\CHM_K(\Q)/\TateM$ induced by $\alpha$ and $\beta$, respectively, via the functor $\GrCorr_K(\Q) \hook \CHM_K(\Q)/\TateM$ (\ref{rmk:embedding GrCorr into CHM/Tate}).
By above it follows that the morphisms $f$ and $g$ also have $f^i = 0$ and $g^i = 0$ for $i < 0$, where $f^i : M_\Q(X) \to M_\Q(Y) \otimes \TateM^{\otimes i}$ and $g^i : M_\Q(Y) \to M_\Q(X) \otimes \TateM^{\otimes i}$ denote the $i$-th components of $f$ and $g$, respectively.
Therefore the proposition (\ref{prop:isomorphism in A/T gives isomorphism in A sometimes}) implies that the morphisms $f^0$ and $g^0$ are mutually inverse isomorphisms $M_\Q(X) \isoto M_\Q(Y)$ in $\CHM_K(\Q)$.
\end{paraproof}

\para
\begin{paracor}
Let $H^* : \Var^\circ \to \GrVec^+_\Q$ be a Weil cohomology theory with coefficients in $\Q$ (\ref{def:Weil cohomology theory}).
Suppose $X$ and $Y$ are smooth projective varieties over $K$ and $F : \Pf(X) \isoto \Pf(Y)$ is a triangulated equivalence.
If the corresponding perfect complex $\sE^\bull$ has support of dimension $n = \dim(X) = \dim(Y)$, then there is an isomorphism
  \[ H^*(X) \iso H^*(Y) \]
of nonnegatively graded $\Q$-vector spaces.
\end{paracor}

\begin{paraproof}
Let $\bar H^* : \CHM_K(\Q) \to \GrVec^+_\Q$ denote the realization functor of the Weil cohomology theory $H^*$ (\ref{def:realization functors}), so that there is a commutative diagram
  \[ \begin{tikzcd}
    \Var_K^\circ \arrow{r}{M_\Q}\arrow{d}{H^*}
      & \CHM_K(\Q) \arrow[dashed]{dl}{\bar H^*}\\
    \GrVec^+_\Q
      & 
  \end{tikzcd} \]
By (\ref{thm:isomorphism of motives when support of dimension n}) one has an isomorphism $M_\Q(X) \shiso M_\Q(Y)$, so the realization functor $\bar H^*$ induces an isomorphism $H^*(X) \shiso H^*(Y)$ in $\GrVec^+_\Q$.
\end{paraproof}

%%%%%%%%%%%%%%%%%%%%%%%%%%%%%%%%%%%%%%%%%%%%%%%%%%%%%%%%%%%%%%%%%%%%%%
\section{Perfect correspondences and noncommutative Chow motives}
%%%%%%%%%%%%%%%%%%%%%%%%%%%%%%%%%%%%%%%%%%%%%%%%%%%%%%%%%%%%%%%%%%%%%%

\para
Let $K$ be a field and $\Lambda$ a commutative ring.
Recall that $\SmPrDGQe_K$ (resp. $\SmPrDGMo_K$) denotes the homotopy category of smooth proper DG categories up to quasi-equivalence (resp. Morita equivalence) (\ref{def:GeoDGQe and SmPrDGQe}), and that $\GeoDGQe_K \sub \SmPrDGQe_K$ (resp. $\GeoDGMo_K \sub \SmPrDGMo_K$) denotes the full subcategory whose objects are DG categories of the form $\PfDG(X)$ for some $X \in \Var_K$.
Recall also that $\NChMot_K(\Lambda)$ denotes the category of noncommutative Chow motives over $K$ with coefficients in $\Lambda$, and that we have a canonical functor
  \[ U = U_\Lambda : \SmPrDGMo_K \too \NChMot_K(\Lambda) \]
associating to a smooth proper DG category $\AA$ its noncommutative Chow motive $U(\AA)$ (\ref{def:noncommutative Chow motives}).
Finally recall that for a variety $X \in \Var_K$ we write $NM(X)$ for its noncommutative Chow motive $U(\PfDG(X))$.

\para[rmk:functor PfCorr to NChMot]
Recall that the category of perfect correspondences $\PfCorr_K$ is canonically equivalent to the category $\GeoDGQe_K$ (\ref{def:GeoDGQe}).
Composing with the canonical functor $\GeoDGQe_K \to \GeoDGMo_K$ (\ref{rmk:quasi-equivalences are Morita equivalences}), the inclusion $\GeoDGMo_K \hook \SmPrDGMo_K$, and the functor $U : \SmPrDGMo_K \to \NChMot_K(\Lambda)$, one gets a canonical functor
  \begin{equation} \label{eq:functor PfCorr to NChMot}
    NM : \PfCorr_K \too \NChMot_K(\Lambda)
  \end{equation}
associating to a perfect correspondence $\sE^\bull \in \Pf(X \times Y)$ a morphism $NM(X) \to NM(Y)$ of noncommutative Chow motives.

\para[thm:derived equivalence implies isomorphism of noncommutative motives]
Since isomorphism in $\PfCorr_K$ is nothing but equivalence of triangulated categories of perfect complexes (\ref{rmk:derived equivalence is the same as isomorphism in PfCorr}), the existence of the functor (\ref{eq:functor PfCorr to NChMot}) demonstrates

\begin{parathm}
Let $X$ and $Y$ be smooth projective varieties over a field $K$, and $\Lambda$ a commutative ring.
If the triangulated categories $\Pf(X)$ and $\Pf(Y)$ are equivalent, then the noncommutative Chow motives $NM(X)$ and $NM(Y)$ are isomorphic in $\NChMot_K(\Lambda)$.
\end{parathm}

\para[rmk:functor CHM/Tate to NChMot]
Let $K$ be a field.
Recall that $\CHM_K(\Q)/\TateM$ denotes the category of Chow motives modulo Tate twists (\ref{def:Chow motives modulo Tate twists}), and that there is a canonical equivalence (\ref{rmk:CHM/Tate = Kar(GrCorr)})
  \[ \CHM_K(\Q)/\TateM \isoto \GrChMot_K(\Q). \]
Recall also that there is a canonical equivalence $\KMot_K(\Q) \isoto \GrChMot_K(\Q)$ between the categories of K-motives and graded Chow motives with rational coefficients (\ref{rmk:KCorr(Q) = GrCorr(Q)}).
Finally recall that there is a fully faithful functor $\KCorr_K(\Q) \hook \NChMot_K(\Q)$ (\ref{rmk:morphisms of NC motives}); since $\NChMot_K(\Q)$ is by definition karoubian, this induces a fully faithful functor
  \[ \Kar(\KCorr_K(\Q)) = \KMot(\Q) \hooklong \NChMot_K(\Q). \]

Composing the above functors one obtains a canonical fully faithful functor
  \[ R : \CHM_K(\Q)/\TateM \hooklong \NChMot_K(\Q) \]
which one verifies is in fact $\Q$-linear and symmetric monoidal.

\para[rmk:compatibility of motive functors on PfCorr]
Recall that there is a canonical functor $M : \PfCorr_K \to \CHM_K(\Q)/\TateM$ associating to a perfect correspondence a morphism of Chow motives modulo Tate twists (\ref{def:functor PfCorr to CHM/Tate}).

\begin{paraprop}
Let $K$ be a field.
The diagram
  \[ \begin{tikzcd}
      \PfCorr_K \arrow{d}{M}\arrow{dr}{NM}
      & \\
      \CHM_K(\Q)/\TateM \arrow[hookrightarrow]{r}{R}
      & \NChMot_K(\Q)
  \end{tikzcd} \]
commutes.
\end{paraprop}

\begin{paraproof}
The commutativity is obvious on objects.
Recall that
  \[ \Hom_{\PfCorr_K}(X, Y) = \Iso(\Pf(X \times Y)), \]
by definition (\ref{def:PfCorr}).
Also, there are canonical functorial isomorphisms
  \[ \Hom_{\CHM_K(\Q)/\TateM}(M(X), M(Y)) \iso \CH^*(X \times Y, \Q) \tag{by (\ref{rmk:CHM/Tate = Kar(GrCorr)})} \]
and 
  \[ \Hom_{\NChMot_K(\Q)}(NM(X), NM(Y)) \iso K_0(X \times Y)  \otimes_\Z \Q. \tag{by (\ref{rmk:morphisms of NC motives})} \]
Therefore it is sufficient to show for any two varieties $X, Y \in \Var_K$ that the diagram
  \[ \begin{tikzcd}
    \Iso(\Pf(X \times Y)) \arrow{d}{\mu}\arrow{dr}{\chi}
      & \\
    A^*(X \times Y, \Q) \arrow{r}{\sim}
      & K_0(X \times Y) \otimes_\Z \Q
  \end{tikzcd} \]
commutes, where the morphism $\mu = \mu_{X,Y}$ maps $\sE^\bull \mapsto \ch(\sE^\bull) \inpr \sqrt{\td_{X \times Y}}$ (\ref{def:cycle associated to complex}), the morphism $\chi = \chi_{X \times Y}$ is the Euler characteristic (\ref{def:Euler characteristic of perfect complex}), mapping $\sE^\bull \mapsto \sum_i (-1)^i \cdot [H^i(\sE^\bull)]$, and the horizontal morphism is the inverse of the isomorphism $\sE \mapsto \ch(\sE) \inpr \sqrt{\td_{X \times Y}}$ (\ref{rmk:KCorr(Q) = GrCorr(Q)}).
Then the commutativity follows trivially, since the Chern character of $\sE^\bull$ is by definition (\ref{def:Chern character of a complex}) the Chern character of its Euler characteristic.
\end{paraproof}

\printindex

\bibliographystyle{plainnat-aky}
\bibliography{references}

\end{document}